\documentclass[10pt,leqno,twoside]{amsart}
\usepackage{amssymb, amsmath}
\usepackage{bm}
\numberwithin{equation}{section}
\usepackage{latexsym}
\usepackage{amsmath}
\usepackage{amssymb}
\usepackage{mathrsfs}
\usepackage{graphicx}
\usepackage{ifthen}
\usepackage{color}
\usepackage[T1]{fontenc}
\usepackage[latin1]{inputenc}

\numberwithin{equation}{section}

\usepackage{latexsym}
\usepackage{amsmath}
\usepackage{amssymb}

\newtheorem{defi}{Definition}[section]
\newtheorem{thm}[defi]{Theorem}
\newtheorem{lemma}[defi]{Lemma}

\newtheorem{proposition}[defi]{Proposition}

\newcommand{\re}{\Bbb R}

\newcommand{\N}{\mathbb N}
\newcommand{\beq}{\begin{equation}}
\newcommand{\eeq}{\end{equation}}
\def \dis{\displaystyle}
\def \Om{\Omega}
\def \pt{\partial}
\def \al{\alpha}

\def \ep{\varepsilon}

\def \sg{\sigma}

\def \ds{\;\;}
\newcommand{\cL}{{\mathcal L}}

\def \dive{\mbox{\rm div }}
\def \rot{\mbox{\rm rot }}

\newcommand{\cb}{ \color{black} }

\begin{document}

\title{$L^r$-Helmholtz-Weyl decomposition for three dimensional  exterior domains}

\author{Matthias Hieber}
\address{Fachbereich Mathematik, Technische Universit\"at Darmstadt,\\
Schlossgartenstr. 7, Darmstadt 64289, Germany}
\email{hieber@mathematik.tu-darmstadt.de}

\author{Hideo Kozono}
\address{Department of Mathematics, Waseda University\\
169-8555 Tokyo,  Japan \\
Research Alliance Center of Mathematical Sciences, \\
Tohoku University  \\
980-8578 Sendai, Japan}
\email{kozono@waseda.jp}

\author{Anton Seyfert}
\address{Fachbereich Mathematik, Technische Universit\"at Darmstadt,\\
Schlossgartenstr. 7, Darmstadt 64289, Germany}
\email{seyfert@mathematik.tu-darmstadt.de}

\author{Senjo Shimizu}
\address{Graduate School of Human and Environmental Studies, 
Kyoto University \\
606-8501 Kyoto, Japan}
\email{shimizu.senjo.5s@kyoto-u.ac.jp}

\author{Taku Yanagisawa}
\address{Department of Mathematics, 
Nara Women's University \\
Nara 630-8506, Japan}
\email{taku@cc.nara-wu.ac.jp}


\maketitle

\begin{abstract} 
In this article the Helmholtz-Weyl decomposition in three dimensional exterior domains is established within the $L^r$-setting for $1<r<\infty$.  
In fact, given  an $L^r$-vector 
field $\bm{u}$, there exist $\bm{h} \in X^r_{\tiny{\mbox{har}}}(\Omega)$, 
$\bm {w}\in \dot H^{1,r}(\Omega)^3$ with $\dive \bm{w} =0$ and 
$p \in \dot H^{1,r}(\Omega)$ such that 
$\bm{u}$ may be  decomposed uniquely as  
$$
\bm {u} = \bm{h} + \rot \bm{w} + \nabla p.  
$$  
If for the given $L^r$-vector field $\bm{u}$, its harmonic part $\bm{h}$ is chosen from $V^r_{\tiny{\mbox{har}}}(\Omega)$, 
then a decomposition similar to the above one is established, too. 
However, its uniqueness holds in this case only for the case $1 < r < 3$.  
The proof given  relies on an $L^r$-variational inequality allowing  to construct 
$\bm{w} \in \dot H^{1, r}(\Omega)^3$ and $p\in \dot H^{1, r}(\Omega)$ 
for given $\bm{u} \in L^r(\Omega)^3$  
as weak solutions to certain elliptic boundary value problems. 
\end{abstract}

\section{Introduction}

The Helmholtz-Weyl decomposition plays an important role in differential geometry and partial differential equations. 
Historically speaking, such a decomposition{\cb was initiated by the work of Rham, 
Hodge and Kodaira within the setting of general $p$-forms 
on compact Riemannian manifolds and within the $L^2$-setting. 
The case of{\cb Riemannian manifolds with boundaries} was treated first by 
Morrey \cite{Mo} and Abraham, Marsden and Ratiu \cite{AMR88}.
Whereas the classical theory of the Helmholtz-Weyl decomposition is mainly 
concerned with the $L^2$-setting, it is nowadays well understood that the existence of such a decomposition 
within the $L^r$-framework for $1<r<\infty$ is of essential importance in many problems of analysis, and in particular in the analysis of the Navier-Stokes equations. 
\par
Consider in particular the famous Leray problem of finding a solution to the stationary Navier-Stokes equation with inhomogeneous Dirichlet boundary data $g$ on $\partial D$ 
in a multiply connected bounded domain $D$. The total flux condition on $g$ meaning 
$\displaystyle{\int_{\partial D} g \cdot \nu \; d\sigma =0}$ is 
a necessary condition for the solvablity of Leray's problem, and Leray \cite{Leray} himself   
showed that under the restricted flux condition on $g$, i.e.,  
$\displaystyle{\int_{S_j} g \cdot \nu \; d\sigma =0}$ for all 
$j=1,\ldots,L$, where $\partial D = \cup_{j=1}^LS_j$ with $S_j$ denoting disjoint closed smooth surfaces, there exists a solution to this problem. 
He, however, left open the general case, where only the total flux condition is being assumed. For details see e.g. the monograph by Galdi \cite{Gal11}. 
This problem was recently solved by Korobkov, Pileckas and Russo \cite{KPR15} for bounded domains in the two-dimensional setting and for axially symmetric domains in 
$\re^3$, hereby making use of arguments very much related to the Helmholtz-Weyl decomposition for {\em bounded} domains.  
The relationship  between the solvability of Leray's problem and harmonic vector fields on $D$ has been clarified in \cite{KoYa0}: the restricted flux condition on $g$ yields a solenoidal 
extension $\tilde g$ to $D$ with a trivial harmonic part, while the total flux condition yields such an extension with a non-trivial harmonic part. 
%
%It is well known that the $L^r$-decomposition 
%of vector fields is an important tool e.g. in the investigation of  the Navier-Stokes equations. 
%Note that a decomposition theorem for smooth vector fields is not sufficient for  handling nonlinear equations, due to the lack of knowledge of smooth solutions. 
\par  
Solonnikov \cite{Sol} and Fujiawara-Morimoto \cite{FuMo} were the first to prove an $L^r$-Helmholtz decomposition in $n$-dimensional {\it bounded} domains with smooth boundaries. More precisely, they proved 
that every $\bm{u} \in L^r(\Omega)^n$ can be  uniquely decomposed  as  
\begin{equation}\label{eqn:1.1}
\bm{u} = \bm{v} + \nabla p,
\end{equation}          
where $\bm{v} \in L^r_\sg(\Omega)$, i.e., $\dive \bm{v} =0$, $\bm{v}\cdot\bm{\nu}|_{\pt\Om}=0$ and $p \in H^{1, r}(\Omega)$
($\bm{\nu}$ denotes the unit outer normal to $\pt\Om$). For a survey of results known in this direction, see \cite{HiSa18}.   
For the case $n=3$, the second and the fifth author proved in \cite{KoYa1} the Helmholtz-Weyl  decomposition of $\bm{v}$ satisfying (\ref{eqn:1.1}) stating that $\bm{u}$ can de decomposed as  
\begin{equation}\label{eqn:1.2}
\bm{u} = \bm{h}  + \rot \bm{w} + \nabla p,
\end{equation}          
where $\bm{h} \in C^{\infty}(\bar \Omega)^3$ with 
$$
\dive \bm{h}=0, \; \rot \bm{h}=\bm{0} \quad\mbox{and $\bm{h}\cdot \bm{\nu}|_{\pt\Om} = 0$}, 
$$
and  $\bm{w} \in H^{1, r}(\Omega)^3$ with  
$$\dive \bm{w} =0 \mbox{ and } \bm{w}\times \bm{\nu}|_{\pt\Omega} =\bm{0}.$$ 
They also proved a similar decomposition to (\ref{eqn:1.2}) with the harmonic part $\bm{h}$ satisfying another boundary condition $\bm{h}\times \bm{\nu}|_{\pt\Omega}=\bm{0}$. 

Considering the situation of {\it unbounded} domains $\Omega$, Miyakawa \cite{Mi} 
was the first to prove a decomposition of the form  (\ref{eqn:1.1}) 
in three dimensional exterior domains with smooth boundaries. 
Later on, Simader-Sohr \cite{SiSo2} extended Miyakawa's result to the $n$-dimensional setting. 
Indeed, they showed that for
any domain $\Omega \subset \re^n$ with smooth boundary, a decomposition of the form  (\ref{eqn:1.1}) holds true if and only if 
the {\it weak Neumann}  problem for the Poisson equation is uniquely solvable in the sense that for every $\bm{u} \in L^r(\Omega)^n$ there exists a unique 
$p \in L^r_{\tiny{\mbox{loc}}}(\bar\Omega)$ with $\nabla p \in L^r(\Omega)^n$ such that 
\begin{equation}\label{eqn:1.3} 
\int_{\Omega}\nabla p\cdot\nabla \varphi dx = \int_{\Omega}\bm{u}\cdot\nabla \varphi dx 
\quad\mbox{for all $\varphi \in C^\infty_0(\bar \Omega)$}.  
\end{equation}  
Based on (\ref{eqn:1.3}), it was shown by Farwig-Sohr \cite{FaSo}, Heywood \cite{He}, Thun-Miyakawa \cite{TuMi}
and Geissert-Heck-Hieber-Sawada \cite{GeHeHiSa} that a decomposition of the form (\ref{eqn:1.1}) holds true for various types of unbounded domains. 
However, it seems that a Helmholtz-Weyl decomposition of the form (\ref{eqn:1.2}) is not known for unbounded domains. 

It is the aim of this article to  establish a Helmholtz-Weyl decomposition to the form (\ref{eqn:1.2}) for  three dimensional {\it exterior} domains. Due to the results in \cite{HiKoSeShYa}, it is 
known that $X^r_{\tiny{\mbox{har}}}(\Omega)$ as well as $V^r_{\tiny{\mbox{har}}}(\Omega)$ are of finite dimension, where 
\begin{eqnarray}
&& X^r_{\tiny{\mbox{har}}}(\Omega) \equiv \{\bm{h} \in L^r(\Omega)^3; 
\dive \bm{h} =0, \, \rot \bm{h}=\bm{0}, 
\, \bm{h}\cdot\bm{\nu}|_{\pt\Omega}=0\},  \\
&& V^r_{\tiny{\mbox{har}}}(\Omega) \equiv \{\bm{h} \in L^r(\Omega)^3; \dive \bm{h} =0, 
\, \rot \bm{h}=\bm{0}, 
\, \bm{h}\times\bm{\nu}|_{\pt\Omega}=\bm{0}\}.  
\end{eqnarray}       
In contrast to the situation of bounded domains, it is not obvious to determine whether $X^r_{\tiny{\mbox{har}}}(\Omega)$ and $V^r_{\tiny{\mbox{har}}}(\Omega)$ are finite dimensional. The reason for this 
is that Rellich's compact imdedding result $H^{1, r}(\Omega) \subset L^r(\Omega)$ does not hold for exterior domains $\Omega \subset \re^3$. Showing that the behavior of $\bm{h}$ as $|x| \to \infty$ for 
$\bm{h} \in X^r_{\tiny{\mbox{har}}}(\Omega)$ and for  $\bm{h} \in V^r_{\tiny{\mbox{har}}}(\Omega)$ can be controlled in terms of local $L^r$-bounds of $\bm{h}$ in $\Omega$ replaces  Rellich's compactness 
argument in the proof of the fact that the above spaces are finite dimensional. For details, see \cite{HiKoSeShYa}.  Another aspect of the decomposition (\ref{eqn:1.2}) with the fact that 
$\dim X^r_{\tiny{\mbox{har}}}(\Omega)< \infty$ and $\dim V^r_{\tiny{\mbox{har}}}(\Omega)<\infty$ 
is the Hodge theorem in exterior domains which states that the cohomology group on $L^r(\Omega)^3$ is 
isomorphic to $X^r_{\tiny{\mbox{har}}}(\Omega)$ and $V^r_{\tiny{\mbox{har}}}(\Omega)$  
in accordance with the boundary conditions.    

In order to obtain an  $L^r$-decomposition of the (\ref{eqn:1.2}) in exterior domains, 
we search for vector potentials $\bm{w}\in L^{r}_{\tiny{\mbox{loc}}}(\bar\Omega)^3$ 
with $\nabla \bm{w} \in L^r(\Omega)^{3^2}$ 
as weak solutions to elliptic boundary value problems of the form  
\begin{equation}\label{eqn:1.6}
\left\{
\begin{array}{lll}
\rot\rot \bm{w} & = \ds \rot \bm{u} \quad &\mbox{in }  \Omega, \\
\dive \bm{w} & = \ds 0 \quad &\mbox{in } \Omega, \\
\bm{w}\times \bm{\nu} & = \ds \bm{0} \quad&\mbox{on } \pt \Omega.   
\end{array}
\right.
\end{equation} 
To this end, we consider the $L^r$-variational inequality 
\begin{eqnarray}
&& \|\nabla \bm{w}\|_{L^r(\Omega)} \label{eqn:1.7} \\
&\le&  
C\sup\Big\{
\dfrac{\displaystyle{\Big|\int_{\Omega}(\rot \bm{w}\cdot\rot \bm{\psi} + \dive \bm{w} \dive \bm{\psi})dx\Big|}}
{\|\nabla \bm{\psi}\|_{L^{r^\prime}(\Omega)}}; \bm{\psi}\in C^\infty_0(\bar\Omega)^3, 
\bm{\psi}\times\bm{\nu}|_{\pt\Om}=\bm{0}\Big\} 
+ C\|\bm{w}\|_{L^r(D)} \nonumber 
\end{eqnarray} 
with $1/r + 1/r^{\prime} =1$, where $D$ is a {\it compact} subdomain of $\Omega$. 
An inequality of the form (\ref{eqn:1.7}) enables us to apply the generalized Lax-Milgram theorem 
in $L^r(\Omega)$, see \cite{KoYa2},  in order to construct a weak solution $\bm{w}$ to (\ref{eqn:1.6}) in the homogeneous Sobolev space $\dot H^{1, r}(\Omega)^3$.   
Note that the form $\rot \bm{u}$ on the right hand side  of (\ref{eqn:1.6}) is necessarily orthogonal 
to the null space associated with the adjoint equation to (\ref{eqn:1.6}) so that for every $\bm{u} \in L^r(\Omega)^3$ with $1 < r < \infty$ there exists a unique weak solution 
$\bm{w} \in \dot H^{1, r}(\Omega)^3$.   

On the other hand, suppose  we choose the harmonic part $\bm{h}$ from $V^r_{\tiny{\mbox{har}}}(\Omega)$ in (\ref{eqn:1.2}). Then the scalar potential $p$ has to be chosen subject  to 
homogeneous boundary condition $p|_{\pt\Om}=0$ on $\pt\Omega$. This implies that we need to solve the {\it weak Dirichlet} problem, i.e., (\ref{eqn:1.3}) 
with $\varphi \in C^\infty_0(\bar \Omega)$ replaced by $\varphi \in C^\infty_0(\Omega)$. 
It was shown by Simader-Sohr \cite{SiSo1}(see also Kozono-Sohr \cite{KoSo}) that a weak Dirichlet problem of this form is uniquely solvable if and only if $r$ satisfies $3/2 < r < 3$.  
Note that there is a crucial difference between Dirichlet and Neumann problems in exterior domains. The latter is uniquely solvable for all $1<r<\infty$ in both bounded  and exterior domains.
Observe that the functional 
$\varphi \mapsto \int_{\Omega}\bm{u}\cdot \nabla \varphi dx$ 
does not need not to vanish for harmonic function $\varphi$ in $\Omega$ with $\varphi|_{\pt\Omega}=0$.  
In order to circumvent this  difficulty, 
we modify the boundary condition for $p$ on $\pt\Om$ in accordance with $\bm{u} \in L^r(\Omega)^3$ 
in the case $1 < r \le 3/2$, which enables us to obtain a satisfactory direct decomposition even for the choice 
of $\bm{h} \in V^r_{\tiny{\mbox{har}}}(\Omega)$. 
Here we should remark that another approach using $\dot H^{1, r}(\Omega)$ was 
investigated by Simader-Sohr \cite{SiSo3}, Pr\"uss-Simonnet\cite{PuSi} 
and Shibata \cite{Sh}.    
In the case where $3 \le r < \infty$, we obtain a decomposition as in (\ref{eqn:1.2}) 
with  $\bm{h} \in V^r_{\tiny{\mbox{har}}}(\Omega)$, 
while the unique expression holds only up to a modulo one dimensional subspace.   
\par
\bigskip
This paper is organized as follows.  In the following section, we shall state our main results.  Section 3 is devoted to various estimates in the homogeneous Sobolev space $\dot H^{1, r}_0(\Omega)$.
In particular, the $L^r$-variational inequality (\ref{eqn:1.7}) estimating $\nabla \bm{w}$ 
in terms of $\dive \bm{w}$ and $\rot \bm{w}$ is investigated there in detail.  
In Section 4, based on the above  $L^r$-variational inequality, we prove the existence of a weak solution $\bm{w}$ to (\ref{eqn:1.6}), which determines then the vector potential 
in (\ref{eqn:1.2}).  In Section 5, the solvabitity of the weak Dirichlet problem (\ref{eqn:1.3}) with $\varphi \in C^\infty_0(\bar \Omega)$ replaced by $\varphi \in C^\infty_0(\Omega)$ 
is investigated. Finally, in Section 6 we shall prove our main theorems.               

\section{Results}

Throughout this paper, we impose the following assumption on the domain $\Omega$.  
\par
\bigskip
\noindent
{\bf Assumption.}
$\Omega \subset \re^3$ is an exterior domain in $\re^3$ 
with the smooth boundary $\pt\Omega$.
\par
\bigskip
\noindent
We start by introducing  some notations and various function spaces. 
For $1\le r \le \infty$, the function space  $L^r(\Omega)$ stands for all scalar- and vector-valued functions on $\Omega$ 
with the norm $\|\cdot\|_{L^r}$. We denote by $(\cdot, \cdot)$ the duality {\cb pairing} between $L^r(\Omega)$ and $L^{r'}(\Omega)$, where $1/r + 1/r' = 1$ and 
$1<r<\infty$. The spaces 
$\dot H^{1, r}(\Omega)$ and $\dot H^{1, r}_0(\Omega)$ are defined by 
\begin{equation}\label{eqn:2.1}
\dot H^{1, r}(\Omega) \equiv \{[p]; p \in L^r_{\tiny{\rm loc}}(\bar \Omega), \nabla p \in L^r(\Omega)\}, 
\quad 
\dot H^{1, r}_0(\Omega)\equiv \{p\in \dot H^{1, r}(\Omega); p|_{\pt\Om}=0\},    
\end{equation} 
where $[p]$ denotes the equivalent class consisting of all $p$ such that $q\in [p]$ implies that $p-q={\rm const.}$ in $\Omega$. Equipped with the norm $\|p\|_{\dot H^{1, r}}\equiv \|\nabla p\|_{L^r}$, 
both $\dot H^{1, r}(\Omega)$ and $\dot H^{1, r}_0(\Omega)$ are Banach spaces. 
\par
Let $\widehat H^{1, r}_0(\Omega)$ be the completion of $C^{\infty}_0(\Omega)$ 
with respect to the norm $\|\nabla p\|_{L^r}$. 
Obviously, $\widehat H^{1, r}(\Omega) \subset \dot H^{1, r}(\Omega)$ 
for all $1< r < \infty$.  
However, we shall see in Proposition \ref{pr:3.1} below that $r = 3$ is a threshold exponent pointing out  a crucial difference between $\dot H^{1, r}_0(\Omega)$ and 
$\widehat H^{1, r}_0(\Omega)$.    
    
Furthermore, we introduce the space $\dot X^r(\Omega)$, $\dot V^r(\Omega)$, $\dot X^r_{\sigma}(\Omega)$ and  $\dot V^r_{\sigma}(\Omega)$ by 
\begin{eqnarray}\label{eqn:2.2}
&& \dot X^r(\Omega)\equiv\{\bm{u}\in \dot H^{1, r}(\Omega); 
\bm{u}\cdot \bm{\nu}|_{\pt\Om}= 0\}, 
\quad
\dot V^r(\Omega)\equiv\{\bm{u}\in \dot H^{1, r}(\Omega); \bm{u}\times\bm{\nu}|_{\pt\Om}
= \bm{0}\},   \label{eqn:2.2}\\
&& 
\dot X^r_{\sg}(\Omega)\equiv \{\bm{u} \in \dot X^r(\Om); \dive \bm{u} = 0\}, 
\quad
 \dot V^r_{\sg}(\Omega)\equiv \{\bm{u} \in \dot V^r(\Om); \dive \bm{u} = 0\}, \label{eqn:2.3}
\end{eqnarray}
where $\bm{\nu}$ denotes the unit outer normal to $\pt\Om$. 
Finally, let us recall  that the $L^r$-harmonic vector fields $X^r_{\tiny{\mbox{har}}}(\Omega)$ and $V^r_{\tiny{\mbox{har}}}(\Omega)$ are defined by 
\begin{eqnarray}
&& X^r_{\tiny{\mbox{har}}}(\Omega)\equiv \{\bm{h} \in L^r(\Omega); 
\dive \bm{h}=0, \, \rot \bm{h}=0, \, \bm{h}\cdot \bm{\nu}|_{\pt\Om} =0\}, \label{eqn:2.4}\\
&& V^r_{\tiny{\mbox{har}}}(\Omega) \equiv \{\bm{h} \in L^r(\Omega); 
\dive \bm{h}=0, \, \rot \bm{h}=0, \, \bm{h}\times \bm{\nu}|_{\pt\Om} =\bm{0}\}.  
\label{eqn:2.5}
\end{eqnarray}

The following proposition was proved in \cite{HiKoSeShYa}.

\begin{proposition}\label{pr:2.1}\cite[Theorem 2.1]{HiKoSeShYa} 
Let $\Omega $ be as in the Assumption.  Then $X^r_{\tiny{\mbox{\rm har}}}(\Omega)$ and $V^r_{\tiny{\mbox{\rm har}}}(\Omega)$ are both finite dimensional for all $1 < r < \infty$.    
\end{proposition}
\par
\noindent
{\bf Remark.}
The dimensions of the spaces $\dim  X^r_{\tiny{\mbox{har}}}(\Omega)$ and $\dim  V^r_{\tiny{\mbox{har}}}(\Omega)$ are more precisely determined  in \cite[Theorem 2.2]{HiKoSeShYa} in terms of the topological structure 
of the boundary of $\pt\Om$, and in particular by the corresponding {\it Betti numbers}. 
It should be noted that 
$\dim  V^{r_1}_{\tiny{\mbox{har}}}(\Omega)< \dim  V^{r_2}_{\tiny{\mbox{har}}}(\Omega)$ 
provided $1 < r_1 \le 3/2 < r_2 < \infty$.   
See also \cite[Theorems 2.1, 2.2]{HiKoSeShYa0} in two-dimensional exterior domains.    
\par
\bigskip
The two main results of this aricle read as follows. 
\begin{thm}\label{thm:2.2}
Let $\Omega$ satisfy the Assumption, and let $1< r < \infty$.
\par
\noindent  
{\rm (i)} For every $\bm{u} \in L^r(\Omega)$ there exist $\bm{h} \in X^r_{\tiny{\mbox{har}}}(\Omega)$, 
$\bm{w} \in \dot V^r_{\sg}(\Omega)$ 
and $p \in \dot H^{1,r}(\Omega)$ such that 
\begin{equation}\label{eqn:2.6}
\bm{u} = \bm{h} + \rot \bm{w} + \nabla p
\end{equation}
and there exists a constant $C=C(\Omega,r)$ such that  
\begin{equation}\label{eqn:2.7}
\|\bm{h}\|_{L^r} + \|\nabla \bm{w}\|_{L^r} + \|\nabla p\|_{L^r} \le C\|\bm{u}\|_{L^r},
\end{equation}
{\rm (ii)} The above decomposition (\ref{eqn:2.6}) is unique in the sense that if $\bm{u}$ is decomposed as  
\begin{equation}\label{eqn:2.8}
\bm{u} = \tilde {\bm{h}} + \rot \tilde {\bm{w}} + \nabla \tilde p 
\end{equation}
for some $\tilde {\bm{h}} \in X^r_{\tiny{\mbox{har}}}(\Omega)$, 
$\tilde {\bm{w}} \in \dot V^r_{\sg}(\Omega)$ and 
$\tilde p \in \dot H^{1, r}(\Omega)$, then  
\begin{equation}\label{eqn:2.9}
\bm{h} = \tilde {\bm{h}}, \quad \rot \bm{w} = \rot \tilde {\bm{w}}, 
\quad \nabla p = \nabla \tilde p.   
\end{equation}
\end{thm}

Note that in the above theorem, the harmonic part $\bm{h}$ of $\bm{u}$ is an element of $X^r_{\tiny{\mbox{har}}}(\Omega)$ and the assertion for exterior domains parallels the one for  bounded domains;  
see \cite[Threorem 2.1]{KoYa1}. On the other hand, if the harmonic part is chosen from $V^r_{\tiny{\mbox{har}}}(\Omega)$, then the situation is quite different from the one of  bounded domains.

\begin{thm}\label{thm:2.3} 
Let $\Omega$ satisfy the Assumption. 
\par
\noindent
{\rm (i)} Let $1 < r \le 3/2$. Then, for every $\bm{u} \in L^r(\Omega)$, there exist $\bm{h} \in V^r_{\tiny{\mbox{har}}}(\Omega)$, $\bm{w}\in \dot X^r_{\sg}(\Omega)$ and $p \in \dot H^{1, r}_0(\Omega)$ 
such that 
\begin{equation}\label{eqn:2.10}
\bm{u} = \bm{h} + \rot \bm{w} + \nabla p
\end{equation}
and there exists a constant $C=C(\Omega,r)$ such that 
\begin{equation}\label{eqn:2.11}
\|\bm{h}\|_{L^r} + \|\nabla \bm{w}\|_{L^r} + \|\nabla p\|_{L^r} \le C\|\bm{u}\|_{L^r}.
\end{equation}
The above decomposition (\ref{eqn:2.10}) is unique in the sense that if $\bm{u}$ is decomposed as  
\begin{equation}\label{eqn:2.12}
\bm{u} = \tilde {\bm{h}} + \rot \tilde {\bm{w}} + \nabla \tilde p 
\end{equation}
for some $\tilde {\bm{h}} \in V^r_{\tiny{\mbox{har}}}(\Omega)$, $\tilde {\bm{w}} \in \dot X^r_{\sg}(\Omega)$ and $\tilde p \in \dot H^{1, r}_0(\Omega)$, 
then 
\begin{equation}\label{eqn:2.13}
\bm{h} = \tilde{ \bm{h}}, \quad \rot \bm{w} = \rot \tilde {\bm{w}}, \quad \nabla p = \nabla \tilde p.   
\end{equation}

\noindent
{\rm (ii)} Let $3/2 < r < 3$. Then, for every $\bm{u} \in L^r(\Omega)$ there exist $\bm{h} \in V^r_{\tiny{\mbox{har}}}(\Omega)$, $\bm{w}\in \dot X^r_{\sg}(\Omega)$ and 
$p \in \widehat H^{1, r}_0(\Omega)$  
such that $\bm{u}$ may be  decomposed as (\ref{eqn:2.10}) 
including an estimate of the form (\ref{eqn:2.11}).  
The above decomposition (\ref{eqn:2.10}) is unique in the sense that if $\bm{u}$ 
is decomposed as in  (\ref{eqn:2.12}) for some $\tilde {\bm{h}} \in V^r_{\tiny{\mbox{har}}}(\Omega)$, $\tilde {\bm{w}} \in \dot X^r_{\sg}(\Omega)$ and $\tilde p \in \widehat H^{1, r}_0(\Omega)$, 
then (\ref{eqn:2.13}) holds.  

\noindent
{\rm(iii)} Let $3 \le r < \infty$.  Then, for every $\bm{u} \in L^r(\Omega)$, there exist $\bm{h} \in V^r_{\tiny{\mbox{har}}}(\Omega)$, $\bm{w}\in \dot X^r_{\sg}(\Omega)$ and $p \in \dot H^{1, r}_0(\Omega)$  
such that $\bm{u}$ can be  decomposed as in  (\ref{eqn:2.10}) including an  estimate of the form (\ref{eqn:2.11}). 
If $\bm{u}$ is decomposed as in  (\ref{eqn:2.12}) for some $\tilde {\bm{h}} \in V^r_{\tiny{\mbox{har}}}(\Omega)$, $\tilde {\bm{w}} \in \dot X^r_{\sg}(\Omega)$ and $\tilde p \in \dot H^{1, r}_0(\Omega)$, 
then  
\begin{equation}\label{eqn:2.14}
\bm{h}-\tilde {\bm{h}} = \lambda \nabla q_0, 
\quad \rot \bm{w} = \rot \tilde {\bm{w}}, \quad p - \tilde p = \lambda q_0
\end{equation}
for some $\lambda \in \re$, where $q_0$ is the harmonic function in $\dot H^{1, r}_0(\Omega)$ such that $q_0(x) \to 1$ as $|x| \to \infty$.    
\end{thm}
%\par
\bigskip\noindent
{\bf Remarks.} (i) For every $1 < r < \infty$ 
we define $\bm{H}^r(\Omega)\equiv L^r_{\sg}(\Omega)/\{\rot 
\bm{w}; \bm{w}\in \dot V^r_\sg(\Om)\}$ which may be regarded 
as the first de Rham cohomology group on $L^r(\Omega)$.  
It follows from Theorem \ref{thm:2.2} that $\bm{H}^r(\Omega)$ is isomorphic to 
$X^r_{\tiny{\mbox{har}}}(\Omega)$ for all $1<r< \infty$. 
Hence, Theorem \ref{thm:2.2} has an aspect of the Hodge theorem (see, e.g, 
Warner \cite[Theoreom 6.11]{Wa}) in exterior domains.      
\par
%\noindent
(ii) Let $1 < r \le 3/2$ and consider Theorem \ref{thm:2.3}.  
In order to obtain a decomposition as in Theorem \ref{thm:2.3} it is necessary to choose the space $\dot H^{1,r}_0(\Omega)$ for the scalar potential $p$. Indeed, we shall show that there is a vector field 
$\bm{u} \in L^r(\Omega)$ such that for every $p \in \widehat H^{1, r}_0(\Omega)$ 
the decomposition (\ref{eqn:2.10}) is not valid.  
On the other hand, in the case $3/2 < r < 3$, we may choose the smaller 
space $\widehat H^{1, r}_0(\Omega)$ for $p$ in (\ref{eqn:2.10}). 
Notice that $\widehat H^{1, r}_0(\Omega) \subset \dot H^{1, r}_0(\Omega)$ 
for $1 < r < 3$, which is essential to the unique decomposition as in 
(\ref{eqn:2.13}).  This difference may be regarded as an effect of the harmonic vector space $V^r_{\tiny{\mbox{har}}}(\Omega)$ for the threshold value $r=3/2$. 
More precisely, it was proved  in \cite[Theorem 2.2]{HiKoSeShYa} that 
$\dim V^{r_1}_{\tiny{\mbox{har}}}(\Omega)= \dim V^{r_2}_{\tiny{\mbox{har}}}(\Omega) -1 $ provided $1< r_1 \le 3/2 < r_2 < \infty$. 
\par
%\noindent
(iii) In the case $3 \le r < \infty$, the unique expression for the harmonic part $\bm{h}$ and the scalar potential $p$ holds true, of course, modulo the harmonic function $q_0$ in $\Omega$, 
which may be regarded as a crucial difference between bounded and exterior domains.       
\par
%\noindent
(iv) Let 
$L^r_{\tau}(\Omega)\equiv \{\bm{v}\in L^r(\Om); \rot \bm{v}=\bm{0},\; \bm{v}\times \bm{\nu}|_{\pt\Om}=0\}$. 
In contrast to the above (i), we define another cohomology group $\widetilde{\bm{H}}^r(\Omega)$ by 
$$
\widetilde{\bm{H}}^r(\Omega) \equiv
\left\{
\begin{array}{ll}
& L^r_{\tau}(\Omega)/\{\nabla p; p \in \dot H^{1, r}_0(\Omega)\} 
\quad\mbox{for $1 < r \le 3/2$},  \\
& L^r_{\tau}(\Omega)/\{\nabla p; p \in \widehat H^{1, r}_0(\Omega)\} 
\quad\mbox{for $3/2 < r <3$}.  
\end{array}
\right. 
$$ 
Then, it follows from Theorem \ref{thm:2.3} (i) and (ii) 
that $\widetilde{\bm{H}}^r(\Omega)$ is 
isomorphic to $V^r_{\tiny{\mbox{har}}}(\Omega)$ for $1 < r < 3$.  
This may be regarded as a unique harmonic representative 
of the cohomology group in exterior domains.      
\par
%\noindent
(v) It seems an interesting question to ask whether $\bm{w}$ and $p$ have additional regularity properties as $\bm{w} \in H^{m+1, r}(\Omega)$ and $p \in H^{m+1, r}(\Omega)$ 
in (\ref{eqn:2.6})  and (\ref{eqn:2.10}) 
provided $\bm{u} \in H^{m, r}(\Omega)$ for $m=1, 2, \cdots$.  
An affirmative answer to this question for bounded domains $\Omega$ was given in \cite[Theorem 2.4]{KoYa1}. 
\par
\bigskip
\section{Preliminaries}
\subsection{Weak Dirichlet Problem for the Poisson Equation} \mbox{}\\
We first consider weak solutions to the Poisson equations 
in the space $\widehat H^{1, r}_0(\Omega)$.  
To this end, the following characterization of the spaces $\widehat H^{1, r}_0(\Omega)$ and $\dot H^{1, r}_0(\Omega)$ plays an important role.    

\begin{proposition}\label{pr:3.1} 
Let $\Omega$ be as in the Assumption.  
\par
%\noindent
{\rm (i)} Let $1 < r < 3$. Then  
\begin{equation}\label{eqn:3.1}
\widehat H^{1, r}_0(\Omega)=\{p \in \dot H^{1, r}_0(\Omega); p \in L^{r_\ast}(\Omega)\} 
\end{equation}
for $1/r_\ast = 1/r - 1/3$. 
\par
%\noindent
{\rm (ii)} Let  $3 \leqq r < \infty$. Then $\widehat H^{1, r}_0(\Omega) = \dot H^{1, r}_0(\Omega)$. 
Moreover, there is a closed subspace $\widetilde H^{1, r}_0(\Omega)$ of $\dot H^{1, r}_0(\Omega)$ 
such that  
\begin{equation}\label{eqn:3.2}
\dot H^{1, r}_0(\Omega) = \widetilde H^{1, r}_0(\Omega) \oplus \{\lambda q_0; \lambda \in \re\}
\quad\mbox{(direct sum)}, 
\end{equation}
where $q_0$ is the harmonic function in $\dot H^{1, r}_0(\Omega)$ such that $q_0(x)\to 1$ as $|x| \to \infty$.  
\end{proposition} 

For a proof of  proposition \ref{pr:3.1} we refer to  the articles \cite[(7.6)]{SiSo1} by Simader-Sohr  and \cite[Lemma 2.2]{KoSo} by Kozono-Sohr.     

\bigskip\noindent
{\bf Remark.} For $3/2 < r < 3$, we have similarly to (\ref{eqn:3.2}) that 
$$
\dot H^{1, r}_0(\Omega) = \widehat H^{1, r}_0(\Omega) \oplus \{\lambda q_0; \lambda \in \re\}.  
$$
\par
We next introduce for $1<r< \infty$ the generalized Laplacian 
$-\Delta_r: \widehat H^{1,r}_0(\Omega)\to \widehat H^{1, r^{\prime}}_0(\Omega)^{\ast}$ defined by 
\begin{equation}\label{eqn:3.3}
\langle -\Delta_r p, \phi\rangle \equiv (\nabla p, \nabla\phi)
\quad
\mbox{for $p\in \widehat H^{1,r}_0(\Omega)$ and $\phi \in \widehat H^{1, r^{\prime}}_0(\Omega)$}. 
\end{equation}
Here $\langle\cdot, \cdot\rangle$ denotes the duality paring between 
$\widehat H^{1, r^{\prime}}_0(\Omega)^{\ast}$ and 
$\phi \in \widehat H^{1, r^{\prime}}_0(\Omega)$.
\par
Following \cite[Theorem 7.2]{SiSo1} and \cite[Corollary 3.4]{KoSo}, the kernel $\mbox{Ker}(-\Delta_r)$ and the range $\mbox{R}(-\Delta_r)$ of the generalized 
Laplacian $-\Delta_r$ enjoy the following properties.  

\begin{proposition}\label{pr:3.2}\mbox{} 
\par
{\rm (i)} $\mbox{\rm Ker}(-\Delta_r)$ is of finite dimension and $\mbox {\rm R}(-\Delta_r)$ is a closed subspace of $\widehat H^{1,r^{\prime}}_0(\Omega)^\ast$ for all $1<r<\infty$. 
\par
%\noindent
{\rm (ii)} If $1 < r < 3$, then $\mbox{\rm Ker}(-\Delta_r) = \{0\}$.  
\par
%\noindent
{\rm (iii)} If $3/2 < r < \infty$, 
then  $\mbox{\rm R}(-\Delta_r) = \hat H^{1, r^{\prime}}_0(\Omega)^{\ast}$.
\par
%\noindent
{\rm (iv)} If $3 \leqq r < \infty$, then $\mbox{\rm Ker}(-\Delta_r) = \{\lambda q_0; \lambda \in \re\}$.  
\end{proposition}
It follows from Proposition \ref{pr:3.2} and the closed range theorem that 
in the case $1< r \leqq 3/2$,  
for given $f \in \widehat H^{1,r^{\prime}}_0(\Omega)^\ast$ there exists a unique 
$p\in \widehat H^{1, r}_0(\Omega)$ satisfying 
\begin{equation}\label{eqn:3.4}
(\nabla p, \nabla \phi) = \langle f, \phi\rangle 
\quad\mbox{for all $\phi \in \widehat H^{1, r^{\prime}}_0(\Omega)$}
\end{equation}
if and only if $f$ satisfies that $\langle f, q_0\rangle =0$. On the other hand, in the case $3/2 < r < \infty$ such a restriction is redundant so that the existence 
of $p \in \widehat H^{1, r}_0(\Omega)$ is shown for arbitrary 
$f \in \widehat H^{1, r^{\prime}}_0(\Omega)^\ast$, 
while the uniqueness property holds only for $3/2 < r < 3$.     
However, even in the case $1< r \leqq 3/2$, although $f$ does not necessarily satisfy such an orthogonal condition, 
replacing $\widehat H^{1,r^{\prime}}_0(\Omega)$ by 
$\widetilde H^{1, r^{\prime}}_0(\Omega)$ in (\ref{eqn:3.4}), 
we obtain  the following result. 
\begin{lemma}\label{lem:3.3}{\rm (Simader-Sohr \cite[Theorem 7.3]{SiSo1})}.  
Let $\Omega$ be as in the Assumption and let  $1< r \leqq 3/2$.  
Let $\widetilde H^{1, r^{\prime}}_0(\Omega)$ be as in (\ref{eqn:3.2}) with $r$ replaced by $3 \leqq r^{\prime} < \infty$.   
Then for every $g \in \widetilde H^{1, r^{\prime}}_0(\Omega)^\ast$ 
there exists a unique $q \in \widehat H^{1, r}_0(\Omega)$ such that 
\begin{equation}\label{eqn:3.5}
(\nabla q, \nabla \psi) = \langle g, \psi\rangle 
\quad \mbox{for all $\psi \in \widetilde H^{1, r^{\prime}}_0(\Omega)$},
\end{equation}
where $\langle\cdot, \cdot\rangle$ denotes the duality pairing between 
$\widetilde H^{1, r^{\prime}}_0(\Omega)^\ast$ and 
$\widetilde H^{1, r^{\prime}}_0(\Omega)$.  
Moreover, there exists a constant $C=C(\Omega,r)$ such that 
\begin{equation}\label{eqn:3.6}
\|\nabla q\|_{L^r} \le C\|g\|_{\tilde H^{1, r^{\prime}}_0(\Omega)^\ast}. 
\end{equation}
\end{lemma}

It follows from  Proposition \ref{pr:3.2} that 
the generalized Laplacian 
$-\Delta_r:\widehat H^{1,r}_0(\Omega) \to 
\widehat H^{1, r^\prime}_0(\Omega)^\ast = \dot H^{1, r^\prime}_0(\Omega)^\ast$ 
is {\it not} surjective for $1 < r \le 3/2$, 
which means that the target space $\dot H^{1, r^\prime}_0(\Omega)^\ast$ is too large to make $-\Delta_r$ surjective.  
To recover its sujectivity, 
we need to restrict the range of $-\Delta_r$ onto 
$\widetilde H^{1, r^\prime}_0(\Omega)^\ast$,  
which causes necessarily the smaller space 
$\widetilde H^{1, r^\prime}_0(\Omega)$ for the test function $\psi$ in (\ref{eqn:3.5}).   
\par
%\noindent
{\it Proof of Lemma \ref{lem:3.3}}. For the readers convenience, we give here a simplified proof of Lemma \ref{lem:3.3}.   
Let us consider the map 
$-\widetilde \Delta_r:\widehat H^{1,r}_0(\Omega) \to 
\widetilde H^{1, r^\prime}_0(\Omega)^\ast$ 
defined by 
$$
\langle -\widetilde \Delta_r p, \psi\rangle \equiv (\nabla p, \nabla\psi )
\quad\mbox{for $p \in \widehat H^{1, r}_0(\Omega)$ and 
$\psi \in \widetilde H^{1, r^\prime}_0(\Omega)$}.
$$
For the proof of the assertion it suffices to show that $-\widetilde \Delta_r$ is bijective, i.e., that $\mbox{Ker}(-\widetilde \Delta_r) = \{0\}$ 
and $\mbox{R}(-\widetilde \Delta_r) = \widetilde H^{1, r^\prime}_0(\Omega)^\ast$.  
\par
We first show that $\mbox{Ker}(-\tilde \Delta_r) = \{0\}$.  Let $p \in \mbox{Ker}(-\tilde \Delta_r)$.  Then $p \in \widehat H^{1, r}_0(\Omega)$ satisfies 
$$
(\nabla p, \nabla \psi) = 0
\quad\mbox{for all $\psi \in \widetilde H^{1, r^\prime}_0(\Omega)$}. 
$$
Since $q_0 \in \mbox{Ker}(-\Delta_{r^\prime})$, we have $(\nabla p, \nabla q_0)=0$.  Hence, by (\ref{eqn:3.2})  
$$
(\nabla p, \nabla \varphi) = 0 \quad\mbox{for all $\varphi \in \dot H^{1, r^\prime}_0(\Omega)$}. 
$$
Since $1 < r \le 3/2$ and since $C^\infty_0(\Omega)$ is dense in 
$\dot H^{1, r^\prime}_0(\Omega)$ for such $r$, 
it follows from \cite[Theorem A]{KoSo} that $p\equiv 0$ in $\Omega$, 
which yields that $\mbox{Ker}(-\widetilde \Delta_r) = \{0\}$.  
\par
We next prove that $\mbox{R}(-\widetilde \Delta_r) = \widetilde H^{1, r^\prime}_0(\Omega)^\ast$. To this end, we first show that the range $\mbox{R}(-\widetilde \Delta_r)$ 
is closed in 
$\widetilde H^{1, r^\prime}_0(\Omega)^\ast$. 
Suppose that $\{p_m\}_{m=1}^\infty \subset \widehat H^{1, r}_0(\Omega)$ satisfies 
$$
-\widetilde \Delta_rp_m \to g
\quad\mbox{ in $\widetilde H^{1, r^\prime}_0(\Omega)^\ast$ as $m \to \infty$}
$$  
for some $g\in \widetilde H^{1, r^\prime}_0(\Omega)^\ast$.  
Then there exists $f \in \dot H^{1, r^\prime}_0(\Omega)^\ast$ such that 
$$
-\Delta_rp_m \to f
\quad\mbox{ in $ \dot H^{1, r^\prime}_0(\Omega)^\ast$ as $m \to \infty$}.  
$$
Indeed, it follows from (\ref{eqn:3.2}) that for every 
$\varphi \in \widehat H^{1, r^\prime}_0(\Omega)$ 
there exist $\psi \in \widetilde H^{1, r^\prime}_0(\Omega)$ and $\lambda \in \re$ such that 
$\varphi = \psi + \lambda q_0$ with  
$\|\nabla \psi\|_{L^{r^\prime}} \le C\|\nabla \varphi\|_{L^{r^\prime}}$ for some 
$C = C(\Omega, r)$.  Since $q_0 \in \mbox{Ker}(-\Delta_{r^{\prime}})$, we have 
$(\nabla p_m, \nabla q_0) =0$ for all $m\in \N$, and hence 
\begin{eqnarray*}
\frac{|\langle-\Delta_r(p_m - p_l),\varphi\rangle }{\|\nabla \varphi\|_{L^{r^\prime}}}
&=& \frac{|(\nabla p_m -\nabla p_l, \nabla \psi)|}{\|\nabla \varphi\|_{L^{r^\prime}}} \\
&\le&  C^{-1}\frac{|(\nabla p_m -\nabla p_l, \nabla \psi)|}{\|\nabla \psi\|_{L^{r^\prime}}} \\
&\le&  C^{-1}\|-\tilde \Delta_r(p_m - p_l)\|_{\widetilde H^{1, r^\prime}_0(\Omega)^\ast}
\end{eqnarray*}
for all $m, l \in \N$. 
Since $\varphi \in \dot H^{1,r^\prime}_0(\Omega)$ is arbitrary and 
since  $-\widetilde \Delta_rp_m \to g$ 
in $ \widetilde H^{1, r^\prime}_0(\Omega)^\ast$, it follows from the above estimate that 
$-\Delta_r p_m \to f$ in $ \dot H^{1, r^\prime}_0(\Omega)^\ast$  
with some $f \in \dot H^{1, r^\prime}_0(\Omega)^\ast$. 
Since $\mbox{\rm R}(-\Delta_r)$ is closed in  
$\widehat H^{1,r^\prime}(\Omega)^\ast= \dot H^{1,r^\prime}(\Omega)^\ast$, 
implied by Propositions \ref{pr:3.1} (ii) and \ref{pr:3.2} (i), there exits 
$p \in \widehat H^{1, r}_0(\Omega)$ such that $f = -\Delta_r p$.  
It is easy to see that $\langle f, \psi\rangle = \langle g, \psi\rangle $ 
for all $\psi \in \widetilde H^{1, r^\prime}_0(\Omega)$, 
from which we obtain that $-\widetilde \Delta_r p = g$. 
This means that the range $\mbox{R}(-\widetilde \Delta_r)$ is closed in 
$\widetilde H^{1, r^\prime}_0(\Omega)^\ast$.  
\par
Now, it follows from the closed range theorem that 
$\mbox{R}(-\widetilde \Delta_r) = \mbox{Ker}((-\widetilde \Delta_r)^\ast)^\perp$,
and hence we may show that $\mbox{Ker}((-\widetilde \Delta_r)^\ast) = \{0\}$.  
Let $\psi \in \mbox{Ker}((-\widetilde \Delta_r)^\ast)$.  
Then $\psi \in \widetilde H^{1, r^\prime}_0(\Omega)$ satisfies 
$$
0=\langle p, (-\widetilde \Delta_r)^\ast\psi\rangle=\langle-\widetilde \Delta_rp, \psi\rangle 
= (\nabla p, \nabla \psi)= \langle-\Delta_{r^\prime}\psi, p\rangle
\quad\mbox{for all $p \in \widehat H^{1, r}_0(\Omega)$}.
$$
This means that $\psi\in \mbox{Ker}(-\Delta_{r^\prime})$, 
whence $\psi \in \widetilde H^{1, r^\prime}_0(\Omega) \cap
\{\lambda q_0; \lambda\in \re \}$.  
By (\ref{eqn:3.2}), $\psi =0$, which yields that 
$\mbox{Ker}((-\widetilde \Delta_r)^\ast) = \{0\}$. 
This proves Lemma \ref{lem:3.3}.  \qed

\subsection{Characterization of $X^r(D)$ and $V^r(D)$} \mbox{} 
\par
Let us take $R>0$ so large that 
$$
\pt\Om \subset B_R\equiv \{x \in \re^3; |x| < R\}, 
$$ 
and define $D \equiv \Omega \cap B_{R+2}$.  We introduce a cut-off procedure decomposing the vector field $\bm{u}$ on $\Omega$ into $\bm{u} = \eta \bm{u} + (1-\eta)\bm{u} \equiv \bm{u}_1 + \bm{u}_2$ 
for a cut-off function $\eta \in C^{\infty}_0(\re^3)$ satisfying $0 \le \eta(x) \le 1 $ for all $x \in \re^3$, $\eta(x) = 1$ for $|x| \le R+1$ and $\eta(x) = 0$ for $|x|\ge R+2$. 
We may regard $\bm{u}_2$ as the vector field on $\re^3$, while $\bm{u}_1$ should be considered as a compact perturbation of $\bm{u}$ on $D$.
Notice that $D$ is a bounded domain in $\re^3$ with boundary $\pt D = \pt\Omega \cup \{x \in \re^3; |x|=R+2\}$.   
Similarly to (\ref{eqn:2.2}), we define  $X^r(D)$ and $V^r(D)$ by 
\begin{equation}\label{eqn:3.7}  
X^r(D)=\{\bm{u} \in H^{1,r}(D); \bm{u}\cdot\bm{\nu} |_{\pt D}=0\}, \quad V^r(D)=\{\bm{u} \in H^{1,r}(D); \bm{u}\times\bm{\nu}|_{\pt D}=\bm{0}\}, 
\end{equation} 
where $\bm{\nu}$ is the unit outer normal to $\pt D$.  

The following estimates may be regarded as the Poincar\'e and the Sobolev inequalities for functions belonging to  $X^r(D)$ and $V^r(D)$. 

\begin{proposition}\label{pr:3.4} 
Let $X^r(D)$ and $V^r(D)$ be as in (\ref{eqn:3.7}).

\noindent
i)  For every $1 < r < \infty$ there exists a constant $C=C(r,R)$ such that 
\begin{equation}\label{eqn:3.8}
\|\bm{u}\|_{L^r(D)} \le C\|\nabla \bm{u}\|_{L^r(D)}
\end{equation} 
holds for all $\bm{u} \in X^r(D)$ and all $\bm{u} \in V^r(D)$.  

\noindent
ii) For every $1 < r < 3$ there exists a constant $C=C(r,R)$ such that 
\begin{equation}\label{eqn:3.9}
\|\bm{u}\|_{L^{r_*}(D)} \le C\|\nabla \bm{u}\|_{L^r(D)} \quad\mbox{with $\frac{1}{r_*} \equiv \frac 1r - \frac13 $}
\end{equation} 
holds for all $\bm{u} \in X^r(D)$ and all $\bm{u} \in V^r(D)$. 
\end{proposition}

Note that in (\ref{eqn:3.8}) and (\ref{eqn:3.9}) we do not impose on $\bm{u}$ the homogeneous boundary condition $\bm{u}|_{\pt D}=0$.  
It suffices to assume either $\bm{u}\cdot \bm{\nu}|_{\pt D} =0$ or $\bm{u}\times \bm{\nu}|_{\pt D} =\bm{0}$.  

\bigskip\noindent
{\it Proof of Proposition \ref{pr:3.4}}.  
i) Let us first prove (\ref{eqn:3.8}) for $\bm{u} \in X^r(D)$ by a contradiction argument.  Assume that  (\ref{eqn:3.8}) fails.  
Then there is a sequence $\{\bm{u}_m\}_{m=1}^\infty\subset X^r(D)$ such that $\|\bm{u}_m\|_{L^r(D)}\equiv 1$ for all $m=1, 2, \cdots$ and such that $\nabla \bm{u}_m \to \bm{0}$ in $L^r(D)$ as $m\to\infty$.  
Since $D$ is a bounded domain,  Rellich compactness theorem implies that there  is a subsequence of $\{\bm{u}_m\}_{m=1}^\infty$, denoted again by  
 $\{\bm{u}_m\}_{m=1}^\infty$ for simplicity, and a function $\bm{u} \in X^r(D)$ such that 
\begin{equation}\label{eqn:3.10}
\nabla \bm{u}_m \rightharpoonup \nabla \bm{u}
\quad\mbox{weakly in $L^r(D)$}, 
\quad
\bm{u}_m \to \bm{u} 
\quad\mbox{strongly in $L^r(D)$}
\end{equation}
as $m \to \infty$.  Since $\nabla \bm{u}_m \to \bm{0}$ in $L^r(D)$, we have $\nabla \bm{u}=\bm{0}$ in $L^r(D)$, which yields  $\bm{u}=\bm{c}$ in $D$ for some constant vector $\bm{c}\in \re^3$.     
On the other hand, there exists some $x_j\in \pt\Om$ such that $\bm{\nu}_{x_j} = \bm{e}_j$ for $j=1, 2, 3$, where $\bm{\nu}_{x_j}$ denotes the unit outer normal to $\pt\Om$ at $x_j$ 
and where $\{\bm{e}_1, \bm{e}_2, \bm{e}_3\}$ is the canonical basis in $\re^3$.  Since $\bm{u}\cdot \bm{\nu}|_{\pt D}=0$, we have that $\bm{c}\cdot \bm{e}_j =0$ for 
$j=1,2,3$, which implies $\bm{c}=\bm{0}$.  This contradicts $\|\bm{u}\|_{L^r(D)}= 1$.  

Concerning $\bm{u} \in V^r(D)$, we may argue in the same way as above since the constant vector satisfying $\bm{c}\times \bm{e}_j =\bm{0}$ for $j=1, 2, 3$ also yields that 
$\bm{u}\equiv\bm{0}$ in $\Omega$.  

\noindent
ii) The proof is quite similar to the one of i).  Assume that (\ref{eqn:3.9}) fails. Then there is a sequence $\{\bm{u}_m\}_{m=1}^\infty\subset X^r(D)$ such that 
$\|\bm{u}_m\|_{L^{r_*}(D)}\equiv 1$ for all $m=1, 2, \cdots$ and such that $\nabla \bm{u}_m \to \bm{0}$ in $L^r(D)$ as $m\to\infty$.  
Since $D$ is bounded and since $r < r_*$, $\{\bm{u}_m\}_{m=1}^\infty$ is a bounded sequence in $H^{1, r}(D)$. 
Hence, again by the Rellich compactness theorem, there is a subsequence of $\{\bm{u}_m\}_{m=1}^\infty$, which we denote again by $\{\bm{u}_m\}_{m=1}^\infty$, and a function $\bm{u} \in X^r(D)$ 
fulfilling  (\ref{eqn:3.10}). Since $\nabla \bm{u}=0$ in $\Omega$ and since $\bm{u}\cdot \bm{\nu}|_{\pt D}=0 $, we conclude similarly as above that $\bm{u} \equiv \bm{0}$.  
On the other hand, by Sobolev's inequality
$$
\|\bm{u}_m\|_{L^{r_*}(D)} \le C(\|\bm{u}_m\|_{L^r(D)} + \|\nabla \bm{u}_m\|_{L^r(D)}), 
$$ 
for some  constant $C=C(r,D)$ independent of $m \in \N$. Assertion (\ref{eqn:3.10}) thus implies 
$$
\bm{u}_m \to \bm{0}\quad\mbox{strongly in $L^{r_*}(D)$ as $m\to\infty$}.   
$$  
contradictingthe fact  that $\|\bm{u}_m\|_{L^{r_*}(D)} \equiv 1$ for all $m\in \N$. 

The proof of (\ref{eqn:3.9}) for $\bm{u} \in V^r(D)$ is parallel to the one of  $\bm{u} \in X^r(D)$ and hence omitted. 
This completes the proof of Proposition \ref{pr:3.4}. \qed

\bigskip
Oiur next aim is to  investigate variational inequalities on $X^r(D)$ and $V^r(D)$ by  $\dive$ and $\rot$ operations.  
Let us define ${\mathcal X}(D)$ and ${\mathcal V}(D)$ by 
$$
{\mathcal X}(D)\equiv \{\bm{\varphi} \in C^{\infty}(\bar D); 
\bm{\varphi}\cdot\bm{\nu}|_{\pt D}=0\},
\quad
{\mathcal V}(D)\equiv \{\bm{\psi} \in C^{\infty}(\bar D); 
\bm{\psi}\times\bm{\nu}|_{\pt D}=\bm{0}\}. 
$$
Furthermore, we introduce harmonic vector fields $X_{\tiny{\mbox{har}}}(D)$ and  $V_{\tiny{\mbox{har}}}(D)$ by 
$$
X_{\tiny{\mbox{har}}}(D)\equiv \{ \bm{h}\in {\mathcal X}(D); \dive \bm{h}=0, 
\rot \bm{h}=\bm{0}\}, 
\quad 
V_{\tiny{\mbox{har}}}(D)\equiv \{ \bm{k}\in {\mathcal V}(D); 
\dive \bm{k}=0, \rot \bm{k}=\bm{0}\}.  
$$
It is known (see \cite[Theorem 2.1]{KoYa1}) that both $X_{\tiny{\mbox{har}}}(D)$ and $V_{\tiny{\mbox{har}}}(D)$ are finite dimensional vector spaces. 
We thus may write $\dim X_{\tiny{\mbox{har}}}(D) = N$ and $\dim V_{\tiny{\mbox{har}}}(D) = L$ as  
$$
X_{\tiny{\mbox{har}}}(D)=\mbox{Span.}\{\bm{h}_1, \cdots, \bm{h}_N\}, \quad V_{\tiny{\mbox{har}}}(D)=\mbox{Span.}\{\bm{k}_1, \cdots, \bm{k}_L\}.
$$

\begin{proposition}\label{pr:3.5} 
{\rm (i)} For every $1< r < \infty$ there exists a constant $C=C(r,D)$ such that 
\begin{equation}\label{eqn:3.11}
\|\nabla \bm{u}\|_{L^r(D)} \le C\Big(
\sup_{\bm{\varphi}\in {\mathcal X}(D)}
\frac{|(\rot \bm{u}, \rot \bm{\varphi})_D + (\dive \bm{u}, \dive \bm{\varphi})_D|}
{\|\nabla \bm{\varphi}\|_{L^{r'}(D)}}+ \sum_{j=1}^N|(\bm{u}, \bm{h}_j)_D|\Big)
\end{equation}
holds for all $\bm{u} \in X^r(D)$, where $(\cdot, \cdot)_D$ denotes the inner product between $L^r(D)$ and $L^{r'}(D)$.    
\par
%\noindent
{\rm (ii)} Let $\bm{u} \in X^q(D)$ for some $1 < q < \infty$.  If 
\begin{equation}\label{eqn:3.12} 
\sup_{\bm{\varphi}\in {\mathcal X}(D)}\frac{|(\rot \bm{u}, \rot \bm{\varphi})_D + (\dive \bm{u}, \dive \bm{\varphi})_D|}{\|\nabla \bm{\varphi}\|_{L^{r'}(D)}} < \infty  
\end{equation}
for some $1 < r < \infty$, then $\bm{u} \in X^r(D)$ and satisfies (\ref{eqn:3.11}). 
\par
%\noindent
{\rm (iii)} For every $1< r < \infty$ there exists a constant $C=C(r,D)$ such that 
\begin{equation}\label{eqn:3.13}
\|\nabla \bm{u}\|_{L^r(D)} \le C\Big(
\sup_{\bm{\psi}\in {\mathcal V}(D)}
\frac{|(\rot \bm{u}, \rot \bm{\psi})_D + (\dive \bm{u}, \dive \bm{\psi})_D|}
{\|\nabla \bm{\psi}\|_{L^{r'}(D)}} 
+ \sum_{j=1}^L|(\bm{u}, \bm{k}_j)_D|\Big)
\end{equation}
for all $\bm{u} \in V^r(D)$.     
\par
%\noindent
{\rm (iv)} Let $\bm{u} \in V^q(D)$ for some $1 < q < \infty$.  If 
\begin{equation}\label{eqn:3.14} 
\sup_{\bm{\psi}\in {\mathcal V}(D)}
\frac{|(\rot \bm{u}, \rot \bm{\psi})_D + (\dive \bm{u}, \dive \bm{\psi})_D|}
{\|\nabla \bm{\psi}\|_{L^{r'}(D)}} < \infty
\end{equation}
for some $1< r < \infty$, then $\bm{u} \in V^r(D)$ and satisfies (\ref{eqn:3.13}). 
\end{proposition}

\bigskip\noindent
{\it Proof.} The estimates (\ref{eqn:3.11}) and (\ref{eqn:3.13}) are proved in  \cite[Lemma 4.1]{KoYa1}. Indeed, the authors of \cite{KoYa1} showed 
(\ref{eqn:3.11}) and (\ref{eqn:3.13}) for $\bm{u }\in X^r(D)$  and for $\bm{u} \in V^r(D)$ under the condition that  $\dive \bm{u}=0$.  It is straightforward easy to generalize their result 
to (\ref{eqn:3.11}) and (\ref{eqn:3.13}) provided $(\dive \bm{u}, \dive \bm{\varphi})_D$ and  $(\dive \bm{u}, \dive \bm{\psi})_D$ are added.   
It hence suffices to prove  assertions (ii), (iii) and (iv). 
Let us first prove (ii). 
Our proof is based on the one of \cite[Lemma 3.1]{KoYa1} by making use of the identity
\begin{equation}\label{eqn:3.15}
(\nabla \bm{u}, \nabla \bm{\varphi})_D = 
(\rot \bm{u}, \rot \bm{\varphi})_D + (\dive \bm{u}, \dive \bm{\varphi})_D - 
\int_{\pt D}\bm{u}\cdot (\bm{\varphi}\cdot\nabla \bm{\nu} 
+ \bm{\varphi}\times \rot \bm{\nu})dS
\end{equation}
for $\bm{u} \in X^q(D)$ and $\bm{\varphi}\in {\mathcal X}(D)$, and where the unit outer normal $\bm{\nu}$ can be extended continuously in the neighborhood of $\pt D$.  
Since $D$ is bounded, we may assume that $1< q < 3/2$ and $q < r$.  Let us first consider the case $1/q -1/3 \le 1/r < 1/q$.  
Define $1 < q_1 < \infty$ so that $1/q_1 \equiv 1/q - 1/3$. Then  
$
\frac{1}{q_1^\prime} = 1-\left(\frac{1}{q} - \frac{1}{3}\right) \ge  \frac{1}{r^\prime} > \frac{1}{r^\prime} - \frac 13   
$.
Hence, by Sobolev's embeddings  
$$
H^{1,r^\prime}(D) \subset L^{q_1^\prime}(D), \quad H^{1, q}(D) \subset L^{q_1}(D).   
$$   
By  (\ref{eqn:3.9}) we obtain  
\begin{equation}\label{eqn:3.16}
|(\bm{u}, \bm{\varphi})_D| 
\le \|\bm{u}\|_{L^{q_1}(D)}\|\bm{\varphi}\|_{L^{q_1^\prime}(D)}
\le C\|\nabla \bm{u}\|_{L^q(D)}\|\nabla \bm{\varphi}\|_{L^{r^\prime}(D)}.
\end{equation} 
Defining $p \in (1,\infty)$ by
$\frac 1p = \frac 32\frac{1}{q_1} = \frac{1}{q} -\frac12\left(1-\frac 1q\right)$
we have $\frac{1}{p^\prime}= \frac{1}{q_1^\prime} - \frac{1}{2}\left( 1-\frac{1}{q_1^\prime} \right)$. 
Hence, by (\ref{eqn:3.8}), Sobolev's inequality and by the trace theorem  
\begin{eqnarray*}
&& \|\bm{u}\|_{L^p(\pt D)} \le C\|\bm{u}\|_{H^{1-\frac 1q, q}(\pt D)} 
\le C\|\bm{u}\|_{H^{1, q}(D)} \le \|\nabla \bm{u}\|_{L^q(D)} \\
&& C\|\bm{\varphi}\|_{L^{p^\prime}(\pt D)} 
\le C\|\bm{\varphi}\|_{H^{1-\frac {1}{q_1^\prime}, q_1^\prime}(\pt D)} 
\le C\|\bm{\varphi}\|_{H^{1, q_1^\prime}(D)}
\le C\|\bm{\varphi} \|_{H^{1, r^\prime}(D)} 
\le C\|\nabla \bm{\varphi} \|_{L^{r^\prime}(D)}.   
\end{eqnarray*}
%Notice that $q_1^\prime \le r^\prime$, implied by $r \le q_1$.  
The above inequalities yield 
\begin{equation}\label{eqn:3.17}
\left|\int_{\pt D}\bm{u}\cdot 
(\bm{\varphi}\cdot\nabla \bm{\nu} + \bm{\varphi}\times \rot\bm{\nu})dS  \right| 
\le C\|\nabla \bm{u}\|_{L^q(D)}\|\nabla \bm{\varphi} \|_{L^{r^\prime}(D)}
\end{equation}
for all $\bm{\varphi} \in {\mathcal X}(D)$.  
By (\ref{eqn:3.12}) and (\ref{eqn:3.15})--(\ref{eqn:3.17}) we obtain  
\begin{eqnarray}
&& \sup_{\bm{\varphi}\in{\mathcal X}(D)}
\frac{|(\nabla \bm{u}, \nabla \bm{\varphi})_D + (\bm{u}, \bm{\varphi})_D|}
{\|\nabla \bm{\varphi}\|_{L^{r^\prime}(D)} + \|\bm{\varphi} \|_{L^{r^\prime}(D)} } \nonumber \\
&\le& 
\sup_{\bm{\varphi}\in {\mathcal X}(D)}
\frac{|(\rot \bm{u}, \rot \bm{\varphi})_D + (\dive \bm{u}, \dive \bm{\varphi})_D|}
{\|\nabla \bm{\varphi}\|_{L^{r'}(D)}} + C\|\nabla \bm{u}\|_{L^q(D)} \label{eqn:3.18}
<\infty. 
\end{eqnarray}
Hence,  it follows from \cite[Lemma 3.1 (2)]{KoYa1} that $\bm{u} \in X^r(D)$.

We next consider the case $1/q_1 - 1/3 \le 1/r < 1/q_1=1/q -1/3$.  Since $q_1 < r$, we obtain  
$\|\nabla \bm{\varphi}\|_{L^{r^\prime}(D)} \le C\|\nabla \bm{\varphi}\|_{L^{q_1^\prime}(D)}$ for all $\bm{\varphi} \in {\mathcal X}(D)$. 
Hence,  by (\ref{eqn:3.12})  
$$
\sup_{\bm{\varphi}\in {\mathcal X}(D)}
\frac{|(\rot \bm{u}, \rot \bm{\varphi})_D + (\dive \bm{u}, \dive \bm{\varphi})_D|}
{\|\nabla \bm{\varphi}\|_{L^{q_1^\prime}(D)}} < \infty.    
$$
Applying the above argument, yields $\bm{u} \in X^{q_1}(D)$.  Replacing $q$ by $q_1$ we see that  $\bm{u} \in X^r(D)$.

Finally, it remains to consider with the case $0 < 1/r < 1/q_1 -1/3 = 1/q - 2/3 \equiv 1/q_2$.  
By the second case, we have that $\bm{u} \in X^{q_2}(D)$.  Pick now  $q_3$ with $1< q_3 < q_2 < r$ such that 
$1/q_3 = 1/r + 1/3$.  Since $\bm{u}\in X^{q_3}(D)$, we see that  $\bm{u} \in X^r(D)$.      

We next prove (iv). T
he proof is quite similar to the one  of (ii).  
Instead of (\ref{eqn:3.15}) we make use of the identity 
\begin{equation}\label{eqn:3.19}
(\nabla \bm{u}, \nabla \bm{\psi})_D = 
(\rot \bm{u}, \rot \bm{\psi})_D + (\dive \bm{u}, \dive \bm{\psi})_D + 
\int_{\pt D}\bm{u}\cdot (\bm{\psi}\cdot\nabla \bm{\nu} - \bm{\psi}\dive\bm{\nu})dS
\end{equation}
for all $\bm{\psi} \in {\mathcal V}(D)$. Similar arguments as to derive (\ref{eqn:3.17}) yield 
$$
\left|\int_{\pt D}\bm{u}\cdot (\bm{\psi}\cdot\nabla \bm{\nu} - \bm{\psi}\dive\bm{\nu})dS  \right| 
\le C\|\nabla \bm{u}\|_{L^q(D)}\|\nabla \bm{\psi} \|_{L^{r^\prime}(D)}
$$ 
Assertions  (\ref{eqn:3.14}), (\ref{eqn:3.8}) and (\ref{eqn:3.19}) imply 
$$
\sup_{\bm{\psi}\in{\mathcal V}(D)}
\frac{|(\nabla \bm{u}, \nabla \bm{\psi})_D + (\bm{u}, \bm{\psi})_D|}
{\|\nabla \bm{\psi}\|_{L^{r^\prime}(D)} + \|\bm{\psi} \|_{L^{r^\prime}(D)} } < \infty. 
$$
Hence, it follows from \cite[Lemma 3.1 (1)]{KoYa1} that $\bm{u} \in V^r(D)$,  provied $1/q-1/3 \le 1/r < 1/q$. The case $1/r< 1/q-1/3$  can be handled 
in the same way as (ii). 
We omit the details. The proof of Proposition \ref{pr:3.5} is complete.  
\qed
\par
\bigskip  
%%%%%%%%%%%%%%%%%%%%%%%%%%%%%%%%%%%%%%%%%%%%%%%%%%%%%%%%%%%%%%%%%%%%%%%%%%%%%%%%%%%%%%%%%%%%%%%%%%%%%% 
%\subsection{$\dive$-$\rot$ estimate in $\dot H^{1, r}(\re^3)$} 
We next consider the $\dive$-$\rot$ estimates in  $\re^3$. 
To this end, let us define $\dot H^{1, r}(\re^3)$ by 
$$
\dot H^{1, r}(\re^3) \equiv \{[\bm{u}]; \bm{u} \in L^r_{\tiny{\mbox{loc}}}(\re^3); 
\nabla \bm{u} \in L^r(\re^3)\}, \quad 1< r< \infty. 
$$
Equipped with the norm $\|[\bm{u}]\|_{\dot H^{1, r}(\re^3)} \equiv \|\nabla \bm{u}\|_{L^r(\re^3)}$, the space $\dot H^{1, r}(\re^3)$ becomes  a Banach space.   
The following proposition is a variant of Proposition \ref{pr:3.5}. Denote by $(\cdot,\cdot)_{\re^3}$ the canonical duality pairing 
between $L^r(\re^3)$ and $L^{r^\prime}(\re^3)$.  
\par
%\bigskip
\begin{proposition}\label{pr:3.6}
{\rm (i)} For every $1 < r < \infty$ there exists a constant $C = C(r)$ such that 
\begin{equation}\label{eqn:3.20} 
\|\nabla \bm{u}\|_{L^r(\re^3)} \le C\sup_{\bm{\varphi} \in C^\infty_0(\re^3)}
\frac{|(\rot \bm{u}, \rot \bm{\varphi})_{\re^3} + (\dive \bm{u}, \dive \bm{\varphi})_{\re^3}|}
{\|\nabla \bm{\varphi}\|_{L^{r^\prime}(\re^3)}}
\end{equation}
for all $\bm{u} \in \dot H^{1, r}(\re^3)$. 
\par
{\rm (ii)}  Let $\bm{u} \in \dot H^{1, q}(\re^3)$ for some $1 < q < \infty$.  If $\bm{u}$ satisfies 
\begin{equation}\label{eqn:3.21}
\sup_{\bm{\varphi} \in C^\infty_0(\re^3)}
\frac{|(\rot \bm{u}, \rot \bm{\varphi})_{\re^3} + (\dive \bm{u}, \dive \bm{\varphi})_{\re^3}|}
{\|\nabla \bm{\varphi}\|_{L^{r^\prime}(\re^3)}} < \infty
\end{equation} 
for some $1 < r < \infty$, then $\bm{u} \in \dot H^{1, r}(\re^3)$ and estimate (\ref{eqn:3.20}) holds true.   
\end{proposition}

\noindent
{\it Proof}. (i) The Calder\'on-Zygmund inequality yields 
\begin{equation}\label{eqn:3.22}
\|\nabla \pt_j \bm{\psi}\|_{L^{r^\prime}(\re^3)} \le C \|-\Delta \bm{\psi}\|_{L^{r^\prime}(\re^3)}, 
\quad j =1, 2, 3
\end{equation}
for all $\bm{\psi} \in C^\infty_0(\re^3)$ and for a constant $C = C(r)$.  Since $H\equiv \{-\Delta \bm{\psi}; \bm{\psi} \in C^\infty_0(\re^3)\}$ 
is dense in $L^{r^\prime}(\re^3)$, (\ref{eqn:3.22}) and integration by parts yield  
\begin{eqnarray}
&&\sup_{\bm{\varphi} \in C^\infty_0(\re^3)}
\frac{|(\rot \bm{u}, \rot \bm{\varphi})_{\re^3} + (\dive \bm{u}, \dive \bm{\varphi})_{\re^3}|}
{\|\nabla \bm{\varphi}\|_{L^{r^\prime}(\re^3)}} \nonumber 
= \sup_{\bm{\varphi} \in C^\infty_0(\re^3)}
\frac{|(\bm{u}, -\Delta \bm{\varphi})_{\re^3}|}{\|\nabla \bm{\varphi}\|_{L^{r^\prime}(\re^3)}} \nonumber \\
&\ge &
\sup_{\bm{\psi} \in C^\infty_0(\re^3)}
\frac{|(u, -\Delta \pt_j \bm{\psi})_{\re^3}|}
{\|\nabla \pt_j \bm{\psi}\|_{L^{r^\prime}(\re^3)}} \nonumber 
\ge 
C\sup_{\bm{\psi} \in C^\infty_0(\re^3)}
\frac{|(\pt_j \bm{u}, -\Delta \bm{\psi})_{\re^3}|}
{\|-\Delta \bm{\psi}\|_{L^{r^\prime}(\re^3)}} \label{eqn:3.23}\\
&=& 
C\sup_{\bm{f}\in L^{r^\prime}(\re^3)}
\frac{|(\pt_j \bm{u}, \bm{f})_{\re^3}|}
{\|\bm{f}\|_{L^{r^\prime}(\re^3)}}  \nonumber 
= C\|\pt_j \bm{u}\|_{L^r(\re^3)}, \qquad j=1, 2, 3, \nonumber
\end{eqnarray}
with the constant $C = C(r)$. This implies (\ref{eqn:3.20}).  
\par
%\noindent
(ii) Suppose that $\bm{u} \in \dot H^{1, q}(\re^3)$ satisfies (\ref{eqn:3.21}). 
Then we obtain from (\ref{eqn:3.22}) and the above estimate that   
\begin{equation}\label{eqn:3.24}
\sup_{\bm{\psi} \in C^\infty_0(\re^3)}
\frac{|(\pt_j \bm{u}, -\Delta \bm{\psi})_{\re^3}|}
{\|-\Delta \bm{\psi}\|_{L^{r^\prime}(\re^3)}} < \infty, 
\quad
j=1, 2, 3.  
\end{equation}
Since $H=\{-\Delta\bm{\psi}; \bm{\psi} \in C^\infty_0(\re^3)\}$ is dense in $L^{r^\prime}(\re^3)$, it follows from (\ref{eqn:3.24}) that the functional 
$H \ni \bm{\psi}\mapsto (\pt_j \bm{u}, -\Delta \bm{\psi})_{\re^3} \in \re$ can be extended uniquely to a continuous functional on $L^{r^\prime}(\re^3)$. 
Hence, there exits a unique $\bm{v}_j \in L^r(\re^3)$ such that 
\begin{equation}\label{eqn:3.25}
(\bm{v}_j, -\Delta \bm{\psi})_{\re^3} = (\pt_j \bm{u}, -\Delta \bm{\psi})_{\re^3} 
\quad
\mbox{for all $\bm{\psi} \in C^\infty_0(\re^3)$, $j=1, 2, 3$}.  
\end{equation}
Since $\bm{h}_j\equiv \bm{v}_j - \pt_j \bm{u} \in L^r(\re^3)+ L^q(\re^3) \subset L^1_{\tiny{\mbox{loc}}}(\re^3)$, it follows from 
(\ref{eqn:3.25}) and Weyl's lemma that $\bm{h}_j \in C^\infty(\re^3)$ and that $\bm{h}_j$ is harmonic in $\re^3$ in the classical sense. 
Moreover, since $\bm{h}_j \in L^r(\re^3) + L^q(\re^3)$, Liouville's theorem implies 
 $\bm{h}_j \equiv \bm{0}$ in $\re^3$, which yields $\pt_j \bm{u} = \bm{v}_j \in L^r(\re^3)$ for $j=1, 2,3$.  
Hence, $\bm{u} \in \dot H^{1, r}(\re^3)$ and satisfies estimate (\ref{eqn:3.20}). 
\qed    
\par
\subsection{Characterization of $\dot X^r(\Omega)$ and $\dot V^r(\Omega)$} \mbox{}
\par
Let us recall the spaces $\dot X^r(\Omega)$ and $\dot V^r(\Omega)$ defined by (\ref{eqn:2.2}).   
Denoting ${\mathcal X}(\Omega) \equiv\{\bm{\varphi} \in C^{\infty}_0(\bar \Omega); \bm{\varphi}\cdot\bm{\nu}|_{\pt\Om}=0\}$ 
and ${\mathcal V}(\Omega)\equiv\{\bm{\psi} \in C^{\infty}_0(\bar \Omega); \bm{\psi}\times \bm{\nu}|_{\pt\Om}=\bm{0}\}$, 
similarly to the space $\widehat H^{1, r}_0(\Omega)$ in (\ref{eqn:2.1}), 
we define $\widehat X^r(\Omega)$ and $\widehat V^r(\Omega)$ to be 
the closure of ${\mathcal X}(\Om)$ and ${\mathcal V}(\Om)$ 
in $\dot H^{1, r}(\Omega)$, respectively.   
The following proposition may be regarded as the characterization of 
$\widehat X^r(\Omega)$ and $\widehat V^r(\Omega)$ corresponding to that of 
$\widehat H^{1, r}_0(\Omega)$ as in Proposition \ref{pr:3.1}.  
\par
%\bigskip
\begin{proposition}\label{pr:3.7} 
Let $\Omega$ be as in the Assumption.  
\par
{\rm (i)} If $1 < r < 3$, then  
\begin{equation}\label{eqn:3.26}
\widehat X^r(\Omega) = \{\bm{u} \in \dot X^r(\Omega); \bm{u} \in L^{r_\ast}(\Omega)\},   
\quad
\widehat V^r(\Omega) = \{\bm{u} \in \dot V^r(\Omega); \bm{u} \in L^{r_\ast}(\Omega)\}
\end{equation}
for $\frac {1}{r_\ast}= \frac 1r - \frac 13$.  
Moreover, 
\begin{equation}\label{eqn:3.27}
\dot X^r(\Omega)= \widehat X^r(\Omega)\oplus \mbox{Span.}\{\bm{h}_1, \bm{h}_2, \bm{h}_3\}, 
\end{equation}
were $\{\bm{h}_1, \bm{h}_2, \bm{h}_3\} \subset \bigcap_{1< p< \infty} \dot X^p(\Omega)$ 
satisfies that 
$\dive \bm{h}_j =0$, $\rot \bm{h}_j =\bm{0}$, $\bm{h}_j(x) \to \bm{e}_j$ as 
$|x|\to \infty$ for $j=1, 2, 3$, 
and where  $\{\bm{e}_1, \bm{e}_2, \bm{e}_3\}$ denotes the canonical basis of $\re^3$. 
Similarly,
\begin{equation}\label{eqn:3.28}
\dot V^r(\Omega)= \widehat V^r(\Omega)\oplus \mbox{Span.}\{\bm{k}_1, \bm{k}_2, \bm{k}_3\}, 
\end{equation}
where 
$\{\bm{k}_1, \bm{k}_2, \bm{k}_3\} \subset \bigcap_{1< p< \infty} \dot V^p(\Omega)$ 
satisfies that 
$\dive \bm{k}_j =0$, $\rot \bm{k}_j =\bm{0}$, $\bm{k}_j(x) \to \bm{e}_j$ 
as $|x| \to \infty$ for $j =1, 2, 3$.   
\par
%\noindent
{\rm (ii)} If $3 \le r < \infty$, then  
\begin{equation}\label{eqn:3.29}
\dot X^r(\Omega) = \widehat X^r(\Omega), 
\quad
\dot V^r(\Omega) = \widehat V^r(\Omega).  
\end{equation}

\medskip
\end{proposition} 

\noindent
{\it Proof.} 
Since the proof of (\ref{eqn:3.26}) and (\ref{eqn:3.29}) is parallel to that of Proposition \ref{pr:3.1}, we only prove (\ref{eqn:3.27}) and (\ref{eqn:3.28}).  
Let us first show (\ref{eqn:3.27}).  Take arbitrary $\{\bm{w}_1, \bm{w}_2, \bm{w}_3\} \subset \bigcap_{1<p< \infty}\dot X^p(\Omega)$ such that 
$\bm{w}_j(x) \to \bm{e}_j$ as $|x| \to \infty$ and such that $\bm{w}_j(\cdot) - \bm{e}_j \in L^{r_\ast}(\Omega)$ for $j=1, 2, 3$. Note that such a triple $\{\bm{w}_1, \bm{w}_2, \bm{w}_3\}$ 
can be constructed for instance as $\bm{w}_j(x)= (1-\eta(x))\bm{e}_j$ for $j=1, 2, 3$ with the cut-off function $\eta \in C^{\infty}(\re^3)$ satisfying $\eta(x)=1$ for $|x| \le R+1$ and 
$\eta(x)= 0$ for $|x|\ge R+2$. Recall that $R$ is chosen so large that $\pt \Omega \subset B_R = \{x \in \re^3; |x| < R\}$. We shall first show that 
\begin{equation}\label{eqn:3.30}
\dot X^r(\Om) = \widehat X^r(\Omega)\oplus \mbox{Span.}\{\bm{w}_1, \bm{w}_2, \bm{w}_3\}.
\end{equation}
Indeed, let $\bm{w} \in \dot X^r(\Omega)$.  
Since $1< r <3$, by \cite[Corollary 2.2]{GiSo} that exists a constant vector $\bm{c} \in \re^3$ 
 such that $\bm{w}(\cdot)- \bm{c} \in L^{r_\ast}(\Omega)$. 
Since $\bm{c} = \sum_{j=1}^3\lambda_j\bm{e}_j$ with $\lambda_j = \bm{c}\cdot \bm{e}_j$ for $j=1,2,3$, we may set $\bm{\hat {w}}(x):= \bm{w}(x) - \sum_{j=1}^3\lambda_j\bm{w}_j(x)$.  
Since $\bm{\hat {w}}(x) = \bm{w}(x) - c - \sum_{j=1}^3\lambda_j(\bm{w}_j(x) - \bm{e}_j)$, we see that $\bm{\hat {w}}\in L^{r_\ast}(\Omega)$. 
Hence, it follows from (\ref{eqn:3.26}) that $\bm{\hat {w}} \in \widehat X^r(\Omega)$.  
Since $\bm{w} \in \dot X^r(\Omega)$ is arbitrary and 
since $\widehat X^r(\Omega) \cap \mbox{Span.}\{\bm{w}_1, \bm{w}_2, \bm{w}_3\} = \{\bm{0}\}$, 
we obtain (\ref{eqn:3.30}). 
\par
Next, we show that the triple $\{\bm{w}_1, \bm{w}_2, \bm{w}_3\}$ may be replaced by the triplet $\{\bm{h}_1, \bm{h}_2, \bm{h}_3\}$ of harmonic vector fields as in (\ref{eqn:3.27}).  
To this end, consider the following Neumann problem  
\begin{equation}\label{eqn:3.31} 
\left\{
\begin{array}{ll}
& \Delta q_j =0 \quad \mbox{in $\Omega$}, \\
& \dis{\frac{\pt q_j}{\pt \bm{\nu}} = \bm{e}_j\cdot\bm{\nu} } \quad\mbox{on $\pt\Omega$}, \\
& \nabla q_j(x) \to 0 \quad\mbox{as $|x| \to \infty$} \quad\mbox{for $j=1,2,3$}.  
\end{array}
\right.
\end{equation}
Since $\int_{\pt\Om}\bm{e}_j\cdot\bm{\nu} dS =0 $ for $j=1,2,3$, there is a smooth solution $q_j$ of (\ref{eqn:3.31}) satisfying $\nabla q_j, \nabla^2 q_j \in L^p(\Omega)$ for all 
$1< p< \infty$, see, e.g., \cite[Theorem 4.1]{SiSo3}.   
Defining $\bm{h}_j(x) = \bm{e}_j - \nabla q_j(x)$, 
we see that $\bm{h}_j \in \bigcap_{1<p< \infty}\dot X^p(\Omega)$ with 
$\bm{h}_j(\cdot) - \bm{e}_j \in \bigcap_{1<p< \infty}L^p(\Om)$ satisfying  
$\dive \bm{h}_j=0$, $\rot \bm{h}_j =\bm{0}$ and $\bm{h}_j(x) \to \bm{e}_j$ 
as $|x| \to \infty$ for $j=1,2,3$.   
This means that the harmonic triple $\{\bm{h}_1, \bm{h}_2, \bm{h}_3\}$ has the same properties 
as the one  of $\{\bm{w}_1, \bm{w}_2, \bm{w}_3\}$ in (\ref{eqn:3.30}), and hence the proof of 
(\ref{eqn:3.27}) is complete. 
\par
We next show (\ref{eqn:3.28}). Considering the triplet $\{\bm{v}_1, \bm{v}_2, \bm{v}_3\} \subset \bigcap_{1<p<\infty}\dot V^r(\Omega)$ with $\bm{v}_j(\cdot)-\bm{e}_j \in L^{r_\ast}(\Omega)$ such that 
$\bm{v}_j(x) \to \bm{e}_j$ as $|x|\to\infty$ for $j=1,2,3$, similarly to (\ref{eqn:3.30}), we obtain the decomposition
\begin{equation}\label{eqn:3.32} 
\dot V^r(\Omega) = \widehat V^r(\Omega)\oplus\mbox{Span.}\{\bm{v}_1, \bm{v}_2, \bm{v}_3\}.  
\end{equation} 
Hence, for the proof of (\ref{eqn:3.28}) we may construct the harmonic triplet 
$\{\bm{k}_1, \bm{k}_2, \bm{k}_3\} \subset \bigcap_{1<p<\infty} \dot V^p(\Omega)$  
with the same properties as those of $\{\bm{v}_1, \bm{v}_2, \bm{v}_3\}$.  
For that purpose, similarly to (\ref{eqn:3.31}), we consider the following Dirichlet problem
\begin{equation}\label{eqn:3.33} 
\left\{
\begin{array}{ll}
& \Delta \pi_j =0 \quad \mbox{in $\Omega$}, \\
& \pi_j = x\cdot \bm{e}_j \quad\mbox{on $\pt\Omega$}, \\
& \pi_j(x) \to 0 \quad\mbox{as $|x| \to \infty$} \quad\mbox{for $j=1,2,3$}.  
\end{array}
\right.
\end{equation}
Let $\pi$ the solution of of (\ref{eqn:3.33}) satisfying  $\nabla \pi \in L^s(\Omega)$ for all $3/2 < s \le \infty$ and $\nabla^2 \pi \in L^p(\Om)$ for all $1 < p < \infty$; 
see, e.g., \cite[Theorem 4.2]{SiSo2}.  
Defining $\bm{k}_j(x) = \nabla (x\cdot \bm{e}_j - \pi_j(x))$, we see that $\{\bm{k}_1, \bm{k}_2, \bm{k}_3\} \subset \bigcap_{1< p < \infty}\dot V^p(\Om)$ 
with $\bm{k}_j(\cdot) - \bm{e}_j = \nabla \pi_j(\cdot) \in \bigcap_{3/2< s< \infty}L^s(\Omega)$ 
are satisfying $\dive \bm{k}_j =0$, $\rot \bm{k}_j =0$ and $\bm{k}_j(x) \to \bm{e}_j$ as $|x|\to\infty$ for $j=1,2,3$. 
Since $3/2 < r_\ast$, the harmonic triplet $\{\bm{k}_1, \bm{k}_2, \bm{k}_3\}$ has the required properties of $\{\bm{v}_1, \bm{v}_2, \bm{v}_3\}$ in (\ref{eqn:3.32}), 
and hence we obtain the desired decomposition (\ref{eqn:3.28}). 
This proves Proposition \ref{pr:3.7}.  \qed
\par
\bigskip
%\noindent 
We next consider fundamental inequalities on $\dot X^r(\Omega)$ and $\dot V^r(\Omega)$. 
\par
%\bigskip 
\begin{proposition}\label{pr:3.8} 
Let $\Omega$ be as in the Assumption and $D = \Omega\cap B_{R+2}$ for $R>0$ sufficiently large.  If $1 < r < \infty$, then there exists $C = C(\Omega, r, R)$ such that 
\begin{equation}\label{eqn:3.34}
\|\bm{w}\|_{L^r(D)} \le C\|\nabla \bm{w}\|_{L^r}
\end{equation}
for all $\bm{w} \in \dot X^r(\Omega)$ and all $\bm{w} \in \dot V^r(\Omega)$.    
\end{proposition}

\noindent
{\it Proof.} We prove the assertion for $\bm{w} \in \dot X^r(\Omega)$, only. The proof for $\bm{w} \in \dot V^r(\Omega)$ is similar.  
We make use of contradiction argument.  Assume that the estimate (\ref{eqn:3.34}) fails.  
Then there exists a sequence $\{\bm{w}_m\}_{m=1}^\infty \subset \dot X^r(\Omega)$ such that $\|\bm{w}_m\|_{L^r(D)} \equiv 1$ for all $m =1, 2, \ldots$ and such that 
$\nabla \bm{w}_m \to \bm{0}$ in $L^r(\Omega)$ as $m \to \infty$.  It follows from Rellich's theorem for $H^{1, r}(D)$ and the weak compactness of $\dot X^r(\Omega)$ that there exists a 
subsequence of $\{\bm{w}_m\}_{m=1}^\infty$, denote again by $\{\bm{w}_m\}_{m=1}^\infty$, and a function $\bm{w} \in \dot X^r(\Omega)$ 
such that 
\begin{equation*}
\nabla \bm{w}_m \rightharpoonup \nabla \bm{w} \quad\mbox{ weakly in } L^r(\Omega) \mbox{ and } \bm{w}_m \to \bm{w} \mbox{ strongly in } L^r(D)
\end{equation*}
as $m \to \infty$. Then $\nabla \bm{w} =\bm{0}$ in $\Omega$, which implies that $\bm{w} = \bm{c}$ in $\Omega$ for a constant vector $\bm{c} \in \re^3$.  
On the other hand, there exists some $x_j\in \pt\Om$ such that 
$\bm{\nu}_{x_j} = \bm{e}_j$ for $j=1, 2, 3$, where $\bm{\nu}_{x_j}$ denotes the unit outer normal to $\pt\Om$ at $x_j$ and where $\{\bm{e}_1, \bm{e}_2, \bm{e}_3\}$ is the same as in Proposition \ref{pr:3.7}.   
Since $\bm{w}\cdot\bm{\nu}|_{\pt\Om} = 0$, it follows that $\bm{c}=0$, which, however,  contradicts $\|\bm{w}\|_{L^r(D)} =1$.  
This proves Proposition \ref{pr:3.8}.  \qed

\medskip\noindent
We consider next $\dive$-$\rot$ estimate in $\dot X(\Omega)$ and $\dot V(\Omega)$ corresponding to Propositions \ref{pr:3.5} and \ref{pr:3.6}.   
\par
%\medskip
\begin{lemma}\label{lem:3.9} Let $\Omega$ be as in the Assumption,and let $D$ and $R$ be 
as above. 
\par
{\rm (i)} For $1 < r < \infty$, there exists a constant $C = C(\Omega, r, R)$ such that  
\begin{equation}\label{eqn:3.35}
\|\nabla \bm{u}\|_{L^r} \le C(\|\rot \bm{u}\|_{L^r} + \|\dive \bm{u}\|_{L^r} + \|\bm{u}\|_{L^r(D)})
\end{equation}
holds for all $\bm{u} \in \dot X^r(\Omega)$ and all $\bm{u} \in \dot V^r(\Omega)$.  
\par
{\rm (ii)} Let $\bm{u} \in \dot X^q(\Omega)$ {\rm(}resp. $\bm{u} \in \dot V^q(\Omega)${\rm)} for some $1 < q < \infty$.  If $\dive \bm{u} \in L^r(\Omega)$ and $\rot \bm{u} \in L^r(\Omega)$ 
for some $1 < r < \infty$, then $\bm{u} \in \dot X^r(\Omega)${\rm(}resp. $\bm{u} \in \dot V^r(\Omega)${\rm )} 
satisfying  (\ref{eqn:3.35}).  
\end{lemma}

\noindent
{\it Proof.} (i) Let $\eta \in C^{\infty}_0(\re^3)$ be a cut-off function satisfying $0 \le \eta(x) \le 1$ for all $x \in \re^3$, 
$\eta(x) = 1$ for $|x| \le R+1$ and $\eta(x) =0$ for $|x| \ge R+2$.  Decompose $\bm{u} \in \dot X^r(\Om)$ in such a way that 
\begin{equation}\label{eqn:3.36}
\bm{u}(x) = \eta(x)\bm{u}(x) + (1-\eta(x))\bm{u}(x) \equiv \bm{u}_1(x) + \bm{u}_2(x), 
\quad x \in \Omega. 
\end{equation}
We may regard $\bm{u}_1$ and $\bm{u}_2$ as functions defined on $D$ and on $\re^3$, respectively.  We first handle $\bm{u}_1 \in X^r(D)$.  Since 
\begin{equation}\label{eqn:3.37}
\dive \bm{u}_1 = \eta\, \dive \bm{u} + \nabla\eta\cdot \bm{u}, 
\quad
\rot \bm{u}_1 = \eta\,\rot \bm{u} + \nabla \eta\times \bm{u}
\end{equation}
and since $\mbox{supp } \eta \subset D$, $\mbox{supp }\nabla \eta \subset B_{R+2}\setminus B_{R+1}$, 
it follows from (\ref{eqn:3.11}) that 
\begin{eqnarray}
\|\nabla \bm{u}_1\|_{L^r(D)} 
&\le& C(\|\rot \bm{u}_1\|_{L^r(D)} + \|\dive \bm{u}_1\|_{L^r(D)} + \|\bm{u}_1\|_{L^r(D)}) 
\nonumber \\
&\le&  C(\|\rot \bm{u}\|_{L^r} + \|\dive \bm{u}\|_{L^r}+ \|\bm{u}\|_{L^r(D)}),   
\label{eqn:3.38}
\end{eqnarray} 
where $C = C(\Om,r, R)$.  
\par
We next consider $\bm{u}_2 \in \dot H^{1, r}(\re^3)$.  Since  
\begin{eqnarray*}
&& |(\rot \bm{u}_2, \rot \bm{\varphi})_{\re^3} + (\dive \bm{u}_2, \dive \bm{\varphi})_{\re^3}| \\
&=& |((1-\eta)\rot \bm{u} -\nabla\eta\times \bm{u}, \rot \bm{\varphi})_{\re^3} + 
((1-\eta)\dive \bm{u} - \nabla\eta\cdot \bm{u}, \dive \bm{\varphi})_{\re^3}| \\
&\le & C(\|\rot \bm{u}\|_{L^r} + \|\bm{u}\|_{L^r(D)})\|\rot \bm{\varphi}\|_{L^{r^\prime}(\re^3)} 
+C(\|\dive \bm{u}\|_{L^r}+ \|\bm{u}\|_{L^r(D)})\|\dive \bm{\varphi}\|_{L^{r^\prime}(\re^3)}  \\
&\le & C(\|\rot \bm{u}\|_{L^r} + \|\dive \bm{u}\|_{L^r} + \|\bm{u}\|_{L^r(D)})
\|\nabla \bm{\varphi}\|_{L^{r^\prime}(\re^3)}
\end{eqnarray*}
for all $\bm{\varphi} \in C^\infty_0(\re^3)$, it follows  from (\ref{eqn:3.20}) that 
\begin{equation}\label{eqn:3.39}
\|\nabla \bm{u}_2\|_{L^r(\re^3)} 
\le  C(\|\rot \bm{u}\|_{L^r} + \|\dive \bm{u}\|_{L^r}+ \|\bm{u}\|_{L^r(D)}). 
\end{equation}
Since $\mbox{supp }\bm{u}_1 \subset D \subset \Om$, $\mbox{supp }\bm{u}_2 \subset \Om$, 
it follows from (\ref{eqn:3.36}), (\ref{eqn:3.38}) and (\ref{eqn:3.39}) that 
\begin{eqnarray*}
\|\nabla \bm{u}\|_{L^r} &\le& \|\nabla \bm{u}_1\|_{L^r(D)} + \|\nabla \bm{u}_2\|_{L^r(\re^3)}\\
&\le & C(\|\rot \bm{u}\|_{L^r} + \|\dive \bm{u}\|_{L^r}+ \|\bm{u}\|_{L^r(D)}) 
\end{eqnarray*}
with $C=C(\Omega, r, R)$, which yields (\ref{eqn:3.35}) for $\bm{u}\in \dot X^r(\Omega)$. The proof of (\ref{eqn:3.35}) for $\bm{u}\in \dot V^r(\Omega)$ is quite similar and hence omitted.  
\par
%\noindent
{\rm (ii)} The proof uses again the cut-off method.  
Given $\bm{u} \in \dot X^q(\Omega)$,  
the essential point is to derive estimate $\|\bm{u}\|_{L^r(D)} \le C\|\bm{u}\|_{H^{1, q}(D)}$.  
To this end, we consider three cases for $q$, namely for $3 \le q < \infty$, $3/2 \le q < 3$ 
and $1 < q < 3/2$.  
\par
%\noindent
{\it Case 1.} Let $3 \le q < \infty$. Proposition \ref{pr:3.8} and Sobolev's inequality imply 
$\bm{u} \in L^r(D)$ and  
\begin{equation}\label{eqn:3.40}
\|\bm{u}\|_{L^r(D)} \le C\|\bm{u}\|_{H^{1, q}(D)} \le C\|\nabla \bm{u}\|_{L^q}
\end{equation}
for some  $C = C(\Omega, r, R)$. We write $\bm{u} = \bm{u}_1 + \bm{u}_2$ as is (\ref{eqn:3.36}).     
Since $\dive \bm{u}\in L^r(\Omega)$ and $\rot \bm{u} \in L^r(\Omega)$, and since $\mbox{supp }\eta \subset D$ and $\mbox{supp }\nabla \eta \subset B_{R+2}\setminus B_{R+1}$, 
we obtain by (\ref{eqn:3.37}) and (\ref{eqn:3.40}) that $\dive \bm{u}_1 \in L^r(D)$ and $\rot \bm{u}_1 \in L^r(D)$.  
Hence, by Proposition \ref{pr:3.5} (ii)  
\begin{equation}\label{eqn:3.41}
\bm{u}_1 \in X^r(D).
\end{equation}     
Since 
\begin{equation}\label{eqn:3.42}
\dive \bm{u}_2 = (1-\eta) \dive \bm{u} - \nabla\eta\cdot \bm{u}, 
\quad
\rot \bm{u}_2 = (1-\eta)\rot \bm{u} - \nabla \eta\times \bm{u}, 
\end{equation}
we obtain  by (\ref{eqn:3.40}) that 
\begin{eqnarray*}
\|\dive \bm{u}_2\|_{L^r(\re^3)} 
&\le& \|(1-\eta)\dive \bm{u}\|_{L^r(\re^3)} + \|\nabla \eta \cdot \bm{u}\|_{L^r(\re^3)}  
\le  C(\|\dive \bm{u}\|_{L^r} + \|\bm{u}\|_{L^r(D)}) \\
&\le&  C(\|\dive \bm{u}\|_{L^r} + \|\nabla \bm{u}\|_{L^q}) \\
\|\rot \bm{u}_2\|_{L^r(\re^3)} 
&\le& \|(1-\eta)\rot \bm{u}\|_{L^r(\re^3)} + \|\nabla \eta\times \bm{u}\|_{L^r(\re^3)}  
\le  C(\|\rot \bm{u}\|_{L^r} + \|\bm{u}\|_{L^r(D)}) \\
&\le & C(\|\rot \bm{u}\|_{L^r} + \|\nabla \bm{u}\|_{L^q}).   
\end{eqnarray*}
Hence, by Proposition \ref{pr:3.6} (ii)  
\begin{equation}\label{eqn:3.43}
\bm{u}_2 \in \dot H^{1, r}(\re^3).
\end{equation}     
Combining  (\ref{eqn:3.41}) with (\ref{eqn:3.43}) we  conclude that $\bm{u} \in \dot X^r(\Omega)$.  
\par
%\noindent
{\it Case 2}.  Let $3/2 \le q <3$.  
If $r \le q$, then $\bm{u} \in L^r(D)$ with estimate (\ref{eqn:3.40}), and hence the argument given in Case 1 is valid  to derive $\bm{u} \in \dot X^r(\Omega)$.  
We hence consider the case  $q < r$. Assume first that  $1/q_1\equiv 1/q - 1/3 \le 1/r < 1/q$.   
Sobolev's inequality implies  (\ref{eqn:3.40}) and hence $\bm{u} \in \dot X^r(\Omega)$.

We next deal with the case $1/r < 1/q_1$.  Since $\dive \bm{u} \in L^q(\Omega)\cap L^r(\Omega)$ and 
$\rot \bm{u} \in L^q(\Omega)\cap L^r(\Omega)$, and since $q < q_1 < r$ we see that $\dive \bm{u} \in L^{q_1}(\Om)$ and $\rot \bm{u} \in L^{q_1}(\Om)$. 
Hence, by the above argument, $\bm{u} \in \dot X^{q_1}(\Omega)$.  
Since $q_1 \ge 3$ we obtain (\ref{eqn:3.40}) with $q$ replaced by $q_1$ so that we may reduce our case to  Case 1 
to obtain $\bm{u} \in \dot X^r(\Omega)$.  
\par
%\noindent
{\it Case 3}. Let $1 < q < 3/2$.  If $r \le q$, then (\ref{eqn:3.40}) holds and $\bm{u} \in \dot X^r(\Omega)$. If $q < r \le q_1$, Sobolev's inequality implies (\ref{eqn:3.40})  
and hence $\bm{u} \in \dot X^r(\Omega)$. It remains to consider the case $q_1 < r < \infty$. The same argument as in the Case 2 is valid and yields  
$\bm{u} \in \dot X^{q_1}(\Om)$.  If $r\le q_2$ with $1/q_2\equiv 1/q_1-1/3 = 1/q - 2/3$, then (\ref{eqn:3.40}) holds with $q$ replaced by $q_1$, which yields  $\bm{u} \in \dot X^r(\Omega)$.  
Finally, consider the case  $q_2 < r < \infty$.  Since $\bm{u} \in \dot X^{q_1}(\Om)$, $\dive \bm{u} \in L^{q_2}(\Omega)$ and 
$\rot \bm{u} \in L^{q_2}(\Omega)$, (\ref{eqn:3.40}) with $r$ and $q$ replaced by $q_2$ and $q_1$, implies $\bm{u} \in \dot X^{q_2}(\Omega)$.  
Defining $1 < q_3 < q_2$ by $1/q_3 \equiv 1/r + 1/3$ we have $u \in H^{1, q_3}(D)$.  
Hence by (\ref{eqn:3.34}) and by Sobolev's inequality $\bm{u} \in L^r(D)$ satisfying  
$$
\|\bm{u}\|_{L^r(D)} \le C\|\bm{u}\|_{H^{1, q_3}(D)} 
\le C\|\bm{u}\|_{H^{1. q_2}(D)} \le C\|\nabla \bm{u}\|_{L^{q_2}}. 
$$ 
The same argument as in the Case 1 is valid now, because (\ref{eqn:3.40}) holds with $q$ replaced by $q_2$. Hence, $\bm{u} \in \dot X^r(\Omega)$ and the proof of ii) for 
$\bm{u} \in \dot X^q(\Om)$ is complete. The proof for $\bm{u} \in \dot V^q(\Om)$ is similar and hence omitted. 
\qed 
\par
\bigskip
We now define the harmonic vector spaces $\dot X^r(\Omega)$ and $\dot V^r(\Omega)$ for $1 < r < \infty$ by
\begin{eqnarray}
&& \dot X^r_{\tiny{\mbox{har}}}(\Omega)\equiv \{\bm{h}\in \dot X^r(\Omega); 
\dive \bm{h}=0,\, \rot \bm{h}=0\}, 
\label{eqn:3.44}  \\
&& 
\dot V^r_{\tiny{\mbox{har}}}(\Omega)\equiv \{\bm{h}\in \dot V^r(\Omega); 
\dive \bm{h}=0,\, \rot \bm{h}=0\}.
\label{eqn:3.45} 
\end{eqnarray}
The following lemma is a consequence of Lemma \ref{lem:3.9}.
\par
%\bigskip
\begin{lemma}\label{lem:3.10} 
Let $\Omega$ be as in the Assumption. If $1 < r < \infty$, then  
\begin{equation*}
\dot X^r_{\tiny{\mbox{\rm har}}}(\Omega)= \dot X^2_{\tiny{\mbox{\rm har}}}(\Omega)\equiv \dot X_{\tiny{\mbox{\rm har}}}(\Omega) \quad \mbox{ as well as } \quad
 \dot V^r_{\tiny{\mbox{\rm har}}}(\Omega)= \dot V^2_{\tiny{\mbox{\rm har}}}(\Omega)\equiv \dot V_{\tiny{\mbox{\rm har}}}(\Omega)   
\end{equation*}
and $\dot X_{\tiny{\mbox{\rm har}}}(\Omega)$ and $\dot V_{\tiny{\mbox{\rm har}}}(\Omega)$ are finite dimensional vector spaces.  
\end{lemma}  

\medskip\noindent
{\it Proof.} The first two assertions follow from Lemma \ref{lem:3.9} (ii). 
In order to prove that $\dim\dot X^r_{\tiny{\mbox{har}}}(\Omega) < \infty$ note that   
since $\dot X^r_{\tiny{\mbox{har}}}(\Omega)$ is a closed subspace in $\dot X^r(\Omega)$, it suffices to show that the unit sphere of $\dot X^r_{\tiny{\mbox{har}}}(\Omega)$ 
is compact in the norm of $\dot X^r(\Omega)$. To this end, suppose that $\{\bm{u}_m\}_{m=1}^\infty \subset \dot X_{\tiny{\mbox{\rm har}}}(\Omega)$ 
satisfies  $\|\nabla \bm{u}_m\|_{L^2}\equiv 1$ for all $m=1, 2, \cdots$. 
By the weak compactness, there exist a subsequence of $\{\bm{u}_m\}_{m=1}^\infty$, denoted again by $\{\bm{u}_m\}_{m=1}^\infty$, and a function $\bm{u} \in \dot X_{\tiny{\mbox{har}}}(\Omega)$ such that 
\begin{equation}\label{eqn:3.46}
\nabla \bm{u}_m \rightharpoonup \nabla \bm{u} 
\quad
\mbox{weakly in $L^2(\Om)$ as $m\to \infty$}. 
\end{equation}
On the other hand, we see by Proposition \ref{pr:3.8} that $\{\bm{u}_m\}_{m=1}^\infty$ is bounded in $H^{1, 2}(D)$. Since $D$ is a bounded domain in $\re^3$, by Rellich's  
theorem we may assume that 
$$
\bm{u}_m \to \bm{u} 
\quad
\mbox{strongly in $L^r(D)$ as $m\to\infty$}.       
$$
Since $\dive \bm{u}_m=0$ and $\rot \bm{u}_m = 0$ for all $m=1, 2, \cdots$, 
Lemma \ref{lem:3.9} (i) implies that 
$\{\nabla \bm{u}_m\}_{m=1}^\infty$ is a Cauchy sequence in $L^2(\Omega)$ and 
by (\ref{eqn:3.46}) it follows that 
$$
\nabla \bm{u}_m \to \nabla \bm{u} 
\quad
\mbox{strongly in $L^2(\Om)$ as $m\to \infty$}. 
$$
Since $\{\bm{u}_m\}_{m=1}^\infty$ is an arbitrary sequence of the unit sphere in $\dot X_{\tiny{\mbox{har}}}(\Omega)$, we conclude that the dimension of $\dot X_{\tiny{\mbox{har}}}(\Omega)$ is finite.  
The proof of $\dim \dot V_{\tiny{\mbox{har}}}(\Omega) < \infty$ is parallel to the one given above and hence omitted.   \qed

\medskip\noindent
Following Lemma \ref{lem:3.10}, we may assume that $\dim \dot X_{\tiny{\mbox{har}}}(\Omega) = N$ and $\dim \dot V_{\tiny{\mbox{har}}}(\Omega)=L$.  
Hence, there are a basis $\{\bm{\varphi}_1, \cdots, \bm{\varphi}_N\}$ of $\dot X_{\tiny{\mbox{har}}}(\Omega)$ 
and a basis $\{\bm{\psi}_1, \cdots, \bm{\psi}_L\}$ of $\dot V_{\tiny{\mbox{har}}}(\Omega)$ such that 
\begin{equation}\label{eqn:3.47}
(\nabla\bm{\varphi}_j, \nabla\bm{\varphi}_k ) = \delta_{jk}, \quad j, k =1, \cdots N, 
\quad
(\nabla \bm{\psi}_l, \nabla\bm{\psi}_m ) = \delta_{lm}, \quad l, m =1, \cdots L,    
\end{equation}
where $\{\delta_{jk}\}_{1\le j, k\le N}$ and $\{\delta_{lm}\}_{1\le l, m \le L}$ denote 
Kronecker's symbols. 
Then for every $1< r < \infty$, there exist  closed subspaces $\widetilde X^r(\Omega)$ and 
$\widetilde V^r(\Omega)$ 
of $\dot X^r(\Omega)$ and $\dot V^r(\Omega)$, respectively, such that 
\begin{equation}\label{eqn:3.48}
\dot X^r(\Omega) = \dot X_{\tiny{\mbox{har}}}(\Omega) \oplus \widetilde X^r(\Omega), 
\quad
\dot V^r(\Omega) = \dot V_{\tiny{\mbox{har}}}(\Omega) \oplus \widetilde V^r(\Omega)
\quad
\mbox{(direct sum)}.  
\end{equation}
The following lemma corresponds to Lemma \ref{lem:3.9} (i).
\par
%\medskip
\begin{lemma}\label{lem:3.11} Let $\Omega$ be as in the Assumption and let  $1 < r < \infty$. \\ 
{\rm (i)} There exists a constant $C = C(\Omega,r)$ such that   
\begin{eqnarray} 
&& \|\nabla \bm{u}\|_{L^r} 
\le C(\|\rot \bm{u}\|_{L^r} + \|\dive \bm{u}\|_{L^r} 
+ \sum_{j=1}^N|(\nabla \bm{u}, \nabla \bm{\varphi}_j)|)
\quad
\mbox{ for all $\bm{u} \in \dot X^r(\Omega)$},  \label{eqn:3.49} \\
&& \|\nabla \bm{u}\|_{L^r} 
\le C(\|\rot \bm{u}\|_{L^r} + \|\dive \bm{u}\|_{L^r} 
+ \sum_{l=1}^L|(\nabla \bm{u}, \nabla \bm{\psi}_l)|)
\quad
\mbox{ for all $\bm{u} \in \dot V^r(\Omega)$}.   \label{eqn:3.50}
\end{eqnarray}
{\rm (ii)} There exists a constant $C = C(\Omega,r)$ such that 
\begin{equation}\label{eqn:3.51}
\|\nabla \bm{u}\|_{L^r} \le C(\|\rot \bm{u}\|_{L^r} + \|\dive \bm{u}\|_{L^r})
\quad
\mbox{for all $\bm{u} \in \tilde X^r(\Omega)$ and for all $\bm{u} \in \tilde V^r(\Omega)$}. 
\end{equation}
\end{lemma}
%\medskip\noindent
{\it Proof.} (i) In order to prove (\ref{eqn:3.49}), we make use of a contradiction argument.  
Assume that the estimate (\ref{eqn:3.49}) fails.  Then there is a sequence $\{\bm{u}_m\}_{m=1}^\infty$ in $\dot X^r(\Omega)$ such that 
\begin{eqnarray*}
&& \|\nabla \bm{u}_m\|_{L^r} \equiv 1 \quad\mbox{for all $m \in \N$}, \\
&& \|\rot \bm{u}_m\|_{L^r} + \|\dive \bm{u}_m\|_{L^r} 
+ \sum_{j=1}^N|(\nabla \bm{u}_m, \nabla\bm{\varphi}_j)| \to 0 
\quad\mbox{as $m \to \infty$}.  
\end{eqnarray*} 
By the weak compactness, there exists a subsequence of $\{\bm{u}_m\}_{m=1}^\infty$, denoted again  by $\{\bm{u}_m\}_{m=1}^\infty$, and a function $\bm{u} \in \dot X^r(\Omega)$ such that 
$$
\nabla \bm{u}_m \rightharpoonup \nabla \bm{u}
\quad
\mbox{weakly in $L^r(\Omega)$}.  
$$  
By Proposition \ref{pr:3.8}, $\{\bm{u}_m\}_{m=1}^\infty$ is bounded in $H^{1, r}(D)$. By Rellich's theorem we may assume that 
$$
\bm{u}_m \to \bm{u} 
\quad\mbox{strongly in $L^r(D)$ as $m \to \infty$}.  
$$ 
By Lemma \ref{lem:3.9} (i), $\{\nabla \bm{u}_m\}_{m=1}^\infty$ is a Cauchy sequence in $L^r(\Omega)$, which implies that 
$\nabla \bm{u}_m \to \nabla \bm{u}$  strongly in $L^r(\Omega)$ as $m\to \infty$. 
Hence, $\dive \bm{u}=0$ and $\rot \bm{u}=\bm{0}$, and Lemma \ref{lem:3.10} implies  $\bm{u}\in \dot X_{\tiny{\mbox{har}}}(\Omega)$. 
Moreover, since $(\nabla \bm{u}, \nabla \bm{\varphi}_j) =0$ for $j =1, \cdots, N$, $\bm{u}=\bm{0}$, 
which contradicts $\|\nabla \bm{u}\|_{L^r} =1$.  The proof of (\ref{eqn:3.50}) is similar to that of (\ref{eqn:3.49}) and hence omitted.  
\par
%\noindent
(ii) Assume that (\ref{eqn:3.51}) fails. Then there is a sequence $\{\bm{u}_m\}_{m=1}^\infty$ in $\tilde X^r(\Omega)$ such that 
$$
\|\nabla \bm{u}_m\|_{L^r} \equiv 1\quad\mbox{ for all $m \in \N$, and }
\quad
\lim_{m\to \infty} (\|\rot \bm{u}_m\|_{L^r} + \|\dive \bm{u}_m\|_{L^r}) = 0. 
$$
Then by the same compactness argument as in the above, we may assume that there is a function $\bm{u} \in \dot X^r(\Omega)$ such that 
$$
\lim_{m\to\infty}\|\nabla \bm{u}_m - \nabla \bm{u}\|_{L^r} =0.   
$$    
Since $\widetilde X^r(\Omega)$ is closed in $\dot X^r(\Omega)$, 
we obtain $\bm{u} \in \widetilde X^r(\Omega)$.  
On the other hand, $\dive \bm{u}=0$ and $\rot \bm{u}=\bm{0}$, which means that $\bm{u} \in \dot X_{\tiny{\mbox{har}}}(\Omega)$.   
Since $\widetilde X^r(\Omega) \cap \dot X_{\tiny{\mbox{har}}}(\Omega) = \{\bm{0}\}$, 
we see that $\bm{u}=\bm{0}$, which contradicts $\|\nabla \bm{u}\|_{L^r}=1$. 
The proof of (\ref{eqn:3.51}) for $\bm{u} \in \tilde V^r(\Omega)$ follows 
similarly from the above lines, and hence details may be omitted.    \qed
\par
\bigskip 
\subsection {Functionals on $\re^3$ supported in $\Omega$} \mbox{}
\par 
%Finally in this section, we establish a norm preserving estimate 
%in such a way that the continuous functional 
%supported in $\Omega$ can be extended to that on the whole space $\re^3$.    
Assume that $\Omega_0$ is a subdomain of $\Omega$ satisfying $\overline{\Omega_0} \subset \Omega$. Consider a further  subdomain $\Omega_1 \subset \Omega$ in such a way 
that $\overline{\Omega_0} \subset \Omega_1$ and that $K\equiv \Omega\setminus\overline{\Omega_1}$ is a bounded domain in $\re^3$.  
Take a cut-off function $\zeta \in C^\infty(\re^3)$ with $0 \le \zeta(x) \le 1$ such that 
\begin{equation}\label{eqn:3.52}
\zeta(x) =
\left\{
\begin{array}{ll}
& 1 \quad\mbox{for $x \in \Omega_1$}, \\
& 0  \quad\mbox{for $x \in \re^3\setminus\Omega$}.  
\end{array}
\right.
\end{equation} 
We show that the space of linear functionals on $\dot X^{r^\prime}(\Omega)$ and $\dot V^{r^\prime}(\Omega)$ supported on $\overline{\Omega_0} \subset \Omega$ 
are continuously embedded into those on $\widehat H^{1, r}(\re^3)$ which is the closure of $C^\infty_0(\re^3)$ with respect to the 
homogeneous norm $\|\nabla \bm{u}\|_{L^r(\re^3)}$. Similarly to Proposition \ref{pr:3.1}, we have the following concrete characterization of 
$\widehat H^{1, r}(\re^n)$ as 
$$
\widehat H^{1, r}(\re^3) = 
\left\{
\begin{array}{ll}
& \{[\bm{u}]; \bm{u} \in L^{r_\ast}(\re^3); \nabla \bm{u} \in L^r(\re^3)\}
\quad\mbox{for $1 < r < 3$}, \\
& 
\{[\bm{u}] ; \bm{u} \in  L^r_{\tiny{\mbox{loc}}}(\re^3),  \nabla \bm{u} \in L^r(\re^3)\}
\quad\mbox{for $3 \le r < \infty$},
\end{array}
\right.
$$
where $[\bm{u}]$ denotes the set of equivalent class of $\bm{u}$ such that 
$\bm{v} \in [\bm{u}]$ implies that 
$\bm{u}-\bm{v}\equiv \mbox{const}$ in $\re^3$. 
Hence, 
\begin{equation}\label{eqn:3.53}
\dot H^{1, r}(\re^3) = \widehat H^{1, r}(\re^3)
\quad\mbox{ for all $1 < r < \infty$.} 
\end{equation}
Indeed, by \cite[Corollary 2.2 (ii)]{GiSo} if $1 < r < 3$, for $\bm{u} \in L^r_{\tiny{\mbox{loc}}}(\re^3)$ with $\nabla \bm{u} \in L^r(\re^3)$ 
there is a constant vector $\bm{c} \in \re$ such that $\bm{u} + \bm{c} \in L^{r_\ast}(\re^3)$, which yields (\ref{eqn:3.53}). 
\begin{proposition}\label{pr:3.12} Let $\Omega$ be as in the Assumption.  
Suppose that $\Omega_0$ is a subdomain satisfying $\overline{\Omega_0} \subset \Omega$ and that 
$\zeta$ is a function as in (\ref{eqn:3.52}).  
%Assume that $f \in \dot X^{r^{\prime}}(\Omega)^\ast$ with $\mbox{supp} f \subset \overline{\Om_0}$. 
Let $f \in \dot X^{r^{\prime}}(\Omega)^\ast$ or $f \in \dot V^{r^{\prime}}(\Omega)^\ast$ for $1 < r < \infty$.  
Then the functional $\tilde f$ defined by 
\begin{equation}\label{eqn:3.54}
\langle \tilde f, \bm{\Phi} \rangle_{\re^n}\equiv \langle f, \zeta\bm{\Phi}\rangle 
\quad\mbox{ for $\bm{\Phi} \in \hat H^{1, r^{\prime}}(\re^3)$}
\end{equation}   
may be regarded as an element of $\widehat H^{1, r^{\prime}}(\re^3)^\ast$ satisfying  
\begin{equation}\label{eqn:3.55}
\|\tilde f\|_{\hat H^{1, r^{\prime}}(\re^3)^\ast} 
\le 
\left\{ 
\begin{array}{ll}
& C \|f\|_{\dot X^{r^{\prime}}(\Omega)^\ast}
\quad\mbox{for $f \in \dot X^{r^{\prime}}(\Omega)^\ast$},  \\
& C\|f\|_{\dot V^{r^{\prime}}(\Omega)^\ast}
\quad\mbox{for $f \in \dot V^{r^{\prime}}(\Omega)^\ast$},   
\end{array}
\right.
\end{equation}
for some  constant $C = C(\Omega,\Omega_0,r)$ independent of $f$.  
\end{proposition}
\par
To be precise let us note that in (\ref{eqn:3.54}) the bracket $\langle\cdot, \cdot\rangle_{\re^3}$ on the left hand side  denotes the duality pairing between 
$\widehat H^{1, r^{\prime}}(\re^3)^\ast$ and $\widehat H^{1, r^{\prime}}(\re^3)$, 
while the bracket $\langle\cdot, \cdot\rangle$  on the right hand side denotes the one  between 
$\dot X^{r^{\prime}}(\Omega)^\ast${\rm(}resp.$\dot V^{r^{\prime}}(\Omega)^\ast${\rm )} 
and $\dot X^{r^{\prime}}(\Omega)${\rm(}resp.$\dot V^{r^{\prime}}(\Omega)${\rm )}.  
\par
\vspace{2mm}
%\bigskip
\noindent
{\it Proof of Proposition \ref{pr:3.12}.} We give here only a detailed proof  for $f \in \dot X^{r^\prime}(\Omega)^\ast$; the proof for $f \in\dot V^{r^\prime}(\Omega)^\ast $ is similar. Note that 
$\zeta \bm{\Phi} \in \dot X^{r^\prime}(\Omega)\cap \dot V^{r^\prime}(\Omega)$ for all $\bm{\Phi} \in \widehat H^{1, r^\prime}(\re^3)$. 
\par 
Consider first the case when $3/2 < r < \infty$. 
It  follows from (\ref{eqn:3.53}) and Sobolev's inequality that 
\begin{equation}\label{eqn:3.56}
\|\bm{\Phi}\|_{L^q(\re^3)} \le C\|\nabla \bm{\Phi}\|_{L^{r^\prime}(\re^3)} 
\quad\mbox{for all $\bm{\Phi}\in \widehat H^{1, r^\prime}(\re^3)$}, 
\end{equation}
where $1/q \equiv 1/r^\prime - 1/3$. Since $\mbox{supp} \nabla \zeta \subset K$ and since $K$ is a bounded domain in $\re^3$, it follows from (\ref{eqn:3.54}) and (\ref{eqn:3.56}) that  
\begin{eqnarray}
|\langle \tilde f, \bm{\Phi}\rangle_{\re^3}|
&\le& \|f\|_{\dot X^{r^\prime}(\Omega)^\ast}\|\zeta \bm{\Phi}\|_{\dot X^{r^\prime}(\Omega)} 
\nonumber \\ 
&=&  \|f\|_{\dot X^{r^\prime}(\Omega)^\ast}\|\nabla (\zeta \bm{\Phi})\|_{L^{r^\prime}}\nonumber\\
&\le&  \|f\|_{\dot X^{r^\prime}(\Omega)^\ast}(\|\nabla \zeta\cdot\bm{\Phi}\|_{L^{r^\prime}(K)} 
+ \|\nabla \bm{\Phi}\|_{L^{r^\prime}}) \nonumber\\ 
&\le& C \|f\|_{\dot X^{r^\prime}(\Omega)^\ast}(\|\bm{\Phi}\|_{L^q(K)} 
+ \|\nabla \bm{\Phi}\|_{L^{r^\prime}(\re^3)} )
\nonumber \\
&\le&  C\|f\|_{\dot X^{r^\prime}(\Omega)^\ast}\|\nabla \bm{\Phi}\|_{L^{r^\prime}(\re^3)} 
\label{eqn:3.57}
\end{eqnarray}
for all $\bm{\Phi} \in \widehat H^{1, r^\prime}(\re^3)$. 
This implies  $\tilde f \in \widehat H^{1, r^\prime}(\re^3)^\ast$ as well as (\ref{eqn:3.55}).  
\par
We consider next the case when $1 < r \leqq 3/2$. It follows from \cite[Lemma 2.5]{KoSo} that the space $S_K$ defined by 
$S_K \equiv \left\{\bm{\varphi} \in C^\infty_0(\re^3); \int_K\bm{\varphi}(x)dx = \bm{0}\right\}$  
is {\it dense} in $\widehat H^{1, r^\prime}(\re^3)$. Poinc\'are's inequality implies  
$$
\|\bm{\varphi} \|_{L^{r^\prime}(K)} \le C\|\nabla \bm{\varphi}\|_{L^{r^\prime}(K)} 
\le C\|\nabla \bm{\varphi}\|_{L^{r^\prime}(\re^3)}
\quad\mbox{for all $\varphi \in S_K$}, 
$$  
for some  $C = C(r, K)$. Hence, similarly to (\ref{eqn:3.57}) we obtain  
$$
|\langle \tilde f, \bm{\varphi}\rangle_{\re^3}| 
\le C\|f\|_{\dot X^{r^\prime}(\Omega)^\ast}(\|\bm{\varphi}\|_{L^{r^\prime}(K)} 
+ \|\nabla \bm{\varphi}\|_{L^{r^\prime}(\re^3)}) 
\le C\|f\|_{\dot X^{r^\prime}(\Omega)^\ast}\|\nabla \bm{\varphi}\|_{L^{r^\prime}(\re^3)}
$$
for all $\bm{\varphi} \in S_K$. 
Since $S_K$ is dense in $\widehat H^{1, r^\prime}(\re^3)$, 
we obtain  $\tilde f \in \widehat H^{1, r^\prime}(\re^3)^\ast$ 
as well as  (\ref{eqn:3.55}).   \qed

\section{Construction of the vector potential}
In this section to prove the existence of the vector potential $\bm{w}$ in (\ref{eqn:2.6}) and (\ref{eqn:2.10}).

\begin{thm}\label{thm:4.1} Let $\Omega$ be as in the Assumption and let  $1 < r < \infty$. 
\par
%\noindent
{\rm (i)} For every $\bm{u} \in L^r(\Omega)$ there exists $\bm{w} \in \dot X^r_\sg(\Omega)$ such that 
\begin{equation}\label{eqn:4.1}
(\rot \bm{w}, \rot \bm{\Phi} ) = (\bm{u}, \rot \bm{\Phi})
\end{equation}
holds for all $\bm{\Phi} \in \dot X^{r^\prime}(\Omega)$.  Moreover, there exists a constant $C = C(\Omega,r)$ such that  
\begin{equation}\label{eqn:4.2}
\|\nabla \bm{w}\|_{L^r} \le C\|\bm{u}\|_{L^r}.
\end{equation}
\par
{\rm (ii)}  For every $\bm{u} \in L^r(\Omega)$ there exists $\bm{w} \in \dot V^r_\sg(\Omega)$ such that 
\begin{equation}\label{eqn:4.3}
(\rot \bm{w}, \rot \bm{\Psi} ) = (\bm{u}, \rot \bm{\Psi})
\end{equation}
holds for all $\bm{\Psi} \in \dot V^{r^\prime}(\Omega)$.   Moreover, there exists a constant $C = C(\Omega,r)$ such that 
\begin{equation}\label{eqn:4.4}
\|\nabla\bm{w}\|_{L^r} \le C\|\bm{u}\|_{L^r}. 
\end{equation}
\end{thm}
\par
The proof of Theorem \ref{thm:4.1} is subdivided into several steps. 
\subsection{Bilinear forms on $\dot X^r(\Omega)$ and $\dot V^r(\Omega)$ } \mbox{}
\par
In order to prove Theorem \ref{thm:4.1}, we first solve (\ref{eqn:4.1}) and (\ref{eqn:4.3}) for general $\bm{w} \in \dot X^r(\Omega)$ and $\bm{w} \in \dot V^r(\Omega)$, 
respectively.  Then, in the next step, we show that such $\bm{w}$ necessarily satisfies  $\dive \bm{w}=0$.  
For the validity of  $\bm{w} \in \dot X^r_\sg(\Omega)$ in (\ref{eqn:4.1}) it is essential that the  right hand side  is written in the form $(\bm{u}, \rot \bm{\Phi})$.  
%So is the case in (\ref{eqn:4.3}).   
Consider the bilinear forms $a_X(\cdot, \cdot)$ and $a_V(\cdot, \cdot)$ defined by 
\begin{eqnarray}
a_{X}(\bm{v}, \bm{\Phi})& =& (\rot\bm{v}, \rot \bm{\Phi}) + (\dive \bm{v}, \dive \bm{\Phi}) 
\quad
\mbox{for $\bm{v} \in \dot X^r(\Omega)$ and $\bm{\Phi} \in \dot X^{r^\prime}(\Omega)$}, \label{eqn:4.5} \\
a_{V}(\bm{w}, \bm{\Psi})& =& (\rot \bm{w}, \rot \bm{\Psi}) + (\dive \bm{w}, \dive \bm{\Psi}) 
\quad
\mbox{for $\bm{w} \in \dot V^r(\Omega)$ and $\bm{\Psi} \in \dot V^{r^\prime}(\Omega)$}.  \label{eqn:4.6}
\end{eqnarray}
The following proposition gives regularity properties for  solutions associated to the bilinear forms (\ref{eqn:4.5}) and (\ref{eqn:4.6}).  

\medskip
\begin{proposition}\label{pr:4.2} Let $\Omega$ be as in the Assumption and let  $1< q < \infty$ and $3/2 < r < \infty$.  
\par
%\noindent
{\rm (i)} Suppose that $f \in \dot X^{q^\prime}(\Om)^\ast \cap \dot X^{r^\prime}(\Omega)^\ast$. 
If $\bm{v} \in \dot X^q(\Omega)$ satisfies 
\begin{equation}\label{eqn:4.7} 
a_X(\bm{v}, \bm{\varphi}) =\langle f, \bm{\varphi} \rangle 
\quad
\mbox{for all $\bm{\varphi} \in {\mathcal X}(\Omega)$}, 
\end{equation}  
then $\bm{v} \in \dot X^r(\Omega)$.  Here the bracket $\langle\cdot, \cdot \rangle$ denotes the duality pairing between $\dot X^{q^\prime}(\Om)^\ast$ and $\dot X^{q^\prime}(\Om)$.    
\par
\noindent
{\rm (ii)} Suppose that $f \in \dot V^{q^\prime}(\Om)^\ast \cap \dot V^{r^\prime}(\Omega)^\ast$. If $\bm{w} \in \dot V^q(\Omega)$ satisfies 
\begin{equation}\label{eqn:4.8} 
a_V(\bm{w}, \bm{\psi}) =\langle f, \bm{\psi} \rangle 
\quad
\mbox{for all $\bm{\psi} \in {\mathcal V}(\Omega)$}, 
\end{equation}  
then  $\bm{w} \in \dot V^r(\Omega)$.  Here the bracket $\langle\cdot, \cdot \rangle$ denotes the duality paring between 
$\dot V^{q^\prime}(\Om)^\ast$ and $\dot V^{q^\prime}(\Om)$. 
\end{proposition}
\par
\medskip
%\noindent
{\it Proof}. Again we only prove (i). Define two domains $\Omega_1$ and $\Omega_2$ by $\Omega_1\equiv D$ and $\Omega_2 \equiv \re^3$, respectively and define further two cut-off functions 
$\zeta_1$ and $\zeta_2$ by $\zeta_1(x) \equiv \eta(x)$ and $\zeta_2(x) \equiv 1- \eta(x)$, respectively, where $\eta\in C^\infty_0(\re^3)$ is the same cut-off function as in (\ref{eqn:3.36}). 
Since 
\begin{equation}\label{eqn:4.9}
\bm{v}(x) = \zeta_1(x)\bm{v}(x) + \zeta_2(x)\bm{v}(x)\equiv \bm{v}_1(x) + \bm{v}_2(x), 
\quad x \in \Omega,
\end{equation}   
we may regard $\bm{v}_1 \in X^q(D)$ and $\bm{v}_2 \in \dot H^{1, q}(\re^3)$.  By (\ref{eqn:4.7}) we obtain 
\begin{equation}
\int_{\Omega_i}(\rot \bm{v}_i\cdot\rot \bm{\varphi}_i + \dive \bm{v}_i\,  \dive \bm{\varphi}_i)dx              =  \langle f, \zeta_i\bm{\varphi}_i\rangle 
- \int_D(2\nabla \zeta_i\cdot\nabla \bm{v} + \Delta \zeta_i\bm{v})\cdot \bm{\varphi}_i dx, 
\quad i=1, 2 \label{eqn:4.10}
\end{equation}
for all $\bm{\varphi}_1 \in {\mathcal X}(D)$ and for all $\bm{\varphi}_2 \in C^\infty_0(\re^3)$.              
Notice that $\mbox{supp }\zeta_1 \subset D$ and that 
$\mbox{supp }\nabla \zeta_i, \mbox{supp }\Delta \zeta_i \subset B_{R+2}\setminus B_{R+1} \subset D$ for $i=1, 2$. 
Since $f \in \dot X^{r^\prime}(\Omega)^\ast$, it follows from (\ref{eqn:3.8}) that 
\begin{eqnarray} 
|\langle f, \zeta_1\bm{\varphi}_1\rangle |
&\le& \|f\|_{\dot X^{r^\prime}(\Omega)^\ast}\|\nabla (\zeta_1\bm{\varphi}_1)\|_{L^{r^\prime}} \nonumber \\
&\le&  \|f\|_{\dot X^{r^\prime}(\Omega)^\ast}(\|\zeta_1\nabla \bm{\varphi}_1\|_{L^{r^\prime}} 
+ \|\nabla \zeta_1\cdot\bm{\varphi}_1\|_{L^{r^\prime}}) \nonumber \\
&\le &C\|f\|_{\dot X^{r^\prime}(\Omega)^\ast}(\|\nabla \bm{\varphi}_1\|_{L^{r^\prime}(D)} 
+ \|\nabla \bm{\varphi}_1\|_{L^{r^\prime}(D)}) \nonumber \\
&\le& C \|f\|_{\dot X^{r^\prime}(\Omega)^\ast}\|\nabla \bm{\varphi}_1\|_{L^{r^\prime}(D)} \label{eqn:4.11}
 \end{eqnarray}
for all $\bm{\varphi}_1 \in {\mathcal X}(D)$.  By Proposition \ref{pr:3.12}
\begin{equation}\label{eqn:4.12}
|\langle f, \zeta_2\bm{\varphi}_2\rangle| 
\le C\|f\|_{\dot X^{r^\prime}(\Omega)^\ast}\|\nabla \bm{\varphi}_2\|_{L^{r^\prime}(\re^3)}
\end{equation} 
for all $\bm{\varphi}_2 \in C^\infty_0(\re^3)$.  
In order to treat  the second term in the  integral on $D$ on the right hand side  of (\ref{eqn:4.10}), we use a technique similar to the one given in the proof of Lemma \ref{lem:3.9} (ii) and 
consider the following three cases. 
\par
%\noindent
{\it Case 1}. Let $1/q - 1/3 \le 1/r < 2/3$. 
Take $1 < s < \infty$ satisfying  $1/s \equiv 1/r + 1/3$ and hence $q^\prime \le s^\prime$. 
Since $D$ is a bounded domain, it follows from (\ref{eqn:3.9}) and (\ref{eqn:3.34}) that 
\begin{eqnarray}
\left|\int_D(2\nabla \zeta_1\cdot\nabla \bm{v} + \Delta \zeta_1\bm{v})\cdot \bm{\varphi}_1 dx\right| 
&\le& C(\|\nabla \bm{v}\|_{L^q(D)} + \|\bm{v}\|_{L^q(D)})\|\bm{\varphi}_1\|_{L^{q^\prime}(D)} \nonumber \\
&\le & C(\|\nabla \bm{v}\|_{L^q} + \|\bm{v}\|_{L^q(D)}))\|\bm{\varphi}_1\|_{L^{s^\prime}(D)}\nonumber \\
&\le & C\|\nabla \bm{v}\|_{L^q}\|\nabla \bm{\varphi}_1\|_{L^{r^\prime}(D)} \label{eqn:4.13}
\end{eqnarray}  
for all $\bm{\varphi}_1 \in {\mathcal X}(D)$ and for some  $C=C(\Omega,r,q,R)$. Similarly, the Sobolev embedding 
$\dot H^{1, r^{\prime}}(\re^3) \subset L^{s^\prime}(\re^3)$ implies  
\begin{eqnarray}
\left| \int_D(2\nabla \zeta_2\cdot\nabla \bm{v} + \Delta \zeta_2 \bm{v})\cdot \bm{\varphi}_2 dx\right| 
&\le& C(\|\nabla \bm{v}\|_{L^q(D)} + \|\bm{v}\|_{L^q(D)})\|\bm{\varphi}_2\|_{L^{s^\prime}(D)} 
\nonumber \\
&\le& C\|\nabla \bm{v}\|_{L^q}\|\nabla \bm{\varphi}_2\|_{L^{r^\prime}(\re^3)}\label{eqn:4.14}
\end{eqnarray} 
for all $\bm{\varphi}_2 \in C^{\infty}_0(\re^3)$ for some  $C=C(\Omega, r, q, R)$. 
It follows from (\ref{eqn:4.10}), (\ref{eqn:4.11}), (\ref{eqn:4.13}) and 
Proposition \ref{pr:3.5} (ii) that 
\begin{equation}\label{eqn:4.15}
\bm{v}_1 \in X^r(D).  
\end{equation} 
Similarly, it follows from (\ref{eqn:4.10}), (\ref{eqn:4.12}), (\ref{eqn:4.14}) 
and Proposition \ref{pr:3.6} (ii) that
\begin{equation}\label{eqn:4.16}
\bm{v}_2 \in \dot H^{1, r}(\re^3).  
\end{equation} 
We now deduce from  (\ref{eqn:4.9}), (\ref{eqn:4.15}) and (\ref{eqn:4.16}) 
that $\bm{v} \in \dot X^r(\Omega)$ provided $1/q -1/3 \le 1/r < 2/3$.    
\par
%\noindent
{\it Case 2}. Let $1/q -2/3 \le 1/r < 1/q - 1/3 \equiv 1/q_1$.    Since $f \in \dot X^{q^\prime}(\Om)^\ast \cap \dot X^{r^\prime}(\Omega)^\ast$ 
and since $q < q_1 < r$, we obtain $f \in \dot X^{q^\prime}(\Om)^\ast \cap \dot X^{q_1^\prime}(\Omega)^\ast$.  Hence, the argument of case 1 with $r$ replaced by $q_1$ implies  
$\bm{v} \in \dot X^{q_1}(\Omega)$.  Since $f \in \dot X^{q_1^\prime}(\Om)^\ast \cap \dot X^{r^\prime}(\Omega)^\ast$, the argument of case 1 with $q$ replaced by $q_1$, implies  
$\bm{v} \in \dot X^r(\Omega)$. 
\par
%\noindent
{\it Case 3}. Let $0 < 1/r < 1/q - 2/3 \equiv 1/q_2$. Take $q_1 < q_3 < q_2$ so that $1/q_3 \equiv 1/r + 1/3$. The arguments of  case 1 and case 2 imply 
$\bm{v} \in \dot X^{q_3}(\Omega)$.  Since $f \in \dot X^{q_3^\prime}(\Om)^\ast \cap \dot X^{r^\prime}(\Omega)^\ast$ we conclude that $\bm{v} \in \dot X^r(\Omega)$.   \qed
\par
\medskip
%\noindent
Based on the bilinear forms $a_X(\cdot, \cdot)$ in (\ref{eqn:4.5}) and $a_V(\cdot, \cdot)$ in (\ref{eqn:4.6}), 
we define now $S_r:\dot X^r(\Omega) \to \dot X^{r^\prime}(\Omega)^\ast$ and $T_r:\dot V^r(\Omega) \to \dot V^{r^\prime}(\Omega)^\ast$ by 
\begin{eqnarray}
\langle S_r\bm{v}, \bm{\Phi}\rangle &\equiv& a_X(\bm{v}, \bm{\Phi}) 
\quad
\mbox{for $\bm{v} \in \dot X^r(\Omega)$ and $\bm{\Phi} \in \dot X^{r^\prime}(\Omega)$}, \label{eqn:4.17}\\
\langle T_r\bm{w}, \bm{\Psi}\rangle &\equiv& a_V(\bm{w}, \bm{\Psi}) 
\quad
\mbox{for $\bm{w} \in \dot V^r(\Omega)$ and $\bm{\Psi} \in \dot V^{r^\prime}(\Omega)$}. 
\label{eqn:4.18}
\end{eqnarray}
Here the bracket $\langle \cdot, \cdot\rangle$ denotes the duality pairing of $\dot X^{r^\prime}(\Omega)^\ast$ and $\dot X^{r^\prime}(\Omega)$ in (\ref{eqn:4.17}), 
and that of $\dot V^{r^\prime}(\Omega)^\ast$ and $\dot V^{r^\prime}(\Omega)$ in (\ref{eqn:4.18}), respectively.  
Obviously, $S_r \in \cL(\dot X^r(\Omega), \dot X^{r^\prime}(\Omega)^\ast)$ and $T_r\in \cL(\dot V^r(\Omega), \dot V^{r^\prime}(\Omega)^\ast)$. 
Here $\cL(Y, Z)$ denotes the set of bounded linear operators from $Y$ to $Z$.  
The following proposition gives estimates of $\nabla \bm{v}$ 
by means of $\rot \bm{v}$ and $\dive \bm{v}$ in $L^r(\Omega)$.  
\par
\begin{proposition}\label{pr:4.3} 
Let $\Omega$ be as in the Assumption and  let $1 < r < \infty$. 
There exists  a constant $C = C(\Omega, r, R)$ such that 
\begin{eqnarray}
\|\nabla \bm{v}\|_{L^r} &\le& 
C(\|S_r\bm{v}\|_{\dot X^{r^\prime}(\Omega)^\ast} + \|\bm{v}\|_{L^r(D)})
\quad\mbox{for all $\bm{v} \in \dot X^r(\Omega)$}, \label{eqn:4.19} \\
\|\nabla \bm{w}\|_{L^r} &\le& 
C(\|T_r\bm{w}\|_{\dot V^{r^\prime}(\Omega)^\ast} + \|\bm{w}\|_{L^r(D)})
\quad\mbox{for all $\bm{w} \in \dot V^r(\Omega)$}. \label{eqn:4.20}
\end{eqnarray}
\end{proposition}

\noindent
{\it Proof.} We first prove (\ref{eqn:4.19}).  For $\bm{v} \in \dot X^r(\Omega)$ set $S_r\bm{v} =f$. By (\ref{eqn:4.17}) 
\begin{equation}\label{eqn:4.21}
a_X(\bm{v}, \bm{\Phi}) =  (\rot \bm{v}, \rot \bm{\Phi}) + (\dive \bm{v}, \dive \bm{\Phi}) 
 =\langle f, \bm{\Phi}\rangle
\quad
\mbox{for all $\bm{\Phi} \in \dot X^{r^\prime}(\Omega)$}.  
\end{equation}
Decomposing  $\bm{v}$ as $\bm{v} = \bm{v}_1 + \bm{v}_2$ as in (\ref{eqn:4.9}), we obtain by  (\ref{eqn:4.10}) and integration by parts  
\begin{equation}
\int_{\Omega_i}(\rot \bm{v}_i\cdot\rot \bm{\varphi}_i + \dive \bm{v}_i\, \dive \bm{\varphi}_i)dx 
= \langle f, \zeta_i\bm{\varphi}_i\rangle 
+ \int_D(2\nabla \zeta_i\cdot\nabla \bm{\varphi}_i + \Delta \zeta_i\bm{\varphi}_i)\cdot \bm{v} dx, 
\quad i=1, 2 \label{eqn:4.22}
\end{equation}
for all $\bm{\varphi}_1 \in {\mathcal X}(D)$ and for all $\bm{\varphi}_2 \in C^\infty_0(\re^3)$. 
Recall that $\Omega_1 = D$ and $\Omega_2 = \re^3$.   By (\ref{eqn:3.8}) 
\begin{eqnarray}
\left|\int_D(2\nabla \zeta_1\cdot\nabla \bm{\varphi}_1 + \Delta \zeta_1\bm{\varphi}_1)\cdot \bm{v} dx\right|
&\le& C(\|\bm{v}\|_{L^r(D)}(\|\nabla \bm{\varphi}_1\|_{L^{r^\prime}(D)} 
+ \|\bm{\varphi}_1\|_{L^{r^\prime}(D)})
\nonumber \\
&\le& C\|\bm{v}\|_{L^r(D)}\|\nabla \bm{\varphi}_1\|_{L^{r^\prime}(D)} \label{eqn:4.23}
\end{eqnarray}
for all $\bm{\varphi}_1 \in {\mathcal X}(D)$.   It follows from (\ref{eqn:4.22}), (\ref{eqn:4.11}), (\ref{eqn:4.23}) and 
Proposition \ref{pr:3.5} that $\nabla \bm{v}_1 \in L^r(D)$ as well as  
\begin{eqnarray}
\|\nabla \bm{v}_1\|_{L^r(D)} 
&\le& 
C\sup_{ \bm{\varphi}_1 \in {\mathcal X}(D) }
\frac{ |(\rot \bm{v}_1, \rot \bm{\varphi}_1)_D + (\dive \bm{v}_1, \dive \bm{\varphi}_1)_D| }
{ \|\nabla \bm{\varphi}_1\|_{ L^{r^\prime}(D)} }+ C\sum_{j=1}^N|(\bm{v}_1, \bm{h}_j)_D| \nonumber \\
&\le& C(\|f\|_{\dot X^{r^\prime}(\Om)^\ast} + \|\bm{v}\|_{L^r(D)}). \label{eqn:4.24}
\end{eqnarray} 
Next, in order to treat  $\bm{v}_2 \in \dot H^{1, r}(\re^3)$, consider the functional $g_{\bm{v}}$ defined by 
$$
g_{\bm{v}}: \dot X^{r^\prime}(\Omega) \ni \bm{\Phi} \mapsto 
\int_D(2\nabla \zeta_2\cdot\nabla \bm{\Phi} + \Delta \zeta_2\bm{\Phi})\cdot \bm{v} dx \in \re.
$$
By Proposition \ref{pr:3.8}
\begin{eqnarray*}
\left|\int_D(2\nabla \zeta_2\cdot\nabla \bm{\Phi} + \Delta \zeta_2\bm{\Phi})\cdot \bm{v} dx\right| 
&\le &C\|\bm{v}\|_{L^r(D)}(\|\nabla\bm{\Phi}\|_{L^{r^\prime}(D)} + \|\bm{\Phi}\|_{L^{r^\prime}(D)}) \\
&\le &C\|\bm{v}\|_{L^r(D)}\|\nabla \bm{\Phi}\|_{L^{r^\prime}}
\end{eqnarray*}
for all $\bm{\Phi} \in \dot X^{r^\prime}(\Om)$ and for some  $C = C(\Omega, r, R)$. This means  
\begin{equation}\label{eqn:4.25}
g_{\bm{v}} \in \dot X^{r^\prime}(\Om)^\ast 
\quad 
\mbox{and  $\|g_{\bm{v}}\|_{\dot X^{r^\prime}(\Om)^\ast} \le C\|\bm{v}\|_{L^r(D)}$.   }
\end{equation}
Taking $\zeta \in C^\infty(\re^3)$ in such a way that $\zeta(x)= 1$ for $|x| \ge R+1$ and that $\zeta(x) = 0$ for $|x| \le R + 1/2$, 
we obtain by (\ref{eqn:4.25}) and Proposition \ref{pr:3.12}  
\begin{eqnarray*}
\left|\int_D(2\nabla \zeta_2\cdot\nabla \bm{\varphi}_2 + \Delta \zeta_2\bm{\varphi}_2)\cdot \bm{v} dx\right|
&=&|\langle g_{\bm{v}}, \zeta\bm{\varphi}_2\rangle| 
= |\langle \tilde{g_{\bm{v}}}, \bm{\varphi}_2\rangle_{\re^3}| 
\le  \|\tilde {g_{\bm{v}}}\|_{\widehat H^{1, r^\prime}(\re^3)^\ast}
\|\nabla \bm{\varphi}_2\|_{L^{r^\prime}(\re^3)} \\
&\le&  
C\|g_{\bm{v}}\|_{\dot X^{r^\prime}(\Om)^\ast}\|\nabla \bm{\varphi}_2\|_{L^{r^\prime}(\re^3)} 
\le  C\|\bm{v}\|_{L^r(D)}\|\nabla \bm{\varphi}_2\|_{L^{r^\prime}(\re^3)}
\end{eqnarray*}
for all $\bm{\varphi}_2\in C^{\infty}_0(\re^3)$. 
This implies together with  (\ref{eqn:4.12}), (\ref{eqn:4.22}) and Proposition \ref{pr:3.6} that 
$\bm{v}_2 \in \dot H^{1, r}(\re^3)$ as well as  
\begin{eqnarray}
\|\nabla \bm{v}_2\|_{L^r(\re^3)} 
&\le& C \sup_{\bm{\varphi}_2 \in C^\infty_0(\re^3)}
\frac{|(\rot \bm{v}_2, \rot \bm{\varphi}_2)_{\re^3} + (\dive \bm{v}_2, \dive \bm{\varphi}_2)_{\re^3}|}
{\|\nabla \bm{\varphi}_2\|_{L^{r^\prime}(\re^3)}}  \nonumber \\
&\le& 
C(\|f\|_{\dot X^{r^\prime}(\Omega)^\ast} + \|\bm{v}\|_{L^r(D)}).   \label{eqn:4.26}
\end{eqnarray}
Finally, the desired estimate (\ref{eqn:4.19}) is a consequence of (\ref{eqn:4.9}), (\ref{eqn:4.24}) and (\ref{eqn:4.26}). The proof of (\ref{eqn:4.20}) is similar to that of (\ref{eqn:4.19}) and hence omitted.    \qed 
\par
\medskip
%\noindent
The following is an immediate consequence of Proposition \ref{pr:4.3}.
\par
%\medskip 
\begin{proposition}\label{pr:4.4} 
Let $\Omega$ be as in the Assumption and let $1<r<\infty$.  
\par
{\rm (i)} The kernel $\mbox{\rm Ker}(S_r)$ is a finite dimensional subspace in $\dot X^r(\Omega)$ and the range $\mbox{\rm R}(S_r)$ is a closed subspace in $\dot X^{r^\prime}(\Omega)^\ast$.  
\par
%\noindent
{\rm (ii)} The kernel $\mbox{\rm Ker}(T_r)$ is a finite dimensional subspace in $\dot V^r(\Omega)$ and the range $\mbox{\rm R}(T_r)$ is a closed subspace in $\dot V^{r^\prime}(\Omega)^\ast$.  
\end{proposition}
%\medskip\noindent
{\it Proof.} By Proposition \ref{pr:3.8} and Rellich's theorem, both $\dot X^r(\Omega)$ and $\dot V^r(\Omega)$ are compactly embedded into $L^r(D)$.  
Hence, the assertion follows from (\ref{eqn:4.19}), (\ref{eqn:4.20}) and \cite[Lemma 3]{Pe}. 
\qed
\par
\bigskip
Concerning the kernels $\mbox{Ker}(S_r)$ and  $\mbox{Ker}(T_r)$, 
we have the following  precise characterization.  
\begin{proposition}\label{pr:4.5} Let $\Omega$ be as in the Assumption and let $1<r<\infty$. 
Then $\mbox{\rm Ker}(S_r) = \dot X_{\tiny{\mbox{\rm har}}}(\Omega)$ and 
$\mbox{\rm Ker}(T_r) = \dot V_{\tiny{\mbox{\rm har}}}(\Omega)$.
\end{proposition}

\medskip\noindent
{\it Proof}. Obviously  $\dot X_{\tiny{\mbox{har}}}(\Omega) \subset \mbox{Ker}(S_r)$ and hence we only prove  that $\mbox{Ker}(S_r) \subset \dot X_{\tiny{\mbox{har}}}(\Omega)$.  
Let $\bm{v} \in \mbox{Ker}(S_r)$. By (\ref{eqn:4.17}) 
\begin{equation}\label{eqn:4.27}
a_X(\bm{v}, \bm{\varphi})= (\rot \bm{v}, \rot \bm{\varphi}) + (\dive \bm{v}, \dive \bm{\varphi})=0
\quad\mbox{for all $\bm{\varphi} \in {\mathcal X}(\Omega)$}.
\end{equation}
Hence, it follows from Proposition \ref{pr:4.2} that $\bm{v} \in \dot X^2(\Omega)$.  
Since ${\mathcal X}(\Omega)$ is dense in $\widehat X^2(\Omega)$, 
we see that (\ref{eqn:4.27}) remains valid for all $\bm{\varphi} \in \widehat X^2(\Omega)$. 
Moreover, we may choose $\bm{\varphi} \in  \dot X^2(\Omega)$ in (\ref{eqn:4.27}). 
Indeed, by Proposition \ref{pr:3.7} (i), every $\bm{\varphi} \in \dot X^2(\Omega)$ 
may be  expressed as 
$$
\bm{\varphi} = \bm{\psi} + \sum_{j=1}^3\lambda_j\bm{h}_j 
\quad
\mbox{with some $\bm{\psi} \in \widehat X^2(\Omega)$ 
and $\lambda_1, \lambda_2, \lambda_3 \in \re$}.  
$$ 
Since $\dive \bm{h}_j =0$, $\rot \bm{h}_j =\bm{0}$, $j=1, 2, 3$ and since (\ref{eqn:4.27}) holds for $\bm{\varphi} = \bm{\psi} \in \widehat X^2(\Omega)$, we obtain 
$$
a_X(\bm{v}, \bm{\varphi}) = a_X(\bm{v}, \bm{\psi}) + \sum_{j=1}^3\lambda_ja_X(\bm{v}, \bm{h}_j) =0.  
$$
Choosing $\bm{\varphi} =\bm{v}$ yields $\dive \bm{v}=0$ and $\rot \bm{v}=\bm{0}$.  This means that $\bm{v} \in \dot X_{\tiny{\mbox{har}}(\Omega)}$.  
Since $\bm{v} \in \mbox{Ker}(S_r)$ is arbitrary, we obtain $\mbox{Ker}(S_r) \subset \dot X_{\tiny{\mbox{har}}}(\Omega)$.
Taking into account assertion  (\ref{eqn:3.28}), the proof of the second assertion parallels the first one.  \qed
\par
\medskip
%\noindent
The following lemma plays an essential role for the proof of Theorem \ref{thm:4.1}. We recall that the spaces $\widetilde X^r(\Omega)$ in $\dot X^r(\Omega)$ and $\widetilde V^r(\Omega)$ 
in $\dot V^r(\Omega)$ were introduced in (\ref{eqn:3.48}).   
%\par
%
%\medskip
%\noindent
\begin{lemma}\label{lem:4.6}  Let $\Omega$ be as in the Assumption and let $1<r < \infty$. \\
{\rm (i)} For given $f \in \dot X^{r^\prime}(\Omega)^\ast$, 
there exists $\bm{v} \in \dot X^r(\Omega)$ 
such that
\begin{equation}\label{eqn:4.28} 
S_r\bm{v} = f
\end{equation}
if and only if $f$ satisfies that 
\begin{equation}\label{eqn:4.29} 
\langle f, \bm{\varphi}\rangle =0
\quad
\mbox{for all $\bm{\varphi} \in \dot X_{\tiny{\mbox{\rm har}}}(\Omega) $}, 
\end{equation}
wehre $\langle\cdot, \cdot\rangle$  denotes the duality pairing of  $\dot X^{r^\prime}(\Omega)^\ast$ and $\dot X^{r^\prime}(\Omega)$.  
\par
%\noindent
{\rm (ii)} For every $f \in \dot X^{r^\prime}(\Omega)^\ast$ satisfying the condition 
(\ref{eqn:4.29}) there exists a unique $\tilde {\bm{v}}\in \widetilde X^r(\Omega)$ satisfying 
\begin{equation}\label{eqn:4.30} 
S_r\tilde {\bm{v}} = f.
\end{equation}
Moreover, there exists a constant $C=C(\Omega, r, R)$ such that 
\begin{equation}\label{eqn:4.31} 
\|\nabla \tilde {\bm{v}}\|_{L^r} \le C\|f\|_{\dot X^{r^\prime}(\Omega)^\ast}.
\end{equation}
\par
%\noindent
{\rm (iii)} For given $g \in \dot V^{r^\prime}(\Omega)^\ast$ there exists 
$\bm{w} \in \dot V^r(\Omega)$ such that
\begin{equation}\label{eqn:4.32} 
T_r\bm{w} = g
\end{equation}
if and only if $g$ satisfies that 
\begin{equation}\label{eqn:4.33} 
\langle g, \bm{\psi}\rangle =0
\quad
\mbox{for all $\bm{\psi} \in \dot V_{\tiny{\mbox{\rm har}}}(\Omega) $}, 
\end{equation}
wehre $\langle\cdot, \cdot\rangle$  denotes the duality paring of  $\dot V^{r^\prime}(\Omega)^\ast$ and $\dot V^{r^\prime}(\Omega)$.  
\par
\noindent
{\rm (iv)} For every $g \in \dot V^{r^\prime}(\Omega)^\ast$ satisfying the condition 
(\ref{eqn:4.33}) there exists a unique $\tilde {\bm{w}}\in \widetilde V^r(\Omega)$ such that 
\begin{equation}\label{eqn:4.34} 
T_r\tilde {\bm{w}} = g.
\end{equation}
Moreover, there exists a constant $C=C(\Omega, r, R)$ such that
\begin{equation}\label{eqn:4.35} 
\|\nabla \tilde {\bm{w}}\|_{L^r} \le C\|g\|_{\dot V^{r^\prime}(\Omega)^\ast}.
\end{equation}
\end{lemma}
%\medskip\noindent
{\it Proof}. (i) By Propositions \ref{pr:4.4} (i) and  \ref{pr:4.5} (i),  and the closed range theorem 
\begin{equation}
\mbox{R}(S_r) = ^{\perp}\mbox{Ker}(S_r^\ast)=^{\perp}\mbox{Ker}(S_{r^\prime}) 
=^{\perp}\dot X_{\tiny{\mbox{\rm har}}}(\Omega) 
= \{f \in \dot X^{r^\prime}(\Omega)^\ast; \langle f, \bm{\varphi}\rangle =0, \forall \bm{\varphi} \in 
\dot X_{\tiny{\mbox{\rm har}}}(\Omega)\}, \label{eqn:4.36}
\end{equation}
which implies assertion (i). 
\par
%\noindent
(ii) Suppose that $f\in \dot X^{r^\prime}(\Omega)^\ast$ satisfies the condition (\ref{eqn:4.29}). By (i), there exists a solution $\bm{v}\in \dot X^r(\Omega)$ of (\ref{eqn:4.28}).  
By (\ref{eqn:3.48}) $\bm{v}$ may be decomposed as $\bm{v}=\bm{v}_0 + \tilde {\bm{v}}$, 
where $\bm{v}_0 \in \dot X_{\tiny{\mbox{\rm har}}}(\Omega)$ and 
$\tilde {\bm{v}} \in \widetilde X^r(\Omega)$.  
Proposition \ref{pr:4.5} implies  $S_r\bm{v}_0 =0$ and hence $S_r\tilde {\bm{v}} =f$. This  yields (\ref{eqn:4.30}). For the proof of the estimate (\ref{eqn:4.31}), it suffices to show that 
\begin{equation}\label{eqn:4.37}
\|\nabla \tilde {\bm{v}}\|_{L^r}\le C\|S_r\tilde {\bm{v}}\|_{\dot X^{r^\prime}(\Omega)^\ast}
\quad
\mbox{for all $\tilde {\bm{v}} \in \tilde X^r(\Omega)$}    
\end{equation} 
with some $C=C(\Omega, r, R)$.  
Assume that (\ref{eqn:4.37}) fails. 
Then there is a sequence $\{\tilde {\bm{v}}_m\}_{m=1}^\infty \subset \widetilde X^r(\Omega)$ such that 
$\|\nabla \tilde {\bm{v}}_m\|_{L^r}\equiv 1$ for all $m\in \N$ and $\lim_{m\to\infty}\|S_r\tilde {\bm{v}}_m\|_{\dot X^{r^\prime}(\Omega)^\ast} =0$. 
Since $\widetilde X^r(\Omega)$ is a closed subspace of $\dot X^r(\Omega)$, 
there is a subsequence of $\{\tilde {\bm{v}}_m\}_{m=1}^\infty$, 
denoted again by $\{\tilde {\bm{v}}_m\}_{m=1}^\infty$ itself,  
and a function $\tilde {\bm{v}} \in \widetilde X^r(\Omega)$ such that 
$$
\nabla \tilde {\bm{v}}_m \rightharpoonup \nabla \tilde {\bm{v}}
\quad
\mbox{weakly in $L^r(\Omega)$ as $m\to\infty$}.   
$$    
Since the injection $\dot X^r(\Omega) \subset L^r(D)$ is compact, the above weak convergence yields 
$$
\lim_{m\to\infty}\|\tilde {\bm{v}}_m - \tilde {\bm{v}}\|_{L^r(D)} =0.  
$$
Hence, by (\ref{eqn:4.19}) $\{\nabla \tilde {\bm{v}}_m\}_{m=1}^\infty$ is a Cauchy sequence in $L^r(\Omega)$, which yields  
$$
\nabla \tilde {\bm{v}}_m \to \nabla \tilde {\bm{v}}
\quad
\mbox{strongly in $L^r(\Omega)$ as $m\to\infty$}.   
$$ 
Since $S_r \in \cL(\dot X^r(\Omega), \dot X^{r^\prime}(\Omega)^\ast)$, we obtain  $S_r\tilde {\bm{v}} =0$,  which implies $\tilde {\bm{v}} \in \mbox{Ker}(S_r)=\dot X_{\tiny{\mbox{\rm har}}}(\Omega)$.
Finally, since $\dot X_{\tiny{\mbox{\rm har}}}(\Omega)\cap \tilde X^r(\Omega)= \{\bm{0}\}$, we obtain $\tilde {\bm{v}}=0$, which contradicts $\|\nabla \tilde {\bm{v}}\|_{L^r}=1$.  
\par
%\noindent
(iii) and (iv) Making use of Propositions \ref{pr:4.4} (ii) and \ref{pr:4.5}, 
the proof of assertions (iii) and (iv) parallels the one of the above (i) and (ii), respectively. 
Details are hence omitted.  \qed 
\par
\medskip
%\noindent
The following lemma is a consequence of Lemma \ref{lem:4.6}
\par
%\smallskip\noindent
\begin{lemma}\label{lem:4.7} Let $\Omega$ be as in the Assumption and let  $1 < r <\infty$.
\par
%\noindent
{\rm (i)} For every $\bm{u} \in L^r(\Omega)$ there exits a unique 
$\bm{v} \in \widetilde  X^r(\Omega)$ such that 
\begin{equation}\label{eqn:4.38}
a_X(\bm{v}, \bm{\Phi})=(\bm{u}, \rot \bm{\Phi})
\quad
\mbox{for all $\bm{\Phi} \in \dot X^{r^\prime}(\Omega)$}.  
\end{equation}
Moreover,  $\bm{v}$ satisfies the estimate 
\begin{equation}\label{eqn:4.39}
\|\nabla \bm{v}\|_{L^r}\le C\|\bm{u}\|_{L^r}, 
\end{equation}
for some constant  $C=C(\Omega,r,R)$.  
\par
%\noindent
{\rm (ii)} For every $\bm{u} \in L^r(\Omega)$ there exits a unique 
$\bm{w} \in \widetilde  V^r(\Omega)$ such that 
\begin{equation}\label{eqn:4.40}
a_V(\bm{w}, \bm{\Psi})=(\bm{u}, \rot \bm{\Psi})
\quad
\mbox{for all $\bm{\Psi} \in \dot V^{r^\prime}(\Omega)$}.  
\end{equation}
Moreover,  $\bm{w}$ satisfies the estimate 
\begin{equation}\label{eqn:4.41}
\|\nabla \bm{w}\|_{L^r}\le C\|\bm{u}\|_{L^r}, 
\end{equation}
for some  $C=C(\Omega, r, R)$. 
\end{lemma}
%\noindent
{\it Proof.} (i) Given $\bm{u} \in L^r(\Omega)$, we define the functional $f_{\bm{u}}$ by 
\begin{equation}\label{eqn:4.42}
\langle f_{\bm{u}}, \bm{\Phi}\rangle \equiv (\bm{u}, \rot \bm{\Phi})
\quad
\mbox{for $\bm{\Phi} \in \dot X^{r^\prime}(\Omega)$}.  
\end{equation}
Obviously, $f_{\bm{u}}\in \dot X^{r^\prime}(\Omega)^\ast$ and  
\begin{equation}\label{eqn:4.43}
\|f_{\bm{u}}\|_{\dot X^{r^\prime}(\Omega)^\ast} \le \|\bm{u}\|_{L^r}.  
\end{equation}
By (\ref{eqn:4.42}), $\langle f_{\bm{u}}, \bm{\varphi}\rangle =0$ for all $\bm{\varphi} \in \dot X_{\tiny{\mbox{\rm har}}}(\Omega)$. This means that condition (\ref{eqn:4.29}) is fulfilled for $f_{\bm{u}}$. 
Hence, it follows from Lemma \ref{lem:4.6} and (\ref{eqn:4.43}) that there exists a unique $\bm{v} \in \tilde X^r(\Omega)$ such that 
\begin{equation}\label{eqn:4.44}
S_r\bm{v} = f_{\bm{u}} 
\end{equation}
satisfying 
\begin{equation}\label{eqn:4.45}
\|\nabla \bm{v}\|_{L^r} \le C\|f\|_{ \dot X^{r^\prime}(\Omega)^\ast} \le C\|\bm{u}\|_{L^r}.
\end{equation} 
By (\ref{eqn:4.42}), (\ref{eqn:4.44}) and  (\ref{eqn:4.45}) we see that $\bm{v}$ solves (\ref{eqn:4.38}) and satisfies  (\ref{eqn:4.39}).  
\par
%\noindent
(ii) By Proposition \ref{pr:4.5} and Lemma \ref{lem:4.6}(iv), 
the above argument of (i) is also applicable  to the proof of assertion ii). 
Again we omit the details. \qed  
\subsection{Proof of Theorem \ref{thm:4.1}} \mbox{}
\par
We base the  proof of Theorem \ref{thm:4.1} on Lemma \ref{lem:4.7}. Indeed, we only need to show that the two solutions $\bm{v} \in \widetilde X^r(\Omega)$ of (\ref{eqn:4.38}) 
and $\bm{w} \in \widetilde V^r(\Omega)$ of (\ref{eqn:4.40}) necessarily satisfy $\dive \bm{v} = \dive \bm{w} =0$. To this end, we investigate the bilinear forms $a_X(\cdot, \cdot)$ on 
$\dot X^2_\sigma(\Omega)\times \dot X^2_\sigma(\Omega)$ and $a_V(\cdot, \cdot)$ on $\dot V^2_\sigma(\Omega)\times \dot V^2_\sigma(\Omega)$,  respectively.        
Similarly to (\ref{eqn:2.3}), we introduce the spaces $\widehat X^r_\sg(\Omega)$ and 
$\widehat V^r_\sg(\Omega)$ defined by 
\begin{equation}\label{eqn:4.46}
\widehat X^r_{\sg}(\Omega)\equiv \{\bm{u} \in \hat X^r(\Om); \dive \bm{u} = 0\}, 
\quad
\widehat V^r_{\sg}(\Omega)\equiv \{\bm{u} \in \hat V^r(\Om); \dive \bm{u} = 0\}.  
\end{equation}
\begin{proposition}\label{pr:4.8} Let $\Omega$ be as in the Assumption. 
\par
%\noindent
{\rm (i)} For every $\bm{\phi} \in \dot X^2(\Omega)$ there exist 
$\bm{\phi}_\sg \in \widehat X^2_{\sigma}(\Omega)$ and $\bm{\phi}_c \in \dot X^2(\Omega)$ 
with $\rot \bm{\phi}_c =\bm{0}$ such that 
$\bm{\phi}$ can be decomposed as 
\begin{equation}\label{eqn:4.47}
\bm{\phi}= \bm{\phi}_\sg + \bm{\phi}_c.  
\end{equation} 
\par
{\rm (ii)} For every $\bm{\psi} \in \dot V^2(\Omega)$ 
there exist $\bm{\psi}_\sg \in \widehat V^2_\sg(\Omega)$ and $\bm{\psi}_c \in \dot V^2(\Omega)$ 
with $\rot \bm{\psi}_c =\bm{0}$ such that 
$\bm{\psi}$ can be  decomposed as 
\begin{equation}\label{eqn:4.48}
\bm{\psi}= \bm{\psi}_\sg + \bm{\psi}_c.  
\end{equation}
\end{proposition}  
%\noindent
{\it Proof.} (i) We consider first the case when $\bm{\phi} \in \widehat X^2(\Omega)$.   
Then there exists $p \in \dot H^{1, 6}(\Omega)$ with $D^2 p \in L^2(\Omega)$ such that 
\begin{equation}\label{eqn:4.49}
\left\{
\begin{array}{ll}
& \Delta p = \dive \bm{\phi} 
\quad
\mbox{in $\Omega$}, \\
& \dfrac{\pt p}{\pt\bm{\nu}} =\bm{\phi}\cdot\bm{\nu} =0 
\quad
\mbox{on $\pt\Omega$}.
\end{array}
\right.
\end{equation}
For the time being, assume the existence of such $p$. Then, setting $\bm{\phi}_{\sg} \equiv \bm{\phi} - \nabla p$ and $\bm{\phi}_c \equiv \nabla p$, 
we obtain by Proposition \ref{pr:3.7} that $\bm{\phi}_\sg \in \widehat X_{\sg}^2(\Omega)$ and 
$\bm{\phi}_c \in \dot X^2(\Omega)$ 
with $\rot \bm{\phi}_c =\bm{0}$. Hence,  the desired decomposition (\ref{eqn:4.47}) follows.  
For general $\bm{\phi} \in \dot X^2(\Omega)$, (\ref{eqn:3.27}) implies 
$$
\bm{\phi}= \hat {\bm{\phi}} + \sum_{j=1}^3\lambda_j\bm{h}_j,
$$  
where $\hat {\bm{\phi}} \in \widehat X^2(\Omega)$ and $\lambda_1, \lambda_2, \lambda_3 \in \re$.  
Since $\hat {\bm{\phi}} \in \widehat X^2(\Omega)$, the above observation yields that  
$\hat {\bm{\phi}}$ may be expressed as 
$$
\hat {\bm{\phi}} = \hat {\bm{\phi}}_\sg + \hat {\bm{\phi}}_c, 
$$
where $\hat {\bm{\phi}}_\sg \in \widehat X^2_{\sg}(\Omega)$ and $
\hat {\bm{\phi}}_c \in \dot X^2(\Omega)$ with $\rot \hat{\bm{\phi}} = \bm{0}$.  
Now defining $\bm{\phi}_\sg \equiv \hat {\bm{\phi}}_\sg$ and 
$\bm{\phi}_c \equiv \hat {\bm{\phi}}_c + \sum_{j=1}^3\lambda_j\bm{h}_j$, 
we see  that $\bm{\phi}_\sg \in \widehat X_\sg^2(\Omega)$ and 
$\bm{\phi}_c \in \dot X^2(\Omega)$ with $\rot \bm{\phi}_c = \bm{0}$, 
and the desired decomposition (\ref{eqn:4.47}) follows again.

It remains to show the existence of a solution $p$ of (\ref{eqn:4.49}) for 
$\bm{\phi}\in \widehat X^2(\Omega)$ within the class $p \in \dot H^{1, 6}(\Omega)$ with 
$D^2p \in L^2(\Omega)$.  
Since ${\mathcal X}(\Omega)$ is dense in $\widehat X^2(\Omega)$, 
there is a sequence $\{\bm{\phi}_m\}_{m=1}^\infty \subset {\mathcal X}(\Omega)$ such that 
\begin{equation}\label{eqn:4.50}
\nabla \bm{\phi}_m \to \nabla \bm{\phi}
\quad
\mbox{in $L^2(\Omega)$}, 
\quad
\bm{\phi}_m \to \bm{\phi}
\quad
\mbox{in $L^6(\Omega)$ as $m\to\infty$}.  
\end{equation}
The solvability of the weak Neumann problem (see, e.g.,  \cite[Theorem 1.3]{SiSo2} and \cite[Theorem 4.1]{SiSo3})  
implies that for each $m\in \N$ there exists a function $p_m \in \dot H^{1, r}(\Omega)$ for all $1 < r < \infty$ satisfying 
\begin{equation}\label{eqn:4.51}
(\nabla p_m, \nabla \varphi) = (\bm{\phi}_m, \nabla \varphi)
\quad
\mbox{for all $\varphi \in C^\infty_0(\bar\Omega)$}.   
\end{equation}   
Hence, by (\ref{eqn:4.50}), (\ref{eqn:4.51}) and the variational inequality we obtain 
\begin{eqnarray*}
\|\nabla p_m - \nabla p_l\|_{L^6} 
&\le& C\sup_{\varphi \in C^\infty_0(\bar\Omega)} 
\frac{|(\nabla p_m - \nabla p_l, \nabla\varphi)|}{\|\nabla \varphi\|_{L^{\frac 65}} } \\
&\le& C\sup_{\varphi \in C^\infty_0(\bar\Omega)} 
\frac{|( \bm{\phi}_m - \bm{\phi}_l, \nabla\varphi)|}{\|\nabla \varphi\|_{L^{\frac 65}} } \\
&\le& C\|\bm{\phi_m} - \bm{\phi}_l\|_{L^6} 
\to 0 
\end{eqnarray*}
as $m, l \to \infty$.  Hence, there exists $p \in \dot H^{1, 6}(\Omega)$ such that 
\begin{equation}\label{eqn:4.52}
\nabla p_m \to \nabla p 
\quad
\mbox{in $L^6(\Omega)$ as $m\to\infty$}.  
\end{equation}
Now, letting $m\to\infty$ in both sides of (\ref{eqn:4.51}), we obtain by (\ref{eqn:4.50}) and (\ref{eqn:4.52}) that 
$$
(\nabla p,\nabla \varphi) = (\bm{\phi}, \nabla\varphi)
\quad
\mbox{for all $\varphi \in C^\infty_0(\bar\Omega)$}.
$$ 
This  means that $p$ is the unique weak solution of (\ref{eqn:4.49}).  
\par
We next show that $D^2p \in L^2(\Omega)$.  Let $\eta\in C^\infty_0(\re^3)$ be the cut-off function in the proof of Lemma \ref{lem:3.9}.  
Similarly to (\ref{eqn:3.36}), we set 
$$
p^{(2)}_m(x)\equiv (1-\eta(x))(p_m(x) - p^D_m)
\quad
\mbox{with $\dis{p^D_m \equiv \frac{1}{|D|}\int_Dp_m(y)dy}$, $m \in \N$}.  
$$
Notice that $p^{(2)}_m \in \dot H^{1, r}(\re^3)$ for all $1 < r < \infty$ and that 
$p^{(2)}_m(x) = p_m(x) - p^D_m$ for $|x|\ge R+2$. By (\ref{eqn:4.51}) 
$$
\Delta p_m^{(2)}= (1-\eta)\dive \bm{\phi}_m -2\nabla \eta\cdot\nabla p_m -\Delta\eta(p_m - p^D_m)
\equiv f_m, 
\quad m \in \N. 
$$
Set $\tilde p_m^{(2)}\equiv \Gamma\ast f_m$, where $\Gamma(x) = -(4\pi|x|)^{-1}$ denotes the Newton kernel. Since $\nabla p^{(2)}_m - \nabla \tilde p_{m}^{(2)} \in L^{q}(\re^3)$ for all $3/2 < q< \infty $ 
and since $ p^{(2)}_m - \tilde p_{m}^{(2)}$ is harmonic in $\re^3$, it follows from 
Liouville's theorem that for each $m\in\N$ there exists a constant $c_m \in \re$ such that $p^{(2)}_m(x) - \tilde p_{m}^{(2)}(x)=c_m$ for all $x \in \re^3$.  
The Calder\'on-Zygmund estimate and Poincar\'e's inequality on $D$ yield 
\begin{eqnarray}
&& \|D^2p^{(2)}_m - D^2p_{l}^{(2)}\|_{L^2(\re^3)} \nonumber \\
&=& \|D^2\tilde p^{(2)}_m - D^2\tilde p_{l}^{(2)}\|_{L^2(\re^3)} \nonumber \\
&\le&  C\|\Delta(\tilde p^{(2)}_m - \tilde p_{l}^{(2)})\|_{L^2(\re^3)} \nonumber \\
&\le&  C\|f_m - f_l\|_{L^2(\re^3)} \nonumber \\
&\le& C(\|\dive (\phi_m - \phi_l)\|_{L^2} + \|\nabla (p_m - p_l)\|_{L^2(D)} 
+ \|p_m-p_l-(p_m^D - p_l^D)\|_{L^2(D)}) \nonumber \\
&\le &C(\|\nabla (\phi_m - \phi_l)\|_{L^2} + \|\nabla(p_m - p_l)\|_{L^2(D)}) \nonumber \\
&\le & C(\|\nabla (\phi_m - \phi_l)\|_{L^2} + \|\nabla (p_m - p_l)\|_{L^6})\label{eqn:4.53}
\end{eqnarray}
for all $m,  l \in \N$. It follows from  (\ref{eqn:4.50}) and (\ref{eqn:4.52})  that there exists  $g \in L^2(\re^3)$ such that $D^2p^{(2)}_m \to g$ in $L^2(\re^3)$.   
Since $p_m(x) = p^{(2)}_m(x)$ for all $|x| \ge R+2$ and all $m \in \N$, we see   by (\ref{eqn:4.52}) that $D^2p(x) = g(x)$ for all $|x|\ge R+2$. Thus   $D^2 p \in L^2(\Omega)$.  
\par
%\noindent
(ii)  We first consider the case where $\psi \in \widehat V^2(\Omega)$. Similarly to (\ref{eqn:4.49}) we  consider the Dirichlet problem 
\begin{equation}\label{eqn:4.54}
\left\{
\begin{array}{ll}
 \Delta p &= \dive \psi 
\quad
\mbox{in $\Omega$}, \\
p &= 0 
\quad
\mbox{on $\pt\Omega$}
\end{array}
\right.
\end{equation}
in the class $p \in \dot H^{1, 6}_0(\Omega)$ with $D^2 p \in L^2(\Omega)$. Setting $\psi_\sg \equiv \psi - \nabla p$ and $\psi_c \equiv \nabla p$, we see  that 
$\psi_\sg \in \widehat V^2(\Omega)$, $\psi_c \in \dot V^2(\Omega)$ 
and that the decomposition (\ref{eqn:4.48}) holds.  
For general $\psi \in \dot V^2(\Om)$, in the same way as above, 
we  make use of (\ref{eqn:3.28}) to reduce the general case to $\psi \in \hat V^2(\Omega)$.     
\par
It remains to show the existence of a solution $p \in \dot H^{1, 6}_0(\Omega)$ to (\ref{eqn:4.54}) with $D^2 p \in L^2(\Omega)$.  
Since ${\mathcal V}(\Omega)$ is dense in $\hat V^2(\Omega)$ there is a sequence 
$\{\psi_m\}_{m=1}^\infty \subset {\mathcal V}(\Omega)$ such that 
\begin{equation}\label{eqn:4.55}
\nabla \psi_m \to \nabla \psi
\quad
\mbox{in $L^2(\Omega)$}, 
\quad
\psi_m \to \psi
\quad
\mbox{in $L^6(\Omega)$ as $m\to\infty$}.  
\end{equation}
By (\ref{eqn:3.2}) and Proposition \ref{pr:3.2}(iii) we know  that $-\Delta_6$ is bijective as a mapping from  $\widetilde H^{1, 6}_0(\Omega)$ onto 
$\widehat H^{1, \frac 65}_0(\Omega)^\ast$, and hence 
for each $m\in\N$, there exists a unique $p_m \in \widetilde H^{1, 6}_0(\Omega)$ such that 
\begin{equation}\label{eqn:4.56}
(\nabla p_m, \nabla \Phi) = (\psi_m, \Phi)
\quad
\mbox{for all $\Phi \in C^\infty_0(\Omega)$}.   
\end{equation}   
Moreover,
\begin{eqnarray*}
\|\nabla p_m - \nabla p_l\|_{L^6} 
&\le& C\sup_{\Phi \in C^\infty_0(\Omega)} 
\frac{|(\nabla p_m - \nabla p_l, \nabla\Phi)|}{\|\nabla \Phi\|_{L^{\frac 65}} } 
\le  C\sup_{\Phi \in C^\infty_0(\Omega)} 
\frac{|( \psi_m - \psi_l, \nabla\Phi)|}{\|\nabla \Phi\|_{L^{\frac 65}} } 
\le C\|\psi_m - \psi_l\|_{L^6} \\
& \to& 0 
\end{eqnarray*}
as $m, l \to \infty$.  Hence, there exists $p \in \widetilde H^{1, 6}_0(\Omega)$ such that 
\begin{equation}\label{eqn:4.57}
\nabla p_m \to \nabla p 
\quad
\mbox{in $L^6(\Omega)$ as $m\to\infty$}.   
\end{equation}
Letting $m\to \infty$ on both sides of (\ref{eqn:4.56}), we obtain by (\ref{eqn:4.55}) and (\ref{eqn:4.57}) that 
\begin{equation}\label{eqn:4.58}
(\nabla p, \nabla \Phi) = (\psi, \nabla \Phi)
\quad
\mbox{for all $\Phi \in C^\infty_0(\Omega)$}.
\end{equation}
This means that $p$ is the weak solution of (\ref{eqn:4.54}).   
\par
We finally show that $D^2 p \in L^2(\Omega)$.  
The argument is similar to the one given above (i).  
Indeed, setting  
$p_m^{(2)}(x)\equiv (1-\eta(x))p_m(x)$, we see that 
\begin{equation}\label{eqn:4.59}
\|D^2p^{(2)}_m - D^2p^{(2)}_l\|_{L^2(\re^3)} \le C(\|\nabla(\psi_m - \psi_l)\|_{L^2} 
+ \|\nabla(p_m - p_l)\|_{L^6}), 
\quad m, l \in \N.  
\end{equation}
In comparison with (i), it is here unnecessary to subtract the mean value $p^D_m$ on $D$ 
since $p_m|_{\pt\Om} = 0$, so that we may make use of Poincar\'e's inequality.  
In fact,  
$\|p_m - p_l\|_{L^2(D)} \le C\|\nabla (p_m - p_l)\|_{L^2(D)}$ for all $m, l \in \N$. By (\ref{eqn:4.59}) there is $g \in L^2(\re^3)$ such that $p^{(2)}_m \to g$ in $L^2(\re^3)$ 
as $m\to\infty$. It follows from  (\ref{eqn:4.57}) that $D^2p(x) = g(x)$ for $|x| \ge R+2$. Hence, $D^2p \in L^2(\Omega)$.     \qed
\par
\bigskip
%\noindent
We recall that $\dot X_\sg^r(\Omega)$ and $\dot V_{\sg}^r(\Omega)$ are 
defined by (\ref{eqn:2.3}). 
Since $\dot X_{\tiny{\mbox{\rm har}}}(\Omega)$ is a finite dimensional subspace of 
$\dot X^2_{\sg}(\Omega)$, 
there is a closed subspace $\widetilde X_{\sg}^2(\Omega)$ of $\dot X^2_{\sg}(\Omega)$ such that 
\begin{equation}\label{eqn:4.60}
\dot X^2_{\sg}(\Omega)= \widetilde X_{\sg}^2(\Omega) \oplus \dot X_{\tiny{\mbox{\rm har}}}(\Omega), 
\quad
\mbox{(direct sum).} 
\end{equation} 
Similarly, there is a closed subspace $\widetilde V_{\sg}^2(\Omega)$ of $\dot V^2_{\sg}(\Omega)$ 
such that 
\begin{equation}\label{eqn:4.61}
\dot V^2_{\sg}(\Omega)= \widetilde V_{\sg}^2(\Omega) \oplus \dot V_{\tiny{\mbox{\rm har}}}(\Omega), 
\quad
\mbox{(direct sum).}
\end{equation} 
%\medskip
\begin{lemma}\label{lem:4.9} 
For every $\bm{u} \in L^2(\Omega)$ there exists a unique $\bm{v} \in \widetilde X_{\sg}^2(\Omega)$ 
such that 
\begin{equation}\label{eqn:4.62}
a_X(\bm{v},\bm{\varphi}) = (\bm{u}, \rot \bm{\varphi})
\quad
\mbox{for all $\bm{\varphi} \in \dot X^2(\Omega)$}.  
\end{equation}
For every $\bm{u} \in L^2(\Omega)$ there is a unique $\bm{w} \in \widetilde V_{\sg}^2(\Omega)$ 
such that 
\begin{equation}\label{eqn:4.63}
a_V(\bm{w}, \bm{\psi}) = (\bm{u}, \rot \bm{\psi})
\quad
\mbox{for all $\bm{\psi} \in \dot V^2(\Omega)$}.  
\end{equation}
\end{lemma}

\medskip
{\it Proof}. Again, we only prove the first assertion.   
We show first that for every $\bm{u} \in L^2(\Omega)$ there exists a unique 
$\bm{v} \in \widetilde X^2_\sg(\Omega)$ such that 
\begin{equation}\label{eqn:4.64}
(\rot \bm{v}, \rot \bm{\phi}) = (\bm{u}, \rot \bm{\phi}) 
\quad
\mbox{for all $\bm{\phi} \in \widetilde X^2_\sg(\Omega)$}.  
\end{equation} 
Indeed, by (\ref{eqn:3.51}), 
$(\bm{v}, \bm{\phi})_{\widetilde X^2_\sg} \equiv (\rot \bm{v}, \rot \bm{\phi})$ 
defines an inner product on the Hilbert space $\widetilde X^2_\sg(\Omega)$.   
Then, for $\bm{u} \in L^2(\Omega)$ consider the functional $f_{\bm{u}}$ on 
$\widetilde X^2_{\sigma}(\Omega)$ defined by 
$\langle f_{\bm{u}}, \phi\rangle = (\bm{u}, \rot \bm{\phi})$ 
for $\bm{\phi} \in \widetilde X^2_\sg(\Omega)$. 
Since 
$$
|\langle f_{\bm{u}}, \bm{\phi}\rangle | \le \|\bm{u}\|_{L^2}\|\rot \bm{\phi}\|_{L^2} \le \|\bm{u}\|_{L^2}\|\bm{\phi}\|_{\widetilde X^2_\sg}
$$
for all $\bm{\phi}\in \widetilde X^2_\sg(\Omega)$, 
we see  that $f_{\bm{u}} \in \widetilde X^2_\sg(\Omega)^\ast$. 
Hence, by the Riesz representation theorem, there exists a unique 
$\bm{v} \in \widetilde X^2_\sg(\Omega)$ such that 
$(\bm{v}, \bm{\phi})_{\widetilde X^2_\sg}= \langle f_{\bm{u}}, \bm{\phi}\rangle$ 
for all $\bm{\phi} \in \widetilde X^2_\sg(\Omega)$. This  implies (\ref{eqn:4.64}).   
\par
We next show that $\bm{v}\in \widetilde X^2_\sg(\Omega)$ satisfies also (\ref{eqn:4.62}). 
By Proposition \ref{pr:4.8} (i), for  $\bm{\varphi} \in \dot X^2(\Omega)$ 
there exist $\bm{\varphi}_\sg \in \widehat X^2_{\sg}(\Omega)$ and $\bm{\varphi}_c \in \dot X^2(\Omega)$ 
with $\rot \bm{\varphi}_c =\bm{0}$ such that $\bm{\varphi} = \bm{\varphi}_\sg + \bm{\varphi}_c$.  
Moreover, since $\bm{\varphi}_\sg \in \widehat X^2_{\sg}(\Omega) \subset \dot X^2_\sg(\Omega)$, (\ref{eqn:4.60})  implies $\bm{\varphi}_\sg = \tilde {\bm{\varphi}}_{\sg} + \bm{\varphi}_{\tiny{\mbox{\rm har}}}$ for 
some $\tilde {\bm{\varphi}}_{\sg} \in \widetilde X^2_{\sg}(\Omega)$ and $\bm{\varphi}_{\tiny{\mbox{\rm har}}} \in \dot X_{\tiny{\mbox{\rm har}}}$.  
Since $\rot \bm{\varphi} = \rot \bm{\varphi}_\sg = \rot \tilde {\bm{\varphi}}_\sg$, 
we deduce from (\ref{eqn:4.64}) that 
\begin{eqnarray*}
a_X(\bm{v}, \bm{\varphi}) = (\rot \bm{v}, \rot \bm{\varphi})
= (\rot \bm{v}, \rot \tilde {\bm{\varphi}}_{\sg}) 
= (\bm{u}, \rot \tilde {\bm{\varphi}}_\sg) 
=  (\bm{u}, \rot \bm{\varphi}). 
\end{eqnarray*}
Since $\bm{\varphi} \in \dot X^2(\Omega)$ is arbitrary, 
the above identity implies (\ref{eqn:4.62}).    \qed
%
%\bigskip
%\noindent
%Since $\tilde X^2_{\sigma}(\Omega) \subset \tilde X^2(\Omega)$, from Lemma \ref{lem:4.9} 
%and the 
%{\it unique} slovability of the equation (\ref {eqn:4.38}) with $r =2$ in Lemma \ref{lem:4.7} 
%we obtain the following lemma.  
\begin{lemma}\label{lem:4.10} 
Let $\Omega$ be as in the Assumption and let $\bm{u} \in L^2(\Omega)$.  
\par
%\noindent
{\rm (i)} Suppose that $\bm{v} \in \widetilde X^2(\Omega)$ satisfies 
\begin{equation}\label{eqn:4.65}
a_X(\bm{v}, \bm{\varphi}) = (\bm{u}, \rot \bm{\varphi})
\quad
\mbox{for all $\bm{\varphi}\in \widetilde X^2(\Omega)$}. 
\end{equation} 
Then $\dive \bm{v} =0$.  
\par
%\noindent
{\rm (ii)} Suppose that $\bm{w} \in \widetilde V^2(\Omega)$ satisfies 
\begin{equation}\label{eqn:4.66}
a_V(\bm{w}, \bm{\psi}) = (\bm{u}, \rot \bm{\psi})
\quad
\mbox{for all $\bm{\psi}\in \widetilde V^2(\Omega)$}. 
\end{equation} 
Then $\dive \bm{w} =0$. 
\end{lemma}
%\medskip\noindent
{\it Proof}.  Observe that (\ref{eqn:4.65}) and (\ref{eqn:4.66}) yield (\ref{eqn:4.62}) and (\ref{eqn:4.63}), respectively.  
Indeed,  $\dive \bm{v} = \dive \bm{w} =0$ follows from Lemma \ref{lem:4.9} with the aid of the unique solvability of 
(\ref{eqn:4.38}) and (\ref{eqn:4.40}) in Lemma \ref{lem:4.7}. \qed   
\par
\medskip
%\noindent
We are finally in the position to prove Theorem \ref{thm:4.1}. 
%\bigskip
\begin{thm}\label{thm:4.11} Let $\Omega$ be as in the Assumption and let  $1 < r < \infty$. 
\par
%\noindent 
{\rm (i)} For every $\bm{u} \in L^r(\Omega)$ there exists a unique 
$\bm{v} \in \widetilde X^r(\Omega)$ satisfying  (\ref{eqn:4.38}) and $\dive \bm{v} =0$. In particular,  
\begin{equation}\label{eqn:4.67}
(\rot \bm{v}, \rot \bm{\Phi}) = (\bm{u}, \rot \bm{\Phi})
\quad
\mbox{for all $\bm{\Phi} \in \dot X^{r^\prime}(\Omega)$}.
\end{equation}
Moreover, $\bm{v}$ is subject to the estimate (\ref{eqn:4.39}).  
\par
%\noindent
{\rm (ii)} For every $\bm{u} \in L^r(\Omega)$ there exists a unique 
$\bm{w} \in \widetilde V^r(\Omega)$ ssatisfying (\ref{eqn:4.40}) and  
$\dive \bm{w} =0$. In particular,  
\begin{equation}\label{eqn:4.68}
(\rot \bm{w}, \rot \bm{\Psi}) = (\bm{u}, \rot \bm{\Psi})
\quad
\mbox{for all $\bm{\Psi} \in \dot V^{r^\prime}(\Omega)$}.
\end{equation}
Moreover, $\bm{w}$ is subject to the estimate (\ref{eqn:4.41}). 
\end{thm}
%\bigskip\noindent
{\it Proof}. (i) Due to Lemma \ref{lem:4.7} it suffices to prove that $\dive \bm{v} =0$.  
Since $C^\infty_0(\Omega)$ is dense in $L^r(\Omega)$, there is a sequence $\{\bm{u}_m\}_{m=1}^\infty$ with  $\bm{u}_m \to \bm{u}$ in $L^r(\Omega)$. Lemma \ref{lem:4.7} (i) yields  
that for each $m\in \N$ there exits a unique $\bm{v}_m \in \widetilde X^r(\Omega)$ with  
\begin{equation}\label{eqn:4.69}
a_X(\bm{v}_m, \bm{\Phi}) = (\bm{u}_m, \rot \bm{\Phi})
\quad
\mbox{for all $\bm{\Phi} \in \dot X^{r^\prime}(\Omega)$}. 
\end{equation}
Moreover,
\begin{equation}\label{eqn:4.70}
\|\nabla \bm{v}_m\|_{L^r} \le C\|\bm{u}_m\|_{L^r},
\quad m=1, 2, \cdots, 
\end{equation}
for some  $C=C(\Om, r, R)$. By (\ref{eqn:4.69})  
\begin{equation}\label{eqn:4.71}
a_X(\bm{v}_m, \bm{\varphi}) = (\bm{u}_m, \rot \bm{\varphi})
\quad
\mbox{for all $\bm{\varphi} \in {\mathcal X}(\Omega)$}. 
\end{equation}  
Since $\dot X^2(\Omega) \ni \bm{\phi}  \mapsto (\bm{u}_m, \rot \bm{\phi}) \in \re$ is a continuous functional, 
it follows from (\ref{eqn:4.71}) and Proposition \ref{pr:4.2}(i) that $\bm{v}_m \in \dot X^2(\Om)$.   
Furthermore, by (\ref{eqn:4.71})  
\begin{equation}\label{eqn:4.72}
a_X(\bm{v}_m, \bm{\phi}) = (\bm{u}_m, \rot \bm{\phi})
\quad
\mbox{for all $\bm{\phi} \in \dot X^2(\Omega)$}. 
\end{equation}  
Indeed, since ${\mathcal X}(\Omega)$ is dense in $\widehat X^2(\Omega)$, 
we observe that (\ref{eqn:4.72}) remains true for $\bm{\phi} \in \widehat X^2(\Omega)$.  
For general $\bm{\phi} \in \dot X^2(\Omega)$, we see by (\ref{eqn:3.27}) 
that $\bm{\phi} = \hat {\bm{\phi}} + \sum_{j=1}^3\lambda_j\bm{h}_i$ with 
$\hat {\bm{\phi}} \in \widehat X^2(\Omega)$ and $\lambda_1, \lambda_2, \lambda_3 \in \re$.  
Since $\rot \bm{\phi} = \rot \hat{\bm{\phi}}$ and since (\ref{eqn:4.72}) holds for 
$\bm{\phi} = \hat {\bm{\phi}}$, 
we verify that (\ref{eqn:4.72}) holds for all $\bm{\phi} \in \dot X^2(\Omega)$.   
\par
On the other hand, by (\ref{eqn:3.48}),  $\bm{v}_m \in \dot X^2(\Omega)$ may be  decomposed as 
$$
{\bm{v}}_m = \tilde {\bm{v}}_m + \bm{g}_m
\quad
\mbox{with $\tilde {\bm{v}}_m \in \widetilde X^2(\Omega)$ and 
$\bm{g}_m \in \dot X_{\tiny{\mbox{\rm har}}}(\Omega)$},
\quad m\in N.  
$$
Since $a_X(\bm{g}_m, \bm{\phi}) = 0$ for all $\bm{\phi} \in \dot X^2(\Omega)$, 
we obtain from (\ref{eqn:4.72}) that 
$$
a_X(\tilde {\bm{v}}_m, \bm{\phi})= a_X(\bm{v}_m, \bm{\phi})= (\bm{u}_m, \rot \bm{\phi})
\quad
\mbox{for all $\bm{\phi} \in \dot X^2(\Omega)$}.  
$$
Hence, by Lemma \ref{lem:4.10}(i),  $\dive \tilde {\bm{v}}_m=0$, which yields that $\dive \bm{v}_m=0$.  
Since the mapping $L^r(\Omega) \ni \bm{u} \mapsto\bm{v} \in \widetilde X^r(\Omega)$ with 
$S_r\bm{v} = f_{\bm{u}}$ is continuous, by (\ref{eqn:4.39}) and since $\bm{u}_m \to \bm{u}$ 
in $L^r(\Omega)$, it follows from (\ref{eqn:4.70}) that $\{\bm{v}_m\}_{m=1}^\infty$ is a Cauchy sequence in $\widetilde X^r(\Omega)$.  Thus there exists $\bm{v} \in \tilde X^r(\Omega)$ such that 
\begin{equation}\label{eqn:4.73}  
\nabla \bm{v}_m \to \nabla \bm{v}
\quad
\mbox{in $L^r(\Omega)$ as $m \to \infty$}.  
\end{equation}
Obviously, $\bm{v}$ is a unique solution of (\ref{eqn:4.38}) with $\dive \bm{v} =0$ and satisfies the estimate (\ref{eqn:4.39}). 
\par
(ii) The proof of (ii) is again  parallel to the one of  (i) and hence omitted.   \qed

\section{Construction of the scalar potential}
In this section we construct the scalar potential $p$ for the given function  $\bm{u}\in L^r(\Omega)$ in (\ref{eqn:2.6}) and (\ref{eqn:2.10}) in Theorems \ref{thm:2.2} and \ref{thm:2.3}, 
respectively. We start by describing  the scalar potential $p \in \dot H^{1, r}(\Omega)$ in (\ref{eqn:2.6}),  
which is determined by the solution of the weak Neumann problem of the Poission equation.
\begin{proposition}\label{thm:5.1}{\rm(}Simader-Sohr \cite[Theorem 1.4]{SiSo2}{\rm)} 
Let $\Omega$ be as in the Assumption and let $1< r < \infty$. 
For every $\bm{u} \in L^r(\Omega)$ there is  a unique $p \in \dot H^{1, r}(\Omega)$ such that 
\begin{equation}\label{eqn:5.1}
(\nabla p, \nabla \psi) = (\bm{u}, \nabla \psi) 
\quad
\mbox{for all $\psi \in \dot H^{1, r^\prime}(\Omega)$}.  
\end{equation} 
Moreover,  $p$ is subject to the estimate 
\begin{equation}\label{eqn:5.2}
\|\nabla p\|_{L^r} \le C\|\bm{u}\|_{L^r}, 
\end{equation}
for some  $C=C(\Omega, r)$.  
\end{proposition}  
We next investigate the scalar potential $p$ in Theorem \ref{thm:2.3}. 
In comparison with $p$ in (\ref{eqn:5.1}), 
we need to solve the weak Dirichlet problem for the Poisson equation. 
To this end, we make use of the Laplacian 
$-\Delta_r: \hat H^{1, r}_0(\Omega) \to \hat H^{1, r^\prime}_0(\Omega)^\ast$ given in (\ref{eqn:3.3}).  
We start with the case where 
$3/2 < r < \infty$. 
%\medskip\noindent
\begin{lemma}\label{lem:5.2} Let $\Omega$ be as in the Assumption. 
\par
{\rm (i)} Let $3/2 < r < 3$. 
For every $\bm{u} \in L^r(\Omega)$ there is a unique $p \in \widehat H^{1, r}_0(\Omega)$ satisfying 
\begin{equation}\label{eqn:5.3}
(\nabla p, \nabla \phi) = (\bm{u}, \nabla \phi) \mbox{ for all } \phi \in \widehat H^{1, r^\prime}_0(\Omega) 
\mbox{ and }  
\end{equation} 
\begin{equation}\label{eqn:5.4}
\|\nabla p\|_{L^r} \le C\|\bm{u}\|_{L^r}, 
\end{equation}
for some $C=C(\Omega, r)$.  
\par
%\noindent
{\rm (ii)} Let $3 \le r < \infty$. 
For every $\bm{u} \in L^r(\Omega)$ there is a unique $p \in \widetilde H^{1, r}_0(\Omega)$ 
satisfying (\ref{eqn:5.3}).  
Moreover,  $p$ is subject to the estimate (\ref{eqn:5.4}).    
\end{lemma}
%\medskip\noindent  
{\it Proof}. (i) For given $\bm{u} \in L^r(\Omega)$ we define the functional 
$f_{\bm{u}}$ by $\langle f_{\bm{u}}, \phi\rangle = (\bm{u}, \nabla \phi)$ 
for $\phi \in \widehat H^{1, r^\prime}_0(\Omega)$.  
Obviously, $f_{\bm{u}} \in \widehat H^{1, r^\prime}_0(\Omega)^\ast$ satisfies 
$\|f_{\bm{u}}\|_{\widehat H^{1, r}_0(\Omega)^\ast} \le \|\bm{u}\|_{L^r}$.  
Hence, the assertion is a consequence of Proposition \ref{pr:3.2} (ii) and (iii).  
\par
%\noindent
(ii) By (\ref{eqn:3.2}) and Proposition \ref{pr:3.2} (iv),  
$-\Delta_r$ is a {\it bijective} operator from $\widetilde H^{1, r}_0(\Omega)$ onto 
$\widehat H^{1, r^\prime}_0(\Omega)^\ast$. 
Hence, we obtain the assertion in the same way as above.  
\qed
\par
\bigskip
We next consider the case where $1 < r \leqq 3/2$. As stated in Proposition \ref{pr:3.2} (iii),  
$-\Delta_r$ is not surjective in this case. Hence, we consider the larger space 
$\dot H^{1, r}_0(\Omega)$ for the solvability of (\ref{eqn:5.3}). 
%\medskip\noindent
\begin{lemma}\label{lem:5.3}{\rm (}Simader-Sohr \cite[Theorem 1.2]{SiSo3}, 
Pr\"uss-Simonett \cite[Theorem 7.4.3]{PuSi}, Shibata \cite[Theorem 3.2]{Sh}{\rm)} 
Let $\Omega$ be as in the Assumption and let  $1 < r \le 3/2$. 
For every $\bm{u} \in L^r(\Omega)$, there exists a unique $p\in \dot H^{1, r}_0(\Omega)$ satisfying  (\ref{eqn:5.3}) and  (\ref{eqn:5.4}).  
\end{lemma}
%
%\medskip\noindent 
%
%
{\it Proof}.  
We give here a different proof from \cite{SiSo3}, \cite{PuSi} and \cite{Sh}.  
Let us first show first that there exists a scalar function 
$\bar p \in \dot H^{1, r}(\Omega)$ with $\bar p|_{\pt \Om} = \mbox{const}.$ satisfying (\ref{eqn:5.3}).    
To this end, we use the cut-off function $\eta \in C^\infty_0(\re^3)$ defined in the proof of Lemma \ref{lem:3.9}. 
For given $\bm{u} \in L^r(\Omega)$, consider the functional $g_{\bm{u}}$ defined by 
$\langle g_{\bm{u}}, \psi\rangle \equiv (\bm{u}, \nabla \psi)$ for 
$\psi \in \widetilde H^{1, r^\prime}_0(\Omega)$. 
Then  $g_{\bm{u}} \in \widetilde H^{1, r^\prime}_0(\Omega)^\ast$ and 
$\|g_{\bm{u}}\|_{\widetilde H^{1, r^\prime}_0(\Omega)^\ast} \le \|\bm{u}\|_{L^r}$. 
By Lemma \ref{lem:3.3}, there exists a unique $\pi \in \widehat H^{1, r}_0(\Omega)$ such that 
\begin{equation}\label{eqn:5.5} 
(\nabla \pi, \nabla \psi) = (\bm{u}, \nabla \psi)
\quad
\mbox{for all $\psi \in \widetilde H^{1, r^\prime}_0(\Omega)$}. 
\end{equation}   
and satisfies 
\begin{equation}\label{eqn:5.6}
\|\nabla \pi\|_{L^r} \le C\|\bm{u}\|_{L^r},
\end{equation} 
for some  $C = C(\Omega, r)$. 
Similarly, now $\bm{u}$ being replaced by $\nabla \eta$, 
there is a unique $\alpha \in \widehat H^{1, r}_0(\Omega)$ 
such that  
\begin{equation}\label{eqn:5.7} 
(\nabla \alpha, \nabla \psi) = (\nabla \eta, \nabla \psi)
\quad
\mbox{for all $\psi \in \widetilde H^{1, r^\prime}_0(\Omega)$},  
\end{equation} 
where $\eta \in C^\infty_0(\re^3)$ is as in the proof of Lemma \ref{lem:3.9}. Defining $q$ by $q(x) \equiv \eta(x) - \alpha(x)$ for $x \in \Omega$, we see that  
\begin{equation}\label{eqn:5.8}
(\nabla q, \nabla q_0) \ne 0, 
\end{equation} 
where $q_0 \in \bigcap_{s>3/2}\dot H^{1, s}_0(\Omega)$ is a harmonic function given by (\ref{eqn:3.2}) 
with $q_0(x) \to 1$ as $|x|\to\infty$.    
\par
Assume (\ref{eqn:5.8}) for the time being. Define 
\begin{equation}\label{eqn:5.9}
\bar p(x) \equiv \pi(x) + \lambda q(x),\quad x \in \Omega\quad\mbox{with $\lambda \equiv \dfrac{(\bm{u}, \nabla q_0)}{(\nabla q, \nabla q_0)}$}.  
\end{equation}
Then $\nabla \bar p \in L^r(\Omega)$ and $\bar p|_{\pt\Om} = \lambda$.  Moreover, 
$\bar p$ satisfies 
\begin{equation}\label{eqn:5.10}
(\nabla \bar p, \nabla \phi) = (\bm{u}, \nabla \phi)
\quad
\mbox{for all $\phi \in \widehat H^{1, r^\prime}_0(\Omega)$}.
\end{equation}
Indeed, since $1 < r \le 3/2$, by (\ref{eqn:3.2}) every $\phi \in \widehat H^{1, r^\prime}_0(\Omega)$ 
may be  decomposed as 
$\phi = \psi + \mu q_0$ for some $\psi \in \widetilde H^{1, r^\prime}_0(\Omega)$ and $\mu \in \re$.   
Since $(\nabla \pi, \nabla q_0) =0$, we obtain from (\ref{eqn:5.5}), (\ref{eqn:5.7}) and  (\ref{eqn:5.9})  
\begin{eqnarray*}
(\nabla \bar p, \nabla \phi) 
&=& (\nabla \pi + \lambda\nabla q, \nabla \psi + \mu\nabla q_0) \\
&=& (\nabla \pi,  \nabla \psi) + \mu(\nabla \pi, \nabla q_0) + \lambda(\nabla q, \nabla \psi) 
+ \lambda\mu(\nabla q, \nabla q_0) \\
&=& (\bm{u}, \nabla \psi) + \lambda(\nabla \eta -\nabla \al, \nabla \psi) + \mu(\bm{u}, \nabla q_0) \\
&=& (\bm{u}, \nabla (\psi + \mu q_0)) \\
&=& (\bm{u}, \nabla \phi), 
\end{eqnarray*}
which implies (\ref{eqn:5.10}). Setting $p(x) \equiv \bar p(x) - \lambda$ implies $p \in \dot H^{1, r}_0(\Omega)$. Since $\nabla p = \nabla \bar p$, (\ref{eqn:5.3}) is implied by  (\ref{eqn:5.10}).  
Moreover, (\ref{eqn:5.6}) and  (\ref{eqn:5.9}) yield 
$$
\|\nabla p\|_{L^r} = \|\nabla \bar p\|_{L^r} \le C\|\bm{u}\|_{L^r},
$$
which implies (\ref{eqn:5.4}).  
\par
Finally, concerning uniqueness, suppose that $\hat p \in \dot H^{1, r}_0(\Omega)$ is another solution of (\ref{eqn:5.3}).  Then $p_\ast \equiv p - \hat p \in\dot H^{1, r}_0(\Omega)$ satisfies 
$(\nabla p_\ast, \nabla \phi) = 0$ for all $\phi \in \widehat H^{1, r^\prime}_0(\Omega)$. Since $1 < r \le 3/2 $, it follows from \cite[Theorem A]{KoSo} that $p_\ast \equiv 0$ on $\Omega$, 
which yields  the desired uniqueness property.      
\par
It only remains to prove (\ref{eqn:5.8}). Since $q_0 \in \mbox{Ker}(-\Delta_{r^\prime})$ 
and since $\al \in \widehat H^{1, r}_0(\Omega)$, 
we see that $(\nabla \al, \nabla q_0) =0$. Hence, 
$$
(\nabla q, \nabla q_0) = (\nabla \eta - \nabla \al, \nabla q_0) = (\nabla \eta, \nabla q_0) = 
\int_{\pt\Om}\frac{\pt q_0}{\pt\bm{\nu}}dS.  
$$
Assume that (\ref{eqn:5.8}) fails. The above identity yields
\begin{equation}\label{eqn:5.11}
\int_{\pt\Om}\frac{\pt q_0}{\pt\bm{\nu}}dS = 0.  
\end{equation}
Since $q_0$ is harmonic in $\Omega$ with $q_0|_{\pt \Om}=0$ and $q_0(x) \to 1$ as $|x| \to \infty$, the maximum principle yields  $0 \le q_0(x) \le 1$ for all $x \in \Omega$ and thus  
$\dfrac{\pt q_0}{\pt\bm{\nu}} \le 0$ on $\pt\Omega$. By (\ref{eqn:5.11})  
$\dfrac{\pt q_0}{\pt\bm{\nu}}\equiv 0$ on $\pt\Omega$, which yields  
$(\nabla q_0, \nabla \Psi) =0$ for all $\Psi \in C^\infty_0(\bar\Omega)$.   
Finally, since $q_0 \in \dot H^{1, r^\prime}_0(\Omega)$, the uniqueness property of the weak Neumann problem for  the Poisson equation implies  
$q_0(x) \equiv \mbox{const}$ on $\Omega$, which causes a contradiction.    
\qed
\section{Proof of the Main Theorems \ref{thm:2.2} and \ref{thm:2.3}} 
\subsection{Generalized Stokes formula} \mbox{}
\par
In this section, we finally prove our main results stated in  Theorems \ref{thm:2.2} and \ref{thm:2.3}. To this end, let us first consider  the generalized Stokes formula in $\Omega$ 
which holds on $E^r_{\tiny{\mbox{\rm div}}}(\Omega)$ and $E^r_{\tiny{\mbox{\rm rot}}}(\Omega)$  defined by 
\begin{eqnarray*}
E^r_{\tiny{\mbox{\rm div}}}(\Omega)&\equiv& \{\bm{u}\in L^r(\Omega); \dive \bm{u} \in L^r(\Omega)\} 
\quad
\mbox{equipped  the norm $\|\bm{u}\|_{E^r_{\tiny{\mbox{\rm div}}}} = \|\bm{u}\|_{L^r} + \|\dive \bm{u}\|_{L^r}$} \\
E^r_{\tiny{\mbox{\rm rot}}}(\Omega)&\equiv& \{\bm{u}\in L^r(\Omega); \rot \bm{u} \in L^r(\Omega)\}
\quad 
\mbox{equipped  the norm $\|\bm{u}\|_{E^r_{\tiny{\mbox{\rm rot}}}} = \|\bm{u}\|_{L^r} + \|\rot \bm{u}\|_{L^r}$}. 
\end{eqnarray*} 
We state the generalized Stokes formula in exterior domains. 
The following generalized Stokes formulae (\ref{eqn:6.1}) and (\ref{eqn:6.2}) are well-known in the case of 
bounded domains $\Omega$. 
There seems to be no proof in the literature for the situation of exterior domains $\Omega$.
%
%\medskip
\begin{lemma}\label{lem:6.1}
Let $\Omega$ be as in the Assumption and let  $1< r < \infty$.
There exist bounded operators $\gamma_{\bm{\nu}}$ and $\tau_{\bm{\nu}}$ on $E^r_{\tiny{\mbox{\rm div}}}(\Omega)$ and on $E^r_{\tiny{\mbox{\rm rot}}}(\Omega)$ with the properties that 
\begin{eqnarray*}
&& \gamma_{\bm{\nu}}:E^r_{\tiny{\mbox{\rm div}}}(\Omega) \ni \bm{u}   \mapsto 
\gamma_{\bm{\nu}} \bm{u} \in H^{1-\frac{1}{r^\prime}, r^\prime}(\pt\Om)^\ast, 
\quad \gamma_{\bm{\nu}}\bm{u} = \bm{u}\cdot \bm{\nu}|_{\pt\Om} \quad\mbox{for $\bm{u} \in C^1_0(\bar\Omega)$},  \\
&& \tau_{\bm{\nu}}:E^r_{\tiny{\mbox{\rm rot}}}(\Omega) \ni \bm{u} \mapsto 
\tau_{\bm{\nu}} \bm{u} \in H^{1-\frac{1}{r^\prime}, r^\prime}(\pt\Om)^\ast, 
\quad \tau_{\bm{\nu}} \bm{u} = \bm{u}\times \bm{\nu}|_{\pt\Om} \quad\mbox{for $\bm{u} \in C^1_0(\bar\Omega)$},  
\end{eqnarray*}
respectively.  Furthermore, the following generalized Stokes formulae hold true
\begin{eqnarray}
&&  (\bm{u}, \nabla q) + (\dive \bm{u}, q) = \langle \gamma_{\bm{\nu}} \bm{u}, \gamma_0 q\rangle_{\pt\Om} 
\quad
\mbox{for all $\bm{u} \in E^r_{\tiny{\mbox{\rm div}}}(\Omega)$ and all $q \in H^{1, r^\prime}(\Omega)$}, 
\label{eqn:6.1} \\
&&  (\bm{u}, \rot \bm{\psi}) = 
(\rot \bm{u}, \bm{\psi}) + \langle \tau_{\bm{\nu}} \bm{u}, \gamma_0 \bm{\psi}\rangle_{\pt\Om} 
\quad
\mbox{for all $\bm{u} \in E^r_{\tiny{\mbox{\rm rot}}}(\Omega)$ and all $\bm{\psi} \in H^{1, r^\prime}(\Omega)$},
\label{eqn:6.2}
\end{eqnarray}
where $\gamma_0$ denotes the usual trace operator from $H^{1, r^\prime}(\Omega)$ onto $H^{1-\frac{1}{r^\prime}, r^\prime}(\pt\Om)$, and 
$\langle\cdot, \cdot\rangle_{\pt\Om}$ denotes the duality pairing between $H^{1-\frac{1}{r^\prime}, r^\prime}(\pt\Om)^\ast$ and $H^{1-\frac{1}{r^\prime}, r^\prime}(\pt\Om)$. 
\end{lemma}

%The generalized Stokes formulae (\ref{eqn:6.1}) and (\ref{eqn:6.2}) are well-known in bounded domains.  
%For the sake of completeness, we give here their proof in our exterior domain $\Omega$. 

\noindent
{\it Proof of Lemma \ref{lem:6.1}}.  
We start by  showing that $C^\infty_0(\bar\Omega)$ is dense in both $ E^r_{\tiny{\mbox{\rm div}}}(\Omega)$ and  $E^r_{\tiny{\mbox{\rm rot}}}(\Omega)$ 
for all $1 < r < \infty$.  Indeed, let $\bm{u} \in E^r_{\tiny{\mbox{\rm div}}}(\Omega)$.  Then there exists an extension $\tilde {\bm{u}} \in L^r(\re^3)$ of $\bm{u}$ satisfying   
$\dive \tilde {\bm{u}} \in L^r(\re^3)$ and  
\begin{equation}\label{eqn:6.3.0}
\bm{u}(x) = \tilde{\bm{u}}(x) 
\quad
\mbox{for $x \in \Omega$}, 
\quad 
\|\tilde{\bm{u}}\|_{L^r(\re^3)} \le C\|\bm{u}\|_{L^r}, 
\quad
\|\dive \tilde{\bm{u}}\|_{L^r(\re^3)} \le C\|\dive \bm{u}\|_{L^r}, 
\end{equation}
for some $C=C(\Omega,r)$ independent of $\bm{u}$. For $\delta >0$ and $k\in \N$, we set $\bm{u}_{k, \delta}(x)\equiv \zeta_k(x)(\rho_\delta\ast\tilde{\bm{u}})(x)$ for $ x\in \re^3$, where the 
mollifier $\rho_\delta$ is defined by $\rho_\delta(x) \equiv \delta^{-3}\rho(x/\delta)$ for some $\rho\in C^{\infty}_0(B_1)$ satisfying  $\rho(x) \ge 0$ for all $x \in \re^3$ and  
$\int_{|x| < 1}\rho(x)dx =1$. Let  $\{\zeta_k\}_{k=1}^\infty$ be a sequence of cut-off functions 
defined by $\zeta_k(x) = \zeta(x/k)$ for $k\in \N$, where $\zeta\in C^\infty_0(B_2)$ satisfies 
$0 \le \zeta(x) \le 1$ for all $x \in \re^3$ and $\zeta(x) =1$ for $x \in B_1$. 
Notice that $\zeta_k(x)=1$ for $x \in B_k$ and $\zeta_k(x) =0$ for $x \in \re^3\setminus B_{2k}$ 
with $\|\nabla \zeta_k\|_{L^\infty(\re^3)}\le Ck^{-1}$ for all $k\in \N$, where  $C = \|\nabla \zeta\|_{L^\infty(\re^3)}$.  Obviously, 
$\bm{u}_{k, \delta} \in C^\infty_0(\bar\Omega)$ for all $\delta>0$ and all $k\in \N$. For every $\ep>0$ there is $\delta_0=\delta_0(\ep) >0$ 
such that $\|\rho_{\delta_0}\ast \tilde{\bm{u}} - \tilde{\bm{u}}\|_{L^r(\re^3)} < \ep/4$. For such a $\delta_0$, there exists $k_0\in \N$ such that 
$\|\zeta_{k_0}(\rho_{\delta_0}\ast \tilde{\bm{u}})-\rho_{\delta_0}\ast \tilde{\bm{u}}\|_{L^r(\re^3)} < \ep/4$.
Hence,  
\begin{eqnarray}
\|\bm{u}_{k_0, \delta_0} - \bm{u}\|_{L^r}
&=& \|\zeta_{k_0}(\rho_{\delta_0}\ast \tilde{\bm{u}}) - \tilde{\bm{u}}\|_{L^r} \nonumber \\
&\le& \|\zeta_{k_0}(\rho_{\delta_0}\ast \tilde{\bm{u}})-\rho_{\delta_0}\ast \tilde{\bm{u}}\|_{L^r(\re^3)} 
+ \|\rho_{\delta_0}\ast \tilde{\bm{u}} - \tilde{\bm{u}}\|_{L^r(\re^3)}\label{eqn:6.4.0} \\
&<& \ep/2.\nonumber  
\end{eqnarray}
Similarly,  
\begin{equation}\label{eqn:6.5.0}
\|\dive \rho_{\delta_0}\ast \tilde{\bm{u}} - \dive \tilde{\bm{u}}\|_{L^r(\re^3)}
= \|\rho_{\delta_0}\ast \dive \tilde{\bm{u}} - \dive \tilde{\bm{u}}\|_{L^r(\re^3)} < \ep/4.
\end{equation}
Since $\|\nabla\zeta_k\cdot(\rho_{\delta}\ast\tilde{\bm{u}})\|_{L^r(\re^3)} \le $ $\|\nabla \zeta_k\|_{L^\infty(\re^3)}\|\rho_{\delta}\ast\tilde{\bm{u}}\|_{L^r(\re^3)}
\le $ $ Ck^{-1}\|\bm{u}\|_{L^r}$
for all $\delta>0$ and all $k\in \N$ with $C=C(\Omega, r)$,  we may assume that 
\begin{equation}\label{eqn:6.6.0}
\|\zeta_{k_0}\dive \rho_{\delta_0}\ast \tilde{\bm{u}} - \dive \rho_{\delta_0}\ast \tilde{\bm{u}}\|_{L^r(\re^3)} 
+ \|\nabla\zeta_{k_0}\cdot\rho_{\delta_0}\ast\tilde{\bm{u}}\|_{L^r(\re^3)} 
< \ep/4
\end{equation}
holds for $k_0$. Furthermore, since 
$\dive \bm{u}_{k_0, \delta_0} = \zeta_{k_0}\dive \rho_{\delta_0}\ast\tilde{\bm{u}} 
+ \nabla \zeta_{k_0}\cdot \rho_{\delta_0}\ast\tilde{\bm{u}}$,  
assertions (\ref{eqn:6.5.0}) and (\ref{eqn:6.6.0}) yield 
\begin{eqnarray}
\|\dive \bm{u}_{k_0, \delta_0} - \dive \bm{u}\|_{L^r} \nonumber 
&\le & \|\dive \bm{u}_{k_0, \delta_0} - \dive \tilde{\bm{u}}\|_{L^r(\re^3)} \nonumber \\
&\le& 
\|\zeta_{k_0}\dive \rho_{\delta_0}\ast\tilde{\bm{u}} - \dive\rho_{\delta_0}\ast \tilde{\bm{u}}\|_{L^r(\re^3)} 
+ \|\dive\rho_{\delta_0}\ast \tilde{\bm{u}} - \dive \tilde{\bm{u}}\|_{L^r(\re^3)} \label{eqn:6.7.0}\\ 
&&  + \|\nabla\zeta_{k_0}\cdot(\rho_{\delta_0}\ast\tilde {\bm{u}})\|_{L^r(\re^3)} 
\nonumber \\
&<& \ep/2. \nonumber  
\end{eqnarray}
By (\ref{eqn:6.4.0}) and (\ref{eqn:6.7.0}), 
$\|\bm{u}_{k_0, \delta_0} -  \bm{u}\|_{E^r_{\tiny{\mbox{\rm div}}}} < \ep$. 
Since $\bm{u} \in E^r_{\tiny{\mbox{\rm div}}}(\Omega)$ was arbitrary, $C^\infty_0(\bar\Omega)$ is dense in $E^r_{\tiny{\mbox{\rm div}}}(\Omega)$. An approximating argument similar to the above one yields that 
$C^\infty_0(\bar\Omega)$ is dense in $E^r_{\tiny{\mbox{\rm rot}}}(\Omega)$.
\par
In order to prove  (\ref{eqn:6.1}) note that for every $\bm{u}\in E^r_{\tiny{\mbox{\rm div}}}(\Omega)$, there exists a sequence $\{\bm{u}_m\}_{m=1}^\infty \subset C^\infty_0(\bar\Omega)$ with  
\begin{equation}\label{eqn:6.8.0}
\bm{u}_m \to \bm{u}, \quad \dive\bm{u}_m \to \dive \bm{u} 
\quad
\mbox{in $L^r(\Omega)$ as $m \to \infty$}.  
\end{equation}
Since the trace operator $\gamma_0: H^{1, r^\prime}(\Omega) \to H^{1-\frac{1}{r^\prime}, r^\prime}(\pt\Omega)$ is surjective, for every $p\in H^{1-\frac{1}{r^\prime}, r^\prime}(\pt\Om)$ there exists  
$q \in H^{1, r^\prime}(\Omega)$ satisfying  $\gamma_0q = p$ and $\|q\|_{H^{1, r^\prime}(\Omega)} \le C\|p\|_{H^{1-\frac{1}{r^\prime}, r^\prime}(\pt\Omega)}$ for some $C=C(\Om, r)$.     
Considering  the functional $X_m(p)$ given by 
$$
X_m(p) \equiv (\nabla q, \bm{u}_m) + (q, \dive \bm{u}_m), \quad m \in \N  
$$
we verify that the value of $X_{m}(p)$ is independent of $q$, as long as $q \in H^{1, r^\prime}(\Omega)$ and $\gamma_0q = p$.  
Indeed, assume that there are $q_1, q_2 \in H^{1, r^\prime}(\Omega)$ with  $\gamma_0q_1= \gamma_0q_2 =p$. We aim to show that
\begin{equation}\label{eqn:6.9.0}  
(\nabla q_1, \bm{u}_m) + (q_1, \dive \bm{u}_m) = (\nabla q_2, \bm{u}_m) + (q_2, \dive \bm{u}_m), \quad m \in \N, 
\end{equation}
which for $\tilde q \equiv q_1 - q_2$ is equivalent to 
$$
(\nabla \tilde q, \bm{u}_m) + (\tilde q, \dive \bm{u}_m) = 0, \quad m \in \N.
$$
Since $\gamma_0\tilde q =0$ and since $\mbox{supp }\bm{u}_m$ is compact in $\re^3$, we obtain  
$$
(\nabla \tilde q, \bm{u}_m) + (\tilde q, \dive \bm{u}_m)= \int_\Omega\dive (\tilde q\bm{u}_m)dx = \int_{\pt\Om}\gamma_0\tilde q \bm{u_m}\cdot\bm{\nu}dS = 0, \quad m \in \N,
$$
from which we obtain (\ref{eqn:6.9.0}). The functional $X_m(p)$ is hence well-defined for every $p \in H^{1-\frac{1}{r^\prime}, r^\prime}(\pt\Om)$.  
Furthermore, by (\ref{eqn:6.8.0}) 
$$
X_{\bm{u}}(p) \equiv \lim_{m\to\infty}X_m(p) = (\nabla q, \bm{u}) + (q, \dive \bm{u})
$$
with $|X_{\bm{u}}(p)|\le C\|\bm{u}\|_{E^r_{\tiny{\mbox{\rm div}}}}\|p\|_{H^{1-\frac{1}{r^\prime}, r^\prime}(\pt\Om)}$. This means that for every $\bm{u} \in E^r_{\tiny{\mbox{\rm div}}}(\Omega)$,  
$X_{\bm{u}}$ is a continuous functional on $H^{1-\frac{1}{r^\prime}, r^\prime}(\pt\Om)$ satisfying  
$\|X_{\bm{u}}\|_{H^{1-\frac{1}{r^\prime}, r^\prime}(\pt\Om)^\ast} \le C\|\bm{u}\|_{E^r_{\tiny{\mbox{\rm div}}}}$. 
Finally, defining $\gamma_{\bm{\nu}}\bm{u}$ by $\langle \gamma_{\bm{\nu}}\bm{u}, p\rangle_{\pt\Om}\equiv X_{\bm{u}}(p)$ for $p \in H^{1-\frac{1}{r^\prime}, r^\prime}(\pt\Om)$, 
we see that $\gamma_{\bm{\nu}}\bm{u} \in H^{1-\frac{1}{r^\prime}, r^\prime}(\pt\Om)^\ast$ 
enjoys the property (\ref{eqn:6.1}). 
We also verify that $\gamma_{\bm{\nu}}\bm{u} = \bm{u}\cdot\bm{\nu}|_{\pt\Om}$ 
provided $\bm{u} \in C^1_0(\Omega)$.  
\par
The validity of (\ref{eqn:6.2}) is proved by smilar argument as above.\qed 
\par
\bigskip 
We note that  Lemma \ref{lem:6.1} implies that 
\begin{eqnarray}
X^r_{\tiny{\mbox{\rm har}}}(\Omega) &=& 
\{\bm{h} \in E^r_{\tiny{\mbox{\rm div}}}(\Omega)\cap E^r_{\tiny{\mbox{\rm rot}}}(\Omega); 
\dive \bm{h} =0, \, \rot \bm{h}=\bm{0},\, \gamma_{\bm{\nu}} \bm{h} = 0\}, \label{eqn:6.3} \\
\quad
V^r_{\tiny{\mbox{\rm har}}}(\Omega)& = &\{\bm{h} \in E^r_{\tiny{\mbox{\rm div}}}(\Omega)\cap E^r_{\tiny{\mbox{\rm rot}}}(\Omega); \dive \bm{h} =0, \, \rot \bm{h}=\bm{0}, \, \tau_{\bm{\nu}} \bm{h} = \bm{0}\}.  
\label{eqn:6.4}
\end{eqnarray}
for all $1 < r < \infty$.  
%
%\bigskip
\subsection{Proof of Theorem \ref{thm:2.2}} \mbox{} 
\par
Let $1 < r < \infty$.  For given $\bm{u} \in L^r(\Omega)$ we take $p \in \dot H^{1, r}(\Omega)$ and $\bm{w} \in \dot V^r_\sg(\Omega)$ given by Theorem \ref{thm:5.1} and Theorem \ref{thm:4.1} (ii) and define 
$\bm{h}$ by 
\begin{equation}\label{eqn:6.5}
\bm{h} \equiv \bm{u} - \rot \bm{w} - \nabla p.
\end{equation} 
Let us first show that $\bm{h} \in X^r_{\tiny{\mbox{\rm har}}}(\Omega)$. Obviously, $\bm{h} \in L^r(\Omega)$. Furthermore, 
$\dive \bm{h} =0$ and $\rot \bm{h}= \bm{0}$ in the sense of distribution in $\Omega$, i.e.,  
\begin{eqnarray}
&& (\bm{h}, \nabla \psi) = 0 
\quad
\mbox{for all $\psi \in C^\infty_0(\Omega)$}, \label{eqn:6.6} \\
&& 
(\bm{h}, \rot \bm{\Psi}) = 0 
\quad
\mbox{for all $\bm{\Psi} \in C^\infty_0(\Omega)$}. \label{eqn:6.7}
\end{eqnarray}
Indeed, by (\ref{eqn:5.1}) and (\ref{eqn:6.2})
$$
(\bm{h}, \nabla \psi) = (\bm{u} - \nabla p, \nabla \psi) - (\rot \bm{w}, \nabla \psi)
= -(\bm{w}, \rot \nabla \psi) = 0
\quad
\mbox{for all $\psi \in C^\infty_0(\Omega)$}, 
$$
which implies (\ref{eqn:6.6}).  By (\ref{eqn:4.3}) and (\ref{eqn:6.1})  
$$
(\bm{h}, \rot \bm{\Psi}) 
= (\bm{u}-\rot \bm{w}, \rot \bm{\Psi}) - (\nabla p, \rot \bm{\Psi}) 
= (p, \dive \rot \bm{\Psi}) = 0
\quad
\mbox{for all $\bm{\Psi} \in C^\infty_0(\Omega)$}, 
$$
which implies (\ref{eqn:6.7}). 
\par
We next show that $\gamma_{\bm{\nu}} \bm{h} = 0$. By (\ref{eqn:5.1}) and (\ref{eqn:6.1})  
$$
0 = (\bm{u} - \nabla p, \nabla q) = \langle \gamma_{\bm{\nu}}(\bm{u} - \nabla p), \gamma_0q\rangle_{\pt\Omega} 
\quad
\mbox{for all $q \in H^{1, r^\prime}(\Omega)$}.  
$$
Since $\gamma_0:H^{1, r^\prime}(\Omega) \to H^{1-\frac{1}{r^\prime}}(\pt\Omega)$ is surjective, the above identity implies that $\gamma_{\bm{\nu}}(\bm{u}-\nabla p)=0$.   
Moreover, since $\bm{w} \in \dot V^r(\Omega)$,  $\bm{w}\times \bm{\nu}|_{\pt \Omega}= 0$ and by (\ref{eqn:6.2}) we deduce  
$$
(\rot \bm{w}, \nabla q) = (\bm{w}, \rot (\nabla q))- \langle \bm{w}\times\bm{\nu}, \nabla q\rangle_{\pt\Omega} =0
\quad
\mbox{for all $q \in H^{2, r^{\prime}}(\Omega)$}.  
$$
Since $H^{2, r^\prime}(\Omega)$ is dense in $H^{1, r^\prime}(\Omega)$, 
from this identity we obtain 
$$
(\rot \bm{w}, \nabla q) = 0 \quad \mbox{for all $q \in H^{1, r^{\prime}}(\Omega)$}.  
$$
0n the other hand, by (\ref{eqn:6.1}) we obtain  
$$
0 = (\rot \bm{w}, \nabla q) = 
-(\dive (\rot \bm{w}), p) + \langle\gamma_{\bm{\nu}} \rot \bm{w}, \gamma_0q\rangle_{\pt\Om} 
= \langle\gamma_{\bm{\nu}} \rot \bm{w}, \gamma_0q\rangle_{\pt\Om}
\quad
\mbox{for all $q \in H^{1, r^{\prime}}(\Omega)$}.  
$$
Again by surjectivity of $\gamma_0$, we obtain $\gamma_{\bm{\nu}} \rot \bm{w} =0$ and hence  $\gamma_{\bm{\nu}}\bm{h} =0$.   
It thus follows from (\ref{eqn:6.3}), (\ref{eqn:6.6}) and (\ref{eqn:6.7}) that $\bm{h} \in X^r_{\tiny{\mbox{\rm har}}}(\Omega)$ satisfies  (\ref{eqn:2.6}). The estimate (\ref{eqn:2.7}) is a consequence 
of (\ref{eqn:4.4}), (\ref{eqn:5.2}) and (\ref{eqn:6.5}).   
\par
It remains to prove the uniqueness of the expression (\ref{eqn:2.6}).  Assume that $\tilde {\bm{h}} \in X^r_{\tiny{\mbox{\rm har}}}(\Omega)$, $\tilde {\bm{w}} \in \dot V^r_\sg(\Omega)$ and 
$\tilde p \in \dot H^{1, r}(\Omega)$ fulfill (\ref{eqn:2.8}).  Set  $\hat {\bm{h}} \equiv \bm{h} - \tilde {\bm{h}}$, $\hat {\bm{w}} \equiv \bm{w} - \tilde {\bm{w}}$ and 
$\hat p \equiv p - \tilde p$.  Then
\begin{equation}\label{eqn:6.8}
\hat {\bm{h}} + \rot \hat {\bm{w}} + \nabla \hat p =0.
\end{equation} 
Since $\gamma_{\bm{\nu}}\hat {\bm{h}} =0$ and $\hat {\bm{w}}\times {\bm{\nu}}|_{\pt\Om} =\bm{0}$, we obtain  by (\ref{eqn:6.1}) and (\ref{eqn:6.2}) 
\begin{eqnarray*}
&&(\nabla \hat p, \nabla \psi) \\
&=&- (\hat {\bm{h}}, \nabla \psi) - (\rot \hat {\bm{w}}, \nabla \psi) 
= (\dive \hat {\bm{h}}, \psi) -\langle \gamma_{\bm{\nu}} \hat {\bm{h}}, \gamma_0 \psi\rangle_{\pt\Om} 
-(\hat {\bm{w}}, \rot(\nabla \psi)) + \langle \hat {\bm{w}}\times \bm{\nu}, \nabla\psi\rangle_{\pt\Om} \\
&=& 0
\end{eqnarray*}
for all $\psi \in C^\infty_0(\bar \Omega)$.  By the uniqueness of the weak Neumann problem of the Poisson equation, $\nabla \hat p\equiv 0$ in $\Omega$, which implies  
\begin{equation}\label{eqn:6.9}
\nabla p = \nabla \tilde p.
\end{equation}  
Furthermore,  $\rot \hat {\bm{w}} = -\hat {\bm{h}}$. Since $\rot \hat {\bm{h}}=\bm{0} $ and $\bm{\Psi}\times\bm{\nu}|_{\pt\Omega} =\bm{0}$ for $\bm{\Psi} \in \dot V^{r^\prime}(\Om)$, 
assertion (\ref{eqn:6.2}) implies 
\begin{eqnarray*}
a_V(\hat {\bm{w}}, \bm{\Psi}) &=& (\rot \hat {\bm{w}}, \rot \bm{\Psi}) + (\dive \hat{\bm{w}}, \dive \bm{\Psi})  
= (-\hat {\bm {h}}, \rot \bm{\Psi}) 
= (-\rot \hat {\bm{h}}, \bm{\Psi}) + \langle \hat {\bm{h}}, \bm{\Psi}\times \bm{\nu}\rangle_{\pt\Om} 
= 0
\end{eqnarray*} 
for all $\bm{\Psi} \in \dot V^{r^\prime}(\Om)$. Hence, $\hat {\bm{w}} \in \mbox{\rm Ker}(T_r)$. Finally, by Proposition \ref{pr:4.5} (ii)  
$\hat {\bm{w}} \in \dot V_{\tiny{\mbox{\rm har}}}(\Omega)$ and $\rot \hat {\bm{w}} =\bm{0}$. i.e., 
\begin{equation}\label{eqn:6.10}
\rot \bm{w} = \rot \tilde {\bm{w}}.
\end{equation}
Assertions (\ref{eqn:6.8}), (\ref{eqn:6.9}) and (\ref{eqn:6.10}) imply $\hat {\bm{h}} =\bm{0}$, i.e., 
\begin{equation}\label{eqn:6.11}
\bm{h} = \tilde {\bm{h}}
\end{equation} 
and the desired unique decomposition (\ref{eqn:2.9}) is a consequence of  (\ref{eqn:6.9}), (\ref{eqn:6.10}) and  (\ref{eqn:6.11}).  
This completes the proof of Theorem \ref{thm:2.2}. \qed
\subsection{Proof of Theorem \ref{thm:2.3}}\mbox{} 
\par
Let $\bm{u} \in L^r(\Omega)$.  For $1< r \le 3/2$ take $p \in \dot H^{1, r}_0(\Omega)$ given in Lemma \ref{lem:5.3}. 
For $3/2 < r < 3$ and for $3 \le r < \infty$ take $p \in \widehat H^{1, r}_0(\Omega)$ and 
$p \in \widetilde H^{1, r}_0(\Omega)$ given in Lemma \ref{eqn:5.2} (i) and (ii), respectively.  
Define $\bm{h}$ as in (\ref{eqn:6.5}). Our first aim is to show that $h \in V^r_{\tiny{\mbox{\rm har}}}(\Omega)$.  
Obviously, $\bm{h} \in L^r(\Omega)$. In the same way as done in (\ref{eqn:6.6}) and (\ref{eqn:6.7}),  we verify that $\dive \bm{h} =0$ and $\rot \bm{h} =\bm{0}$ in the sense of distributions in $\Omega$.  
Hence, we may show that $\tau_{\bm{\nu}}\bm{h} =0$.  
Since $\gamma_0p=0$, assertions (\ref{eqn:6.1}) and (\ref{eqn:6.2}) imply  
$$
\langle\tau_{\bm{\nu}}(\nabla p), \gamma_0\bm{\psi}\rangle_{\pt\Om} = (\nabla p, \rot\bm{\psi}) 
= \langle\gamma_0p, \gamma_{\bm{\nu}}(\rot \bm{\psi})\rangle_{\pt\Om} =0
\quad
\mbox{for all $\bm{\psi} \in H^{1, r^\prime}(\Omega)$},   
$$ 
which yields 
\begin{equation}\label{eqn:6.12}
\tau_{\bm{\nu}} (\nabla p)= 0.  
\end{equation} 
For every $\bm{\psi} \in H^{1, r^\prime}(\Omega)$ we consider the weak Neumann problem in $\dot H^{1, r^\prime}(\Omega)$ for  the Poisson equation 
$$
\left\{
\begin{array}{ll}
& \Delta q = \dive \bm{\psi} \quad \mbox{in $\Omega$}, \\
& \dfrac{\pt q}{\pt\bm{\nu}}=\bm{\psi}\cdot\bm{\nu} \quad\mbox{on $\pt\Om$}.
\end{array}
\right.
$$
It is known that there exits a unique $q \in \dot H^{1, r^\prime}(\Omega)$ satisfying  $(\nabla q, \nabla \varphi) = (\bm{\psi}, \nabla \varphi)$ for all $\varphi \in \dot H^{1, r}(\Omega)$. 
Moreover, similarly as in  (\ref{eqn:4.49}) we deuce  $D^2q \in L^{r^\prime}(\Omega)$.  Setting $\bm{\Phi} \equiv \bm{\psi} - \nabla q$ and we see that  $\bm{\Phi} \in L^{r^\prime}(\Omega)$ 
with $\dive \bm{\Phi}=0$, $\rot \bm{\Phi} \in L^{r^\prime}(\Omega)$ and $\bm{\Phi}\cdot\bm{\nu}|_{\pt\Omega}=0$.  Hence, $\bm{\Phi} \in \dot X^{r^\prime}(\Omega)$.       
Since $\rot\bm{\psi} = \rot \bm{\Phi}$ and since $\rot (\bm{u}- \rot \bm{w}) = \rot(\bm{h} + \nabla p)=\bm{0}$, (\ref{eqn:6.2}) and (\ref{eqn:4.1}) yield 
\begin{eqnarray*}
\langle\tau_{\bm{\nu}}(\bm{u} - \rot \bm{w}), \gamma_0\bm{\psi}\rangle_{\pt\Om}
&=&(\bm{u} - \rot \bm{w}, \rot \bm{\psi}) - (\rot (\bm{u}- \rot \bm{w}), \bm{\psi}) \\
&=&(\bm{u} - \rot \bm{w}, \rot \bm{\psi}) \\
&=&(\bm{u} - \rot \bm{w},  \rot \bm{\Phi}) \\
&=&0.
\end{eqnarray*}
Since $\bm{\psi} \in H^{1, r^\prime}(\Omega)$ is  arbitrary, we obtain   
\begin{equation}\label{eqn:6.13}
\tau_{\bm{\nu}}(\bm{u} - \rot \bm{w})=0.  
\end{equation}
By (\ref{eqn:6.12}) and (\ref{eqn:6.13}), $\tau_{\bm{\nu}} \bm{h}=0$ and  $\bm{h} \in V^r_{\tiny{\mbox{\rm har}}}(\Omega)$ with the decomposition (\ref{eqn:2.10}). 
The estimate (\ref{eqn:2.11}) is then a consequence of 
(\ref{eqn:6.5}), (\ref{eqn:4.2}) and (\ref{eqn:5.4}).
\par
It remains to prove the uniqueness of the decomposition  (\ref{eqn:2.10}). \\   
(i) Let $1 < r \le 3/2$ and suppose that $\tilde {\bm{h}} \in V^r_{\tiny{\mbox{\rm har}}}(\Omega)$, $\tilde {\bm{w}} \in \dot X^r_{\sg}(\Omega)$ and 
$\tilde p \in \dot H^{1, r}_0(\Omega)$ satisfy (\ref{eqn:2.12}).  
Setting  $\hat {\bm{h}} \equiv \bm{h} - \tilde {\bm{h}}$, 
$\hat {\bm{w}} \equiv \bm{w} - \tilde {\bm{w}}$ and $\hat p = p - \tilde p$, 
we obtain (\ref{eqn:6.8}). Moreover, For every $\varphi \in C^\infty_0(\Omega)$ we obtain 
\begin{equation}\label{eqn:6.14}
(\nabla \hat p, \nabla \varphi) = -(\hat {\bm{h}}, \nabla \varphi)-(\rot \hat {\bm{w}}, \nabla\varphi) 
= (\dive \hat {\bm{h}}, \varphi) - (\hat {\bm{w}},  \rot(\nabla\varphi)) =0.
\end{equation}  
Thus $\hat p$ is harmonic in $\Omega$.  Since $\hat p|_{\pt\Omega}=0$ and since $\nabla \hat p \in L^{r}(\Omega)$, it follows from \cite[Theorem A]{KoSo} that 
$\hat p \equiv 0$ in $\Omega$ and hence (\ref{eqn:6.9}).  
\par
Furthermore $\rot \hat {\bm{w}} = -\hat {\bm{h}}$. Since $\hat {\bm{h}}\times {\bm{\nu}}|_{\pt\Om} = \bm{0}$, (\ref{eqn:6.2}) yields 
\begin{eqnarray}
a_X(\hat {\bm{w}}, \bm{\Phi}) &= &(\rot \hat {\bm{w}}, \rot \bm{\Phi}) 
+ (\dive \hat {\bm{w}}, \dive \bm{\Phi}) \nonumber \\ 
&=& (-\hat {\bm{h}}, \rot \bm{\Phi}) \nonumber \\
&=& (-\rot \hat {\bm{h}}, \bm{\Phi}) - \langle \hat {\bm{h}}\times {\nu}, \bm{\Phi} \rangle_{\pt\Om} 
\label{eqn:6.15}\\
&=& 0 \nonumber
\end{eqnarray} 
for all $\bm{\Phi} \in \dot X^{r^\prime}(\Om)$. Thus  $\hat {\bm{w}} \in \mbox{Ker}(S_r)$ and it follows from Proposition \ref{pr:4.5} that 
$\hat {\bm{w}} \in \dot X_{\tiny{\mbox{\rm har}}}(\Omega)$ 
yielding  $\rot \hat {\bm{w}} =0$. Hence, we obtain (\ref{eqn:6.10}). 
The desired unique decomposition  (\ref{eqn:2.13}) is a consequence of  (\ref{eqn:6.9}), (\ref{eqn:6.10}) and  (\ref{eqn:6.11}).
\par
(ii) Let $3/2 < r < 3$. Suppose that 
$\tilde {\bm{h}} \in V^r_{\tiny{\mbox{har}}}(\Omega)$, $\tilde {\bm{w}} \in \dot X^r_{\sg}(\Omega)$ and 
$\tilde p \in \hat H^{1, r}_0(\Omega)$ satisfy (\ref{eqn:2.12}). Note that  the identity (\ref{eqn:6.14}) 
remains true for $\hat p \in H^{1, r}_0(\Omega)$. 
Since $3/2 < r < 3$, $C^\infty_0(\Omega)$ is dense in $\widehat H^{1, r^\prime}_0(\Omega)$ 
it follows from  (\ref{eqn:6.14}) that  
\begin{equation}\label{eqn:6.16}
(\nabla \hat p, \nabla \varphi)=0
\quad
\mbox{for all $\varphi \in \widehat H^{1, r^\prime}_0(\Omega)$}, 
\end{equation}
yielding  $\hat p \in \mbox{Ker}(-\Delta_r)$.  Proposition \ref{pr:3.2} (ii) yields  $\hat p =0$ 
and hence (\ref{eqn:6.9}). 
Since (\ref{eqn:6.15}) remains true even for $3/2 < r <3$, we arrive at  (\ref{eqn:6.10}) and hence  (\ref{eqn:6.11}) and (\ref{eqn:2.13}). 
\par
(iii) Let $3 \le r < \infty$ and suppose that $\tilde {\bm{h}} \in V^r_{\tiny{\mbox{har}}}(\Omega)$, $\tilde {\bm{w}} \in \dot X^r_{\sg}(\Omega)$ and 
$\tilde p \in \dot H^{1, r}_0(\Omega)$ satisfy (\ref{eqn:2.12}).  
Since $C^\infty_0(\Omega)$ is dense in $\widehat H^{1, r^\prime}_0(\Omega)$, (\ref{eqn:6.16}) remains true and hence 
$\hat p \in \mbox{Ker}(-\Delta_r)$. By Proposition \ref{pr:3.2} (iv)  
\begin{equation}\label{eqn:6.17}
\hat p = \lambda q_0
\quad
\mbox{for some $\lambda \in \re$}.  
\end{equation}
The assertions  (\ref{eqn:6.8}) and (\ref{eqn:6.17}) imply  $\rot \hat {\bm{w}} = - \hat {\bm{h}} -\lambda \nabla q_0$. 
Since $\hat {\bm{h}}\times \bm{\nu}|_{\pt\Om} = \nabla q_0 \times {\bm{\nu}}|_{\pt\Om} = \bm{0}$, (\ref{eqn:6.2}) implies that 
\begin{eqnarray*}
a_X(\hat {\bm{w}}, \bm{\Phi}) &=& (\rot \hat {\bm{w}}, \rot \bm{\Phi}) + (\dive \hat {\bm{w}}, \dive \bm{\Phi})  \\ 
&=& (-\hat {\bm{h}} - \lambda \nabla q_0, \rot \bm{\Phi})   \\
&=& (-\rot( \hat {\bm{h}} + \lambda \nabla q_0), \bm{\Phi})
 - \langle (\hat {\bm{h}} + \lambda \nabla q_0)\times {\bm{\nu}}, \bm{\Phi} \rangle_{\pt\Om} \\
&=& 0
\end{eqnarray*} 
for all $\bm{\Phi} \in \dot X^{r^\prime}(\Om)$. Hence, $\hat {\bm{w}} \in \mbox{Ker}(S_r)=\dot X_{\tiny{\mbox{\rm har}}}(\Omega)$ and thus (\ref{eqn:6.10}).  
Finally, by (\ref{eqn:6.8}), (\ref{eqn:6.10}) and (\ref{eqn:6.17})  
\begin{equation}\label{eqn:6.18}
\hat {\bm{h}} = \lambda \nabla q_0  
\end{equation}
with the same $\lambda$ as above. The desired uniqueness assertion in (\ref{eqn:2.14}) is then a consequence of  (\ref{eqn:6.10}), (\ref{eqn:6.17}) and (\ref{eqn:6.18}). 
The proof of Theorem \ref{thm:2.3} is complete. \qed

\bigskip
Finally, some comments concerning the  remark after Theorem \ref{thm:2.3} are in order.  In fact, it is impossible to replace $p \in \dot H^{1, r}_0(\Omega)$ by  
$p \in \widehat H^{1, r}_0(\Omega)$ provided $1 < r \le 3/2$. 
More precisely, let $\eta \in C^\infty_0(\Omega)$ be the cut-off function introduced  in Lemma \ref{lem:3.9}.  
Assume that there exist $\bm{h} \in X_{\tiny{\mbox{\rm har}}}(\Omega)$, $\bm{w} \in \dot X^r_{\sg}(\Omega)$ and $p\in \widehat H^{1, r}_0(\Omega)$  such that 
$$
\nabla \eta = \bm{h} + \rot \bm{w} + \nabla p.
$$
Then, for every $\bm{\Phi}\in \dot X^{r^\prime}(\Omega)$ we have 
\begin{eqnarray*}
a_X(\bm{w}, \bm{\Phi}) &=& (\rot \bm{w}, \rot \bm{\Phi}) + (\dive \bm{w}, \dive \bm{\Phi})   \\
&=& (-\bm{h} + \nabla(\eta - p), \rot \bm{\Phi})   \\
&=& (\rot( -\bm{h} + \nabla(\eta - p)), \bm{\Phi}) + 
  \langle (-\bm{h} + \nabla (\eta -p))\times \bm{\nu}, \bm{\Phi} \rangle_{\pt\Om} \\
&=& 0.
\end{eqnarray*} 
By Proposition \ref{pr:4.5}, $\rot \bm {w} =0$ and hence $\nabla \eta = \bm{h} + \nabla p$.  
Set $g(x) = \eta(x) - p(x)$ and we notice  that $g$ is harmonic in $\Omega$. Setting $\tilde g(x) = g(x) -1$, we verify  that 
$$
\left\{
\begin{array}{ll}
& \Delta \tilde g=0 \quad\mbox{in $\Omega$}, \\
& \tilde g =0 \quad\mbox{on $\pt\Omega$}.  
\end{array}
\right.
$$ 
Since $\nabla \tilde g \in L^r(\Omega)$ for $1< r \le 3/2$, it follows from \cite[Theorem A]{KoSo} that $\tilde g \equiv 0$ on $\Omega$, and hence $g(x) \equiv 1$ for all $x\in \Om$.  
Finally,  $\eta (x)= 0$ for $|x| \ge R+2$ and thus  $p(x) = -1$ for all $|x| \ge R+2$. This  contradicts the fact that $p \in L^{r_\ast}(\Omega)$, where $r_\ast$ is defined by 
(\ref{eqn:3.1}). 
\par
\bigskip
{\bf Acknowledgement.} The authors would like to express their thanks to Prof.Yoshihiro Shibata for valuable 
comments, in particular concerning Lemma \ref{lem:5.3}.

\end{document}